\documentclass[a4paper,leqno]{article}

\usepackage{ifthen}
\usepackage[latin1]{inputenc}
\usepackage[T1]{fontenc}
\usepackage{xspace}
\usepackage{amssymb}
\usepackage{amsmath}
\usepackage{amsthm,amsxtra}
\usepackage[all]{xy}
\usepackage[english,francais]{babel}

\theoremstyle{definition}
\newtheorem{ssect}[subsection]{}
\newtheorem{sssect}[subsubsection]{}
\newtheorem{defn}[subsection]{Définition}

\newtheorem*{remqe}{Remarque}
\newtheorem{remq}[subsection]{Remarque}
\newtheorem{exem}[subsection]{Exemple}

\newtheorem*{escho}{Scholie}

\theoremstyle{plain}
\newtheorem{lemm}[subsection]{Lemme}

\newtheorem{prop}[subsection]{Proposition}

\newtheorem{theo}[subsection]{Théorème}

\newtheorem{coro}[subsection]{Corollaire}

\newcommand{\beq}{\begin{equation}}
\newcommand{\eeq}{\end{equation}}

\DeclareMathSymbol{\oio}{\mathopen}{operators}{"5D} 
\DeclareMathSymbol{\fio}{\mathclose}{operators}{"5B} 

\newcommand{\al}{\alpha}

\newcommand{\eps}{\epsilon}
\newcommand{\Lbd}{\Lambda}
\newcommand{\lbd}{\lambda}

\newcommand{\cF}{\mathcal{F}}
\newcommand{\cG}{\mathcal{G}}
\newcommand{\fH}{\mathcal{H}}
\newcommand{\fL}{\mathcal{L}}
\newcommand{\cO}{\mathcal{O}}

\newcommand{\cU}{\mathcal{U}}
\newcommand{\cV}{\mathcal{V}}
\newcommand{\fm}{\mathfrak{m}}
\newcommand{\bA}{\mathbb{A}}
\newcommand{\F}{\mathbb{F}}
\newcommand{\N}{\mathbb{N}}
\newcommand{\bP}{\mathbb{P}}
\newcommand{\Q}{\mathbb{Q}}

\newcommand{\Z}{\mathbb{Z}}

\newcommand{\lbar}{\overline}

\newcommand{\ot}{\otimes}
\newcommand{\hra}{\hookrightarrow}
\newcommand{\vide}{\emptyset}

\newcommand{\Qlb}{\overline{\Q_l}}
\newcommand{\Dbc}{D^b_c}
\newcommand{\cHom}{\fH\mathit{om}}
\newcommand{\simto}[1][]{\sto[#1]{\sim}}

\newcommand{\ensdr}[2]{\left\{ #1 \,\left|\, #2 \right.\right\}}

\newcommand{\sto}[2][]{\xrightarrow[#1]{#2}} 

\newcommand{\Aut}{\mathrm{Aut}}

\newcommand{\et}{\mathrm{\acute{e}t}}
\newcommand{\End}{\mathrm{End}}
\newcommand{\Fr}{\mathrm{Fr}}

\newcommand{\Gal}{\mathrm{Gal}}
\newcommand{\GL}{\mathrm{GL}}
\newcommand{\Hom}{\mathrm{Hom}}
\newcommand{\Id}{\mathrm{Id}}
\newcommand{\Mod}{\mathrm{Mod}}
\newcommand{\op}{\mathrm{op}}
\newcommand{\red}{\mathrm{red}}

\newcommand{\Rep}{\mathrm{Rep}}

\newcommand{\Tr}{\mathrm{Tr}}

\DeclareMathOperator{\Spec}{Spec}
\DeclareMathOperator{\Ker}{Ker}
\DeclareMathOperator{\Img}{Im}

\newcommand{\algt}{algébriquement\xspace}

\newcommand{\carrecart}{carré cartésien\xspace}

\newcommand{\can}{canonique\xspace}
\newcommand{\cara}{caractéristique\xspace}
\newcommand{\cats}{catégories\xspace}
\newcommand{\cat}{catégorie\xspace}
\newcommand{\Cf}{\emph{Cf.\@} }
\newcommand{\cf}{\emph{cf.\@} }

\newcommand{\comm}{commutatif\xspace}

\newcommand{\compte}{compatibilité\xspace}

\newcommand{\cn}{à croisements normaux\xspace}
\newcommand{\dcn}{diviseur \cn}

\newcommand{\df}{de dimension finie\xspace}
\newcommand{\diag}{diagramme\xspace}
\newcommand{\diagcomm}{diagramme commutatif\xspace}
\newcommand{\dist}{distingué\xspace}

\newcommand{\ens}{ensemble\xspace}

\newcommand{\ev}{espace vectoriel\xspace}
\newcommand{\faisc}{faisceau\xspace}
\newcommand{\faisx}{faisceaux\xspace}

\newcommand{\gent}{génériquement\xspace}

\newcommand{\geoq}{géométrique\xspace}
\newcommand{\geoqt}{géométriquement\xspace}

\newcommand{\hyp}{hypothèse\xspace}
\newcommand{\ie}{\emph{i.~e.}, }
\newcommand{\irrs}{irréductibles\xspace}
\newcommand{\irr}{irréductible\xspace}
\newcommand{\iso}{isomorphisme\xspace}

\newcommand{\lin}{linéaire\xspace}

\newcommand{\mrs}{modérément ramifiés\xspace}
\newcommand{\mr}{modérément ramifié\xspace}

\newcommand{\nre}{non ramifiée\xspace}
\newcommand{\noeths}{noethériens\xspace}
\newcommand{\noeth}{noethérien\xspace}
\newcommand{\Ops}{On peut supposer\xspace}
\newcommand{\ops}{on peut supposer\xspace}

\newcommand{\qs}{quasi-\sep}
\newcommand{\rec}{récurrence\xspace}
\newcommand{\reg}{régulier\xspace}

\newcommand{\repns}{représentations\xspace}

\newcommand{\resp}{resp.\@\xspace} 

\newcommand{\schs}{schémas\xspace}
\newcommand{\sch}{schéma\xspace}

\newcommand{\sep}{séparé\xspace}
\newcommand{\ssi}{si et seulement si\xspace}

\newcommand{\teq}{telle que\xspace}

\newcommand{\tf}{de type fini\xspace}

\newcommand{\tq}{tel que\xspace}
\newcommand{\tranvt}{transitivement\xspace}
\newcommand{\trdist}{triangle distingué\xspace}

\newcommand{\tsq}{tels que\xspace}
\newcommand{\varab}{variété abélienne\xspace}
\newcommand{\ver}{vérifi}

\newcommand{\ibid}[1][]{[\emph{ibid.}\ifthenelse{\equal{#1}{}}{}{, }#1]\xspace} 

\makeatletter

\newenvironment{cond}{\begin{list}{}{%
\setlength{\leftmargin}{0mm}%
\setlength{\rightmargin}{0mm}%
\setlength{\itemindent}{\parindent}%
\setlength{\listparindent}{\parindent}%
\setlength{\parsep}{\parskip}%
}\item[]\ignorespaces}{\unskip\end{list}}

\newenvironment{numcond}[1][subsubsection]{\stepcounter{#1}%
\begin{list}{}{%
\setlength{\leftmargin}{0mm}%
\setlength{\rightmargin}{0mm}%
\setlength{\labelwidth}{0mm}%
\setlength{\labelsep}{\parindent}%
\setlength{\itemindent}{\parindent}%
\setlength{\listparindent}{\parindent}%
\setlength{\parsep}{\parskip}%
}\def\@currentlabel{\csname the#1\endcsname}%
\item[\textbf{\@currentlabel}]\ignorespaces}{\unskip\end{list}}

\makeatother

\makeatletter
\newcommand{\sigle}[1]{\protect\@sigle{#1}}
\newcommand{\@sigle}[1]{\gdef\@sigleparam{#1}\aftergroup\@absorbtwo}
\newcommand{\@absorbtwo}[2]{\ifcat a#1#1#2\else \@sigleparam \fi}
\makeatother

\newcommand{\Zch}{\hat{\Z}}
\newcommand{\gb}{{\bar{g}}}
\newcommand{\kb}{{\bar{k}}}
\newcommand{\Kb}{{\lbar{K}}}
\newcommand{\sbar}{{\bar{s}}}
\newcommand{\wb}{{\lbar{w}}}
\newcommand{\xb}{{\bar{x}}}
\newcommand{\yb}{{\bar{y}}}
\newcommand{\zb}{{\bar{z}}}
\newcommand{\etab}{{\bar{\eta}}}
\newcommand{\Qllb}{{\lbar{\Q_{l_\lbd}}}}
\newcommand{\ip}{{(i)}}

\newcommand{\Dbm}{D^b_m}

\newcommand{\ldI}{{\lbd \in I}}

\makeatletter\AtBeginDocument{\@ifpackageloaded{fourier}%
{\DeclareMathAlphabet{\mathbfit}{T1}{futs}{b}{it}}%
{\DeclareMathAlphabet{\mathbfit}{OT1}{cmr}{bx}{it}}}\makeatother

\newcommand{\Eq}{\mathbfit{Eq}}
\newcommand{\Grfini}{\mathbfit{Gr.fini}}
\newcommand{\bC}{\mathbfit{C}}
\newcommand{\bD}{\mathbfit{D}}
\newcommand{\bF}{\mathbfit{F}}

\newcommand{\adm}{\mathrm{adm}}
\newcommand{\pr}{\mathrm{pr}}
\newcommand{\ord}{\mathrm{ord}}
\newcommand{\Fix}{\mathrm{Fix}}
\newcommand{\Swan}{\mathrm{Swan}}
\newcommand{\Ob}{\mathrm{Ob}}
\newcommand{\pt}{\mathrm{pt}}
\newcommand{\naif}{\text{naïf}}

\newcommand{\altgals}{altérations galoisiennes\xspace}
\newcommand{\eqrt}{équivariant\xspace}

\newcommand{\eqrte}{équivariante\xspace}

\newcommand{\Ecomp}{$E$-compatible\xspace}
\newcommand{\Ecomps}{$E$-compatibles\xspace}
\newcommand{\Ecompte}{$E$-\compte}

\IfFileExists{alafrench.tex}{\usepackage{aeguill}  \DeclareMathSymbol{A}{\mathalpha}{operators}{"41}
  \DeclareMathSymbol{B}{\mathalpha}{operators}{"42}
  \DeclareMathSymbol{C}{\mathalpha}{operators}{"43}
  \DeclareMathSymbol{D}{\mathalpha}{operators}{"44}
  \DeclareMathSymbol{E}{\mathalpha}{operators}{"45}
  \DeclareMathSymbol{F}{\mathalpha}{operators}{"46}
  \DeclareMathSymbol{G}{\mathalpha}{operators}{"47}
  \DeclareMathSymbol{H}{\mathalpha}{operators}{"48}
  \DeclareMathSymbol{I}{\mathalpha}{operators}{"49}
  \DeclareMathSymbol{J}{\mathalpha}{operators}{"4A}
  \DeclareMathSymbol{K}{\mathalpha}{operators}{"4B}
  \DeclareMathSymbol{L}{\mathalpha}{operators}{"4C}
  \DeclareMathSymbol{M}{\mathalpha}{operators}{"4D}
  \DeclareMathSymbol{N}{\mathalpha}{operators}{"4E}
  \DeclareMathSymbol{O}{\mathalpha}{operators}{"4F}
  \DeclareMathSymbol{P}{\mathalpha}{operators}{"50}
  \DeclareMathSymbol{Q}{\mathalpha}{operators}{"51}
  \DeclareMathSymbol{R}{\mathalpha}{operators}{"52}
  \DeclareMathSymbol{S}{\mathalpha}{operators}{"53}
  \DeclareMathSymbol{T}{\mathalpha}{operators}{"54}
  \DeclareMathSymbol{U}{\mathalpha}{operators}{"55}
  \DeclareMathSymbol{V}{\mathalpha}{operators}{"56}
  \DeclareMathSymbol{W}{\mathalpha}{operators}{"57}
  \DeclareMathSymbol{X}{\mathalpha}{operators}{"58}
  \DeclareMathSymbol{Y}{\mathalpha}{operators}{"59}
  \DeclareMathSymbol{Z}{\mathalpha}{operators}{"5A}

\DeclareSymbolFont{psym}{U}{psy}{m}{n}

\DeclareMathSymbol{\alpha}{\mathord}{psym}{97}
\DeclareMathSymbol{\beta}{\mathord}{psym}{98}
\DeclareMathSymbol{\gamma}{\mathord}{psym}{103}
\DeclareMathSymbol{\delta}{\mathord}{psym}{100}
\DeclareMathSymbol{\epsilon}{\mathord}{psym}{101}
\DeclareMathSymbol{\varepsilon}{\mathord}{psym}{101}
\DeclareMathSymbol{\zeta}{\mathord}{psym}{122}
\DeclareMathSymbol{\eta}{\mathord}{psym}{104}
\DeclareMathSymbol{\theta}{\mathord}{psym}{113}
\DeclareMathSymbol{\vartheta}{\mathord}{psym}{74}
\DeclareMathSymbol{\iota}{\mathord}{psym}{105}
\DeclareMathSymbol{\kappa}{\mathord}{psym}{107}
\DeclareMathSymbol{\lambda}{\mathord}{psym}{108}
\DeclareMathSymbol{\mu}{\mathord}{psym}{109}
\DeclareMathSymbol{\nu}{\mathord}{psym}{110}
\DeclareMathSymbol{\xi}{\mathord}{psym}{120}
\DeclareMathSymbol{\omicron}{\mathord}{psym}{111}
\DeclareMathSymbol{\pi}{\mathord}{psym}{112}
\DeclareMathSymbol{\varpi}{\mathord}{psym}{118}
\DeclareMathSymbol{\rho}{\mathord}{psym}{114}
\DeclareMathSymbol{\sigma}{\mathord}{psym}{115}
\DeclareMathSymbol{\varsigma}{\mathord}{psym}{86}
\DeclareMathSymbol{\tau}{\mathord}{psym}{116}
\DeclareMathSymbol{\upsilon}{\mathord}{psym}{117}
\DeclareMathSymbol{\phi}{\mathord}{psym}{102}
\DeclareMathSymbol{\varphi}{\mathord}{psym}{106}
\DeclareMathSymbol{\chi}{\mathord}{psym}{99}
\DeclareMathSymbol{\psi}{\mathord}{psym}{121}
\DeclareMathSymbol{\omega}{\mathord}{psym}{119}
}%
{\usepackage[upright]{fourier}\newcommand{\omicron}{\mathrm{o}}}
\usepackage[unicode]{hyperref}
\numberwithin{equation}{subsection}

\renewcommand{\sharp}{\#}
\newcommand{\tm}[1]{\textup{#1}}

\begin{document}
\title{Sur l'indépendance de $l$ en cohomologie $l$-adique sur les corps locaux}
\author{Weizhe Zheng}
\date{}

\maketitle%
{\renewcommand{\thefootnote}{}%
\footnotetext{Date : 5 mars 2009.\label{version=1e}}}

\begin{abstract}
  Gabber a déduit son théorème d'indépendance de~$l$ de la cohomologie d'intersection
  d'un résultat général de stabilité sur les corps finis.
  Dans cet article, nous démontrons un analogue sur les corps locaux de ce résultat général.
  Plus précisément, nous introduisons une notion d'indépendance
  de~$l$ pour les systèmes
  de complexes de faisceaux $l$-adiques sur les schémas de type
  fini sur un corps local équivariants sous des groupes finis et nous établissons sa stabilité par les six
  opérations de Grothendieck et le foncteur des cycles proches.
  Notre méthode permet d'obtenir une nouvelle démonstration du
  théorème de Gabber. Nous donnons aussi une généralisation aux
  champs algébriques.
\end{abstract}

\section*{Introduction}

Dans les années 1980, Gabber a prouvé l'indépendance de~$l$ de la
cohomologie d'intersection \cite[Th.~1]{Fujiwara} : pour $X$ un
schéma propre sur un corps fini $k=\F_{p^f}$, équidimensionnel et
$i\in \Z$, $\det (1-tF_0^f , IH^i(X_{\kb} ,\Q_l))$ est un polynôme
dans $\Z[t]$, indépendant de $l\neq p$. Ici $\kb$ est une clôture
algébrique de $k$, $F_0 \in\Aut(\kb)$ est le Frobenius
géométrique absolu qui envoie $a$ sur $a^{1/p}$. Il déduit ce
résultat d'un théorème général d'indépendance de~$l$
\ibid[Th.~2] : pour $E$ une extension de~$\Q$, il définit une
notion de $E$-compatibilité pour les systèmes de complexes
$l$-adiques ($l\neq p$) sur les schémas séparés de type fini
sur~$k$ et établit la stabilité de cette $E$-compatibilité par
les six opérations de Grothendieck.

L'objet de cet article est d'étudier l'indépendance de~$l$ sur un
corps local~$K$, par lequel on entend le corps des fractions d'un
anneau de valuation discrète hensélien excellent de corps
résiduel fini. De façon précise, on définit une notion de
$E$-compatibilité pour les systèmes de complexes $l$-adiques sur
les schémas de type fini sur~$K$. On prouve la stabilité de cette
$E$-compatibilité par les opérations usuelles, notamment les six
opérations et le foncteur des cycles proches $R\Psi$.

La $E$-compatibilité est sensible aux morphismes finis. Dans
l'étude de la $E$-compatibilité, il est naturel d'utiliser les
résultats de de Jong sur les altérations équivariantes
\cite{deJong2}, ce qui amène à généraliser le problème au cas
équivariant sous des actions de groupes finis.

Au \S~\ref{s1}, après avoir fixé les notations, on définit la
$E$-compatibilité pour les systèmes de complexes équivariants
au-dessus d'un corps fini ou local et énonce la stabilité par les
six opérations. La démonstration dans le cas d'un corps fini est
donnée au \S~\ref{s2}. En vue du théorème de Gabber de la
stabilité de la $E$-compatibilité par les six opérations, il
suffit d'appliquer une technique de Deligne-Lusztig qui permet de se
débarrasser des actions de groupes finis par descente galoisienne.
Au \S~\ref{s3} on se réduit au cas des courbes. L'ingrédient
essentiel est un raffinement \cite[4.4]{Vidal}, dû à Gabber, de
résultats de de Jong, grâce auquel on se ramène au cas de $Rj_*
L_\lbd$ pour l'inclusion $j: U\to X$ du complémentaire d'un
diviseur à croisements normaux $G$-strict $D$ dans un schéma
régulier~$X$ et des complexes $G$-équivariants $L_\lbd$ sur~$U$ à
faisceaux de cohomologie lisses, modérément ramifiés le long
de~$D$. Dans le cas d'un corps fini, le cas des courbes étant
déjà connu depuis Deligne \cite[9.8]{constantes}, on obtient une
démonstration indépendante du théorème de Gabber. La fin de la
démonstration des énoncés du \S~\ref{s1} dans le cas d'un corps
local est donnée au \S~\ref{s4}. On les déduit, à l'aide des
résultats de de Jong, d'un cas particulier de la stabilité par
$R\Psi$, que l'on traite en utilisant à nouveau la technique de
Deligne-Lusztig. Le cas général de la stabilité par $R\Psi$
découle de ce cas particulier et encore une fois des résultats de
de Jong.

Les schémas équivariants sont mieux compris dans le contexte de
leurs champs quotients associés. Ce point est élucidé au
\S~\ref{s5}. On y donne aussi des généralisations des résultats
d'indépendance de~$l$ aux champs algébriques.

Ce travail fait partie de ma thèse. Je remercie vivement mon
directeur de thèse, L.~Illusie, qui m'a donné d'innombrables
conseils lors de la préparation de l'article. Je remercie aussi
N.~Katz pour sa suggestion d'utiliser la technique de
Deligne-Lusztig. Je remercie enfin S.~Morel et F.~Orgogozo pour
leurs longues listes de commentaires.

\section{Indépendance de $l$ et six opérations}\label{s1}

\begin{ssect}\label{ss.cat}
Pour une catégorie~$\bC$, on note $\Eq(\bC)$ la catégorie des
objets de~$\bC$ munis d'une action d'un groupe fini à droite. Les
objets sont les triplets $(X,G,a)$, où $X$~est un objet de~$\bC$,
$G$~est un groupe fini, $a$~est une action de $G$ sur $X$ à droite
(\ie $a: G^\op \to \Aut_{\bC}(X)$ est un homomorphisme de groupes).
On enlève $a$ de la notation lorsqu'il n'y a pas de confusion à
craindre. Un morphisme $(X_1, G_1) \to (X_2, G_2)$ est un couple
$(f,\al)$, où $\al :G_1\to G_2$ est un homomorphisme de groupes,
$f: X_1 \to X_2$ est un morphisme $G_1$-\eqrt de $\bC$, l'action de
$G_1$ sur $X_2$ étant induite par $\al$. La condition de
$G_1$-équivariance signifie que pour tout $g\in G_1$, le \diag
\[
\xymatrix{X_1\ar[d]^{f} \ar[r]^{a_g} & X_1 \ar[d]^{f}\\
X_2\ar[r]^{a_{\al(g)}} & X_2}
\]
est \comm. Pour $(X,G,a)$ un objet de $\Eq(\bC)$ et $g\in G$, $T_g =
(a_g , h\mapsto g^{-1}h g)$ est un automorphisme de $(X,G, a)$. Ceci
définit une action de $G$ sur $(X,G,a)$ à droite. Si les produits
fibrés sont représentables dans $\bC$, il en est de même dans
$\Eq(\bC)$ :
\[(X_1,G_1)\times_{(X,G)} (X_2,G_2) = (X_1\times_X X_2,
G_1\times_G G_2).\]

On note $\Grfini$ la catégorie des groupes finis. Le foncteur
projection $\Eq(\bC) \to \Grfini$ qui envoie $(X,G)$ sur $G$ est
fibrant. On note $\Eq(\bC)_G$ sa fibre en~$G$. Pour $\al : G\to H$
un homomorphisme de groupes finis,
\[\al^* : \Eq(\bC)_H \to \Eq(\bC)_G\]
est donné par $\al^* (Y,H,a) = (Y,G,a\circ \al^\op)$. Pour $\al : G
\to H$ injectif, l'adjoint à gauche
\[\al_* : \Eq(\bC)_G \to
\Eq(\bC)_H
\]
de $\al^*$ existe dès que les sommes disjointes finies sont
représentables dans $\bC$. Pour l'expliciter, choisissons un
système $S$ de représentants de $\al(G)\backslash H$. Pour $(X,G)\in
\Eq(\bC)$, on a
\beq\label{eq.aleb}
\al_* (X,G) \simeq (\coprod_{s\in S} X_s, H),
\eeq
où $X_s = X$ et l'action de $h\in H$ est donnée par $a_g: X_s \to
X_t$, où $g\in G$ et $t\in S$ vérifient $sh=\al(g) t$.
Rappelons qu'une flèche $(f,\al): (X,G) \to (Y,H)$ dans $\Eq(\bC)$
est \emph{cocartésienne} \cite[VI~10]{SGA1} si pour tout objet
$(Y',H)$ de $\Eq(\bC)_H$, l'application
\begin{align*}
\Hom_{H}((Y,H),(Y',H)) & \to \Hom_\al((X,G),(Y',H))\\
u&\mapsto u \circ (f,\al)
\end{align*}
est bijective.

\stepcounter{subsubsection}
\begin{sssect}\label{sss.cocart}
Pour $(Y,H)\in \Eq(\bC)$, si $Y$ est la somme disjointe d'une
famille finie d'objets $(X_j)_{j\in J}$ de~$\bC$, et s'il existe une
action \emph{transitive} de $H$ sur $J$ à droite \teq pour tout
$h\in H$ et tout $j\in J$, $a_h|X_j : X_j \to Y$ se factorise à
travers l'inclusion $X_{jh} \hra Y$, alors pour tout $j\in J$,
l'inclusion $(X_j, G_j) \hra (Y, H)$ est cocartésienne, où $G_j$
est le stabilisateur de $j$ dans~$H$.
\end{sssect}

Soit $(X,G)\in \Eq(\bC)$. On appelle \emph{quotient} de $X$ par~$G$
un objet $Y$ de~$\bC$ muni d'un morphisme $p:X\to Y$ invariant par
$G$ \tq $(X,G)\to (Y,\{1\})$ soit cocartésien, \ie \tq pour tout
objet $Z$ de~$\bC$, l'application
\begin{align*}
\Hom(Y,Z) & \to \Hom(X,Z)^G\\
f&\mapsto f\circ p
\end{align*}
soit bijective. Pour $\al: G\to H$ un homomorphisme surjectif, si le
quotient $Y$ de $X$ par $\Ker \al$ existe, alors l'action de $G$
sur~$X$ induit une action de $H$ sur~$Y$, et la flèche $(X,G)\to
(Y,H)$ est cocartésienne. De plus, le quotient de $X$ par $\Ker
\al$ existe si et seulement s'il existe une flèche cocartésienne
de source $(X,G)$ au-dessus de~$\al$. Cela est vrai pour $\al: G\to
H$ un homomorphisme quelconque si les sommes disjointes finies sont
représentables dans~$\bC$.

\begin{sssect}\label{sss.fibre}
Soit $F: \bF \to \bC$ un foncteur. Pour $(X,G)$ un objet de
$\Eq(\bC)$, on note $F_{(X,G)}$ la catégorie fibre de $\Eq(F) :
\Eq(\bF) \to \Eq(\bC)$ en $(X,G)$. Pour $(\cF,G)$, $(\cG,G)$ deux
objets de $F_{(X,G)}$, $G$ agit à gauche sur $\Hom_{F_X}(\cF,\cG)$
: pour $c\in \Hom_{F_X}(\cF,\cG)$, $g\in G$, $gc \in
\Hom_{F_X}(\cF,\cG)$ est l'unique flèche rendant commutatif le
diagramme suivant
\[\xymatrix{\cF\ar[r]^{a_g}_\sim\ar[d]^{gc} & \cF\ar[d]^{c}\\
\cG\ar[r]^{a_g}_\sim & \cG}
\]
On a une bijection $\Hom_{F_{(X,G)}}((\cF,G), (\cG,G)) \simeq
\Hom_{F_X}(\cF,\cG)^G$.
\end{sssect}
\end{ssect}

\begin{ssect}\label{ss.scheq}
Soit $S$ un \sch noethérien.
On note $\bC_S$ la catégorie des $S$-schémas \tf et on pose $\Eq/S
= \Eq(\bC_S)$.

Soit $(X,G)\in \Eq/S$. Rappelons que l'action de $G$ est dite
\emph{admissible} \cite[V 1.7]{SGA1} si le quotient de $X$ par $G$
existe dans $\bC_S$, ou, ce qui revient au même, si toute
trajectoire de $G$ dans $X$ soit contenue dans un ouvert affine
\ibid[1.8]. Dans ce cas, le morphisme $X\to X/G$ est fini
\ibid[1.5].

Soit $(X,G)\in \Eq/S$. Pour une partie $Z$ de $X$. On note $G_d(Z)$
le stabilisateur de $Z$. Si $x$ est un point de $X$, on appelle
$G_d(x)$ \emph{groupe de décomposition}. Ce groupe opère
canoniquement à gauche sur le corps résiduel $\kappa(x)$, et le
fixateur de $\kappa(x)$ est appelé \emph{groupe d'inertie}, noté
$G_i(x)$. \ibid[2] On dit que l'action de $G$ est \emph{libre} si
elle est admissible et si $G_i(x)$ est trivial pour tout $x\in X$.
Dans ce cas, l'action de $G$ fait de $X$ un $G$-torseur sur $X/G$.

Soit $(X,G)\in \Eq/S$. Si $G$ agit transitivement sur l'ensemble
$\pi_0(X)$ des composantes connexes de $X$, alors pour tout $Y\in
\pi_0(X)$, l'inclusion $(Y,G_d(Y)) \to (X,G)$ est cocartésienne.
\end{ssect}

\subsection{Formalisme $l$-adique équivariant}\label{ss.Dbc}
Supposons que $S$ soit régulier de
dimension $\le 1$. Soit $l$ un nombre premier inversible sur $S$. On
a un formalisme de \faisx $l$-adiques sur les $S$-\schs séparés
\tf \cite[\S~6]{Ekedahl}. Ce formalisme a un sens pour les \schs \tf
sur $S$ (pas nécessairement séparés), et ce n'est que pour
certaines opérations ($Rf_!$ et $Rf^!$) qu'on a besoin d'une \hyp
de séparation sur les morphismes. (Voir \cite{Laszlo-Olsson2} pour
un formalisme sans hypothèse de séparation.) Le formalisme
$l$-adique dans \cite{Ekedahl} marche aussi dans le cadre
équivariant. On va le rappeler brièvement.

On choisit $F: \bF \to \bC_S$ un foncteur bifibré avec $F_X \simeq
(X_{\et}\sptilde)^\op$. Pour $(X,G)\in \Eq/S$, on définit la
catégorie des faisceaux d'ensembles sur $(X,G)$ par $(X,G)\sptilde
= F_{(X,G)}^\op$ (\ref{sss.fibre}), qui est un topos. Un faisceau
sur $(X,G)$ est donc un faisceau $\cF$ sur $X$ muni d'une action
de~$G$ à gauche, compatible à l'action de $G$ sur $X$. En d'autres
termes, $\cF$ est muni d'une famille d'isomorphismes $(b_g : \cF \to
a_{g*} \cF)_{g\in G}$ \teq pour tous $g,h \in G$, le diagramme
\[\xymatrix{\cF \ar[r]^{b_h} \ar[rrd]_{b_{gh}}& a_{h*}\cF \ar[r]^{a_{h*}b_g}
& a_{h*}a_{g*} \cF\ar@{=}[d]\\
&&(a_{gh})_* \cF}
\]
soit \comm. Un morphisme $\cF_1 \to \cF_2$ de faisceaux sur $(X,G)$
est un morphisme $c: \cF_1 \to \cF_2$ de \faisx sur $X$, $G$-\eqrt,
\ie \tq, pour tout $g\in G$, le \diag
\[\xymatrix{\cF_1 \ar[r]^{b_g}\ar[d]_{c} & a_{g*}\cF_1\ar[d]^{a_{g*} c} \\
\cF_2 \ar[r]^{b_g} & a_{g*} \cF_2}\] soit \comm.

Soit $\Lbd$ un anneau local artinien annulé par une puissance de
$l$.  On note $\Mod(X,G,\Lbd)$ la catégorie des faisceaux de
$\Lbd$-modules sur $(X,G)\sptilde$. Un $\Lbd$-module sur $(X,G)$ est
\emph{constructible} si le $\Lbd$-module sur $X$ sous-jacent est
constructible. On note $\Mod_c(X,G,\Lbd)$ la sous-catégorie pleine
de $\Mod(X,G,\Lbd)$ formée des $\Lbd$-modules constructibles sur
$(X,G)$. On note $D(X,G,\Lbd)$ la catégorie dérivée de
$\Mod(X,G,\Lbd)$, $D_c(X,G,\Lbd)$ la sous-catégorie pleine formée
des complexes à faisceaux de cohomologie dans $\Mod_c(X,G,\Lbd)$.

Soient $E$ un corps extension finie de $\Q_l$, $\cO$ son anneau des
entiers, $\fm$ l'idéal maximal de $\cO$. On considère le système
projectif $\cO_\bullet = (\cO/\fm^n)_{n\in \N}$. On pose
\[\Mod_c(X,G,\cO) = 2\text{-}\varprojlim_{n} \Mod_c(X,G,\cO/\fm^n),\]
qui s'identifie à une sous-catégorie pleine de la catégorie
$\Mod(X,G,\cO_\bullet)$ des systèmes projectifs $(M_n)_{n\in \N}$,
où $M_n\in \Mod(X,G,\cO/\fm^n)$ et $M_n \to M_m$ est une flèche
dans $\Mod(X,G,\cO/\fm^n)$ pour $n \ge m$. On note
$D(X,G,\cO_\bullet)$ la catégorie dérivée de
$\Mod(X,G,\cO_\bullet)$. Un faisceau $\cF\in\Mod(X,G,\cO_\bullet)$
est \emph{essentiellement nul} si pour tout $n\in \N$, il existe
$m\ge n$ tel que la flèche $\cF_m \to \cF_n$ soit nulle. Un
complexe $K \in D(X,G,\cO_\bullet)$ est \emph{essentiellement nul}
si $\fH^i K$ est essentiellement nul pour tout  $i\in \Z$ ; $K$ est
\emph{essentiellement constant} s'il existe un complexe $C$ de faisceaux
de $\cO$-modules de torsion sur $(X,G)\sptilde$ et des flèches $L\to K$,
$L\to (C\otimes_{\cO} \cO/\fm^n)_{n\in \N}$ dans
$D(X,G,\cO_\bullet)$ à cônes essentiellement nuls. On considère
le foncteur
\begin{align*}
D^b(X,G,\cO_\bullet) &\to D^-(X,G,\cO_\bullet)\\
K &\mapsto (\cO/\fm)_{n\in \N} \otimes^L_{\cO_\bullet} K.
\end{align*}
On note $D^b(X,G,\cO)$ le quotient de la sous-catégorie pleine de
$D^b(X,G,\cO_\bullet)$, image inverse des complexes essentiellement
constants, par la sous-catégorie épaisse, image inverse des
complexes essentiellement nuls. Le foncteur
\[(\cO/\fm)\otimes^L_{\cO} R\varprojlim - : D^b(X,G,\cO_\bullet) \to D^+(X,G,\cO/\fm)\]
induit un foncteur conservatif \cite[2.6, 2.7]{Ekedahl}
\[(\cO/\fm)\otimes^L_{\cO} - : D^b(X,G,\cO) \to D^b(X,G,\cO/\fm).\]
On note $\Dbc(X,G,\cO)$ la sous-catégorie pleine de $D^b(X,G,\cO)$,
image inverse de $\Dbc(X,G, \cO/\fm)$. Elle est munie d'une
$t$-structure \can $(D^{\le 0}, D^{\ge 0})$, dont le c\oe ur est
$\Mod_c(X,G,\cO)$ \cite[3.5]{Ekedahl}.

On pose $\Mod_c(X,G,E) = \Mod_c(X,G,\cO) \otimes_\cO E$,
$\Dbc(X,G,E) = \Dbc(X,G,\cO) \otimes_\cO E$. Enfin on pose
$\Mod_c(X,G,\Qlb) = 2\text{-}\varinjlim_{E} \Mod_c(X,G,E)$,
$\Dbc(X,G,\Qlb) = 2\text{-}\varinjlim_{E} \Dbc(X,G,E)$, où $E$
parcourt les extensions finies de $\Q_l$ contenues dans une
clôture algébrique $\Qlb$ de $\Q_l$. On note $K(X,G,\Qlb)$ le
groupe de Grothendieck de $\Mod_c(X,G,\Qlb)$, qui est aussi le
groupe de Grothendieck de $\Dbc(X,G,\Qlb)$.

Soit $R= \Lbd, \cO, E$ ou $\Qlb$ comme plus haut. On choisit $F: \bF
\to \bC_S$ un foncteur bifibré avec $F_X \simeq \Mod_c(X,R)^\op$
(\resp $F_X \simeq \Dbc(X,R)^\op$). Pour $(X,G)\in \Eq/S$, on pose
$\Mod_c(X,G,R)_\naif = F_{(X,G)}^\op$ (\resp $\Dbc(X,G,R)_\naif =
F_{(X,G)}^\op$) (\ref{sss.fibre}). On a une équivalence de
catégories
\[\Mod_c(X,G,R) \simto \Mod_c(X,G,R)_\naif\] (\resp un foncteur
\beq\label{eq.Dbc}
\Dbc(X,G,R) \to \Dbc(X,G,R)_\naif
\eeq
$t$-exact pour les $t$-structures canoniques.)

\begin{prop}
Le foncteur \eqref{eq.Dbc} est une équivalence de catégories si
$R= E$ ou $\Qlb$.
\end{prop}

\begin{proof}
Soient $K,L \in \Dbc(X,G,R)$. On note $f: (\pt, G) \to (\pt,
\{1\})$. Le foncteur $R\Hom_{\Dbc(X,G,R)}(K,-)$ s'identifie au
composé
\[\Dbc(X,G,R) \sto{R\Hom_{\Dbc(X,R)}(K,-)} \Dbc(\pt,G,R) \sto{f_*} \Dbc(\pt,R). \]
La suite spectrale de foncteur composé appliquée à $L$ s'écrit
\[E_2^{pq} = H^p (G, \Hom_{\Dbc(X,R)}(K,L[q]) ) \Rightarrow \Hom_{\Dbc(X,G,R)}(K,L[p+q]).\]
Comme $H^p (G,M) = 0$ pour tout $p > 0$ et tout $R$-espace vectoriel
$M$, on a
\[\Hom_{\Dbc(X,G,R)}(K,L) \simeq \Hom_{\Dbc(X,R)}(K,L)^G.\]
Donc \eqref{eq.Dbc} est pleinement fidèle. Pour $K\in
D^{[a,b]}(X,G,R)_\naif$, montrons par récurrence sur $b-a$ que K
est dans l'image essentielle de \eqref{eq.Dbc}. Pour $b-a\le 0$, $K=
\cF[-a]$ pour $\cF \in \Mod_c(X,G,R)_\naif$. Pour $b-a\ge 1$, $K$
est le cône de $(\tau_{\ge b} K)[-1] \to \tau_{\le b-1}K$. Il
suffit alors d'appliquer l'hypothèse de récurrence à $\tau_{\ge
b} K \in D^{[b,b]}$ et à $\tau_{\le b-1} K \in D^{[a,b-1]}$.
\end{proof}


On renvoie à \ref{prop.omi} pour une interprétation de
$\Dbc(X,G,\Qlb)$ (et des six opérations qu'on va définir) en
termes du champ de Deligne-Mumford $[X/G]$.

\begin{ssect}\label{ss.6op}
Pour $K,L \in \Dbc(X,G,\Qlb)$, $G$ agit sur $K\otimes L$ et $R
\cHom(K,L)$ par
\begin{align*}
 & a_g^*(K\otimes L) \simto a_g^* K \ot a_g^* L \sto{b_g\otimes b_g}
K\ot L,\\
 & a_g^* R\cHom(K,L) \simto R\cHom(a_g^*K, a_g^*
L) \sto{R\cHom(b_g^{-1}, b_g)} R\cHom(K,L), \quad g\in G.
\end{align*}
Cela définit des foncteurs bi-exacts
\begin{align*}
-\ot- &: \Dbc(X,G,\Qlb) \times \Dbc(X,G,\Qlb) \to \Dbc(X,G,\Qlb),\\
R\cHom(-,-) &: \Dbc(X,G,\Qlb)^{\op} \times \Dbc(X,G,\Qlb) \to
\Dbc(X,G,\Qlb).
\end{align*}

Soit $(f,\al): (X_1,G_1) \to (X_2, G_2)$ un morphisme dans $\Eq/S$
(\resp un morphisme dans $\Eq/S$ avec $f$ séparé). Pour $L\in
\Dbc(X_2,G_2,\Qlb)$, $G_1$ agit sur $f^*L$ (\resp $Rf^! L$) par
\begin{align*}
& a_g^* f^* L \simeq f^* a_{\al(g)}^* L \sto{f^* b_g} L\\
\text{(\resp}\quad & a_g^* Rf^! L \simeq Rf^! a_{\al(g)}^* L
\sto{Rf^! b_g} L\text{)}, \quad g\in G_1.
\end{align*}
Cela définit un foncteur exact
\[(f,\al)^* \text{ (\resp $R(f,\al)^!$)} : \Dbc(X_2,G_2,\Qlb) \to
\Dbc(X_1, G_1,\Qlb).
\]
Si $(g,\beta): (X_2,G_2) \to (X_3, G_3)$ est un deuxième morphisme
dans $\Eq/S$ (\resp un deuxième morphisme dans $\Eq/S$ avec $g$
séparé), alors $(gf, \beta \al)^* \simeq (f,\al)^* (g,\beta)^*$
(\resp $R(gf,\beta\al)^! \simeq R(f,\al)^! R(g,\beta)^!$).

On va définir un foncteur exact
\[R(f,\al)_* \text{ (\resp $R(f,\al)_!$)} : \Dbc(X_1,G_1,\Qlb) \to
\Dbc(X_2, G_2,\Qlb),
\]
adjoint à droite (\resp à gauche) de $(f,\al)^*$ (\resp
$R(f,\al)^!$). Traitons d'abord trois cas spéciaux.

(i) Cas $G_1=G_2=G$, $\al=\Id_G$. Pour $K\in \Dbc(X_1,G,\Qlb)$, on
pose $R(f,\Id_G)_*K = Rf_* K$ (\resp $R(f,\Id_G)_!K = Rf_! K$), sur
lequel $G$ agit par
\begin{align*}
& Rf_* K \sto{Rf_* b_g} Rf_* a_{g*}K = a_{g*} Rf_* K\\
\text{(\resp}\quad & Rf_! K \sto{Rf_! b_g} Rf_! a_{g*}K \simeq
a_{g*} Rf_! K \text{)}, \quad g\in G.
\end{align*}
Le foncteur $R(f,\Id_G)_*$ (\resp $R(f,\Id_G)_!$) ainsi défini est
un adjoint à droite (\resp à gauche) de $(f,\Id_G)^*$ (\resp
$R(f,\Id_G)^!$).

(ii) Cas $X_1=X_2=X$, $f=\Id_X$, $\al$ surjectif.
Pour $\cF\in \Mod(X,G_1,\cO/\fm^n)$,
on définit $(\Id_X,\al)_* \cF$ (\resp $(\Id_X,\al)_!\cF $) comme le
sous-faisceau de $\cF$ des invariants (\resp le faisceau quotient
de~$\cF$ des coïnvariants) par $\Ker \al$. Il est muni d'une
action de $G_2$ induite par l'action de $G_1$ sur~$\cF$. Le foncteur
$(\Id_X,\al)_*$ (\resp $(\Id_X,\al)_!$) $: \Mod(X,G_1,\cO_\bullet)
\to \Mod(X,G_2,\cO_\bullet)$ ainsi obtenu est un adjoint à droite
(\resp à gauche) de $(\Id_X,\al)^*$, donc est exact à gauche
(\resp à droite) et se dérive en $R(\Id_X,\al)_* :
D^+(X,G_1,\cO_\bullet) \to D^+(X,G_2,\cO_\bullet)$ (\resp
$L(\Id_X,\al)_! : D^-(X,G_1,\cO_\bullet) \to
D^-(X,G_2,\cO_\bullet)$). Cela induit un foncteur $(\Id_X,\al)_*$
(\resp $(\Id_X,\al)_!$) $: \Dbc(X,G_1,\Qlb) \to \Dbc(X,G_2,\Qlb)$
adjoint à droite (à gauche) de $(\Id_X,\al)^*$. La flèche
canonique $\cF^{\Ker \al} \to \cF_{\Ker \al}$ induit un isomorphisme
de foncteurs $(\Id_X,\al)_* \simto (\Id_X,\al)_! : \Dbc(X,G_1,\Qlb)
\to \Dbc(X,G_2,\Qlb)$, donc ces foncteurs sont $t$-exacts pour les
$t$-structures canoniques.

(iii) Cas $X_1=X_2=X$, $f=\Id_X$, $\al$ injectif. On choisit un
système $S$ de représentants de $G_2/\al(G_1)$. Pour $K\in
\Dbc(X,G_1,\Qlb)$, on considère $L = \bigoplus_{s\in S} a_s^* K$.
Soit $g\in G_2$. Pour $s\in S$, on pose $g s = t_s \al(h_s)$, $t_s\in S$,
$h_s\in G_1$. Alors $g$ agit sur $L$ par
\[a_g^* \bigoplus_{s\in S} a_s^*
K \simeq \bigoplus_{s\in S} a_{t_s}^* a_{h_s}^* K
\sto{\bigoplus_{s\in S}a_{t_s}^* b_{h_s}} \bigoplus_{s\in S}
a_{t_s}^* K.
\]
On pose $(\Id_X,\al)_* K = (\Id_X,\al)_! K = L$. Cette définition
ne dépend pas du choix de $S$ à isomorphisme près. Le foncteur
$t$-exact
\[(\Id_X,\al)_* = (\Id_X,\al)_! : \Dbc(X,G_1,\Qlb) \to
\Dbc(X,G_2,\Qlb)\] ainsi défini est à la fois un adjoint à droite
et un adjoint à gauche de $(\Id_X,\al)^*$.

Dans le cas général, on décompose $(f,\al)$ en
\[(X_1,G_1)
\sto{(f,\Id_{G_1})} (X_2,G_1) \sto{(\Id_{X_2},\al)} (X_2,G_2)
\] et
$\alpha$ en $G_1 \sto{\al_1} \Img \al \sto{\al_2} G_2$. On pose
$R(f,\al)_* = (\Id_{X_2},\al_2)_* (\Id_{X_2},\al_1)_*
R(f,\Id_{G_1})_*$ (\resp $R(f,\al)_! = (\Id_{X_2},\al_2)_!
(\Id_{X_2},\al_1)_! R(f,\Id_{G_1})_!$). Le foncteur $R(f,\al)_*$
(\resp $R(f,\al)_!$) ainsi défini est un adjoint à droite (\resp
à gauche) de $(f,\al)^*$ (\resp $R(f,\al)^!$). Il s'en suit que
$R(gf,\beta\al)_* \simeq R(g,\beta)_*R(f,\al)_*$ (\resp
$R(gf,\beta\al)_! \simeq R(g,\beta)_! R(f,\al)_!$). Voir aussi
\cite[XVII 3.3.2]{SGA4}.

On enlève $\al$ de la notation lorsqu'il n'y a pas de confusion à
craindre. Les six opérations sur $\Dbc$ induisent les opérations
correspondantes sur les groupes de Grothendieck.

Soit $a_X : (X,G) \to (S,\{1\})$. On vérifie que $Ra_X^!\Q_l$ est
globalement défini (le cas général découle du cas séparé par
\cite[3.2.4]{BBD}). On pose $D_X = R\cHom_X(-,Ra_X^!\Q_l)$. On a
$D_X D_X \simeq \Id$. Pour $K,L\in \Dbc(X,G,\Qlb)$, on a
$D_XR\cHom_X(K,L) = K \otimes D_X L$. Si $(f,\al):(X,G) \to (Y,H)$
est un morphisme dans $\Eq/S$ avec $f$ \sep, alors $D_X f^* \simeq
Rf^! D_Y$ et $D_Y R(f,\al)_* \simeq R(f,\al)_! D_X$. En effet, tous
les énoncés sauf le dernier résultent facilement des analogues
sans actions de groupe. Pour le dernier énoncé, on remarque que
pour tout $K\in \Dbc(X,G,\Qlb)$ et tout $L\in \Dbc(Y,H,\Qlb)$, on a
\begin{multline*}
\Hom_Y(D_Y Rf_* K,L) \simeq \Hom_Y(D_Y L,Rf_*K) \simeq \Hom_X(f^*D_Y
L,K)
\simeq \Hom_X(D_X Rf^!L,K)\\
\simeq \Hom_X(D_X K, Rf^!L) \simeq \Hom_Y(Rf_! D_X K, L).
\end{multline*}

Pour $(X,G)\in \Eq/S$, $K\in \Dbc(X,G,\Qlb)$ et $g\in G$, $b_g$
induit un \iso $K \simto T_{g*} K$, où $T_g=(a_g, h\mapsto
g^{-1}hg)$ comme dans \ref{ss.cat}.

\begin{sssect}\label{sss.cocarteb}
Si $(f,\al) : (X,G) \to (Y,H)$ est cocartésien avec $\al$ injectif,
alors $(f,\al)_*$ est une équivalence de catégories. Cela résulte
de \eqref{eq.aleb} et de la définition de $(f,\al)_*$.
\end{sssect}
\end{ssect}

\begin{prop}\label{prop.cb}
Soient $(f,\al) : (X,G) \to (Y,H)$, $(g,\beta): (Y',H') \to (Y,H)$
deux morphismes de même but dans $\Eq/S$ avec $f$ séparé. Pour
$r\in H$, formons le \carrecart
\beq\label{eq.dc.cb}
\xymatrix{(X', G')_r \ar[rr]^{(h,\gamma)_r} \ar[d]^{(f',\al')_r}
&& (X, G)\ar[d]^{(f,\al)}\\
(Y',H')\ar[r]^{(g,\beta)} & (Y,H)\ar[r]^{T_r} &(Y,H)}
\eeq
où $T_r$ est comme plus haut. Le foncteur $R(f',\al')_{r!}
(h,\gamma)_r^*$ ne dépend, à isomorphisme près, que de la double
classe $(\Img \beta) r (\Img \al) \subset H$ et on a
\[(g,\beta)^* R(f,\al)_!
\simeq \bigoplus_{r} R(f',\al')_{r!} (h,\gamma)_r^*,\] où $r$
parcourt un système de représentants de $\Img \beta\backslash
H/\Img \al$. En particulier, $R(f,\al)_!$ commute à tout changement
de base $(g,\beta)$ dans $\Eq/S$ \ver ant $\sharp (\Img
\beta\backslash H/\Img \al) = 1$.
\end{prop}

Le cas où $f=g=\Id_Y$ avec $Y$ le spectre d'un corps séparablement
clos est une version de la formule de Mackey \cite[7.3]{Serre}.

On renvoie à \ref{prop.Mackeychamps} pour une interprétation en
termes de champs.

\begin{proof}
On remarque d'abord que $R(f,\al)_!$ commute au changement de base
$(g,\beta)$ dans chacun des trois cas suivants : (i) $\al=\Id$ ;
(ii) $f=\Id$ et $\al$ surjectif ; (iii) $f=\Id$, $\al$~injectif et
$\beta$ surjectif. Cela résulte du théorème de changement de base
usuel \cite[XVII 5.2.6]{SGA4} dans le cas (i), et de la définition
de $(\Id,\al)_!$ dans les cas (ii) et (iii).

On factorise $\al$ en $G\sto{\al_1} \Img \al \sto{\al_2} H$, $\beta$
en $H'\sto{\beta_1} \Img \beta \sto{\beta_2} H$, et \eqref{eq.dc.cb}
en un diagramme à carrés cartésiens
\[
\xymatrix{(X',G')_r \ar[rrr]^{(h,\gamma)_r}\ar[d] &&& (X,G)\ar[d]^{(f,\al_1)}\\
(Y',I_r')\ar[r]^{(g,\sigma_r)} \ar[d]^{(\Id_{Y'},\al'_{2,r})}
& (Y,I_r)\ar[rr]^{(a_r,\tau_r)}\ar[d]^{(\Id_Y,\rho_r)} && (Y,\Img \al)\ar[d]^{(\Id_Y,\al_2)}\\
(Y',H') \ar[r]^{(g,\beta_1)} & (Y,\Img \beta)
\ar[r]^{(\Id_Y,\beta_2)} & (Y,H)\ar[r]^{T_r} & (Y,H)}
\]
On a
\begin{alignat*}{2}
R(f',\al')_{r!} (h,\gamma)_r^* &\simeq (\Id_{Y'},\al'_{2,r})_!
(g,\sigma_r)^* (a_r,\tau_r)^* R(f,\al_1)_! &&\qquad \text{d'après (i) et (ii)}\\
&\simeq (g,\beta_1)^* (\Id_Y,\rho_r)_! (a_r,\tau_r)^* R(f,\al_1)_!
&&\qquad \text{d'après (iii)}.
\end{alignat*}

Donc pour montrer la proposition, \ops que $f=g=\Id_Y$ et que $\al$
et $\beta$ sont des inclusions. On peut prendre pour
\eqref{eq.dc.cb} le diagramme suivant :
\[\xymatrix{(Y,H'\cap rGr^{-1}) \ar[rr]^{(a_r,\gamma_r)}\ar[d]^{(\Id,\al'_r)} &&
(Y,G)\ar[d]^{(\Id,\al)}\\
(Y,H')\ar[r]^{(\Id,\beta)} & (Y,H)\ar[r]^{T_r} &(Y,H)}
\]
où $\al'_r$ est l'inclusion, $\gamma_r$ est donné par $g\mapsto
r^{-1} g r$. Soit $K \in \Dbc(Y,G,\Qlb)$. Soit $R$ un système de
représentants de $H'\backslash H/G$. Pour $r\in R$, soit $S_r$ un
système de représentants de $H'/(H'\cap rGr^{-1})$. Posons $L_r =
\bigoplus_{s\in S_r} a_s^* a_r^*K \simeq
(\Id,\al'_r)_!(a_r,\gamma_r)^*K$. Comme $\bigcup_{r\in R} S_r r$ est
un système de représentants de $H/G$, on a
\beq\label{eq.prop.cb}
\bigoplus_{r\in R} L_r \simeq \bigoplus_{r\in R}\bigoplus_{s\in S_r}
a_{sr}^* K \simto (\Id,\beta)^* (\Id,\al)_! K.
\eeq
Pour $h\in H'$, si on pose $hs = t_s h_s$, $t_s\in S_r$, $h_s\in
H'\cap rGr^{-1}$, alors $h$ agit sur $L_r$ par
\[a_h^*\bigoplus_{s\in S_r} a_s^* a_r^* K \simeq
\bigoplus_{s\in S_r} a_{t_s}^* a_{h_s}^* a_r^* K
\sto{\bigoplus_{s\in S_r} a_{t_s}^* b_{a_r^* K, h_s}}
\bigoplus_{s\in S_r} a_{t_s}^* a_r^* K,
\]
où $b_{a_r^* K, h_s}$ est donné par
\[a_{h_s}^* a_r^* K \simeq a_r^* a_{\gamma_r(h_s)}^* K \sto{a_r^* b_{\gamma_r(h_s)}} a_r^* K.\]
Comme $h s r = t_s r r^{-1} h_s r$, $t_s r \in S_r r$, $r^{-1} h_s r
\in G$, $h$ agit sur $(\Id,\beta)^* (\Id,\al)_! K$ par
\[a_h^* \bigoplus_{r\in R} \bigoplus_{s\in S_r}a_{sr}^* K
\simeq \bigoplus_{r\in R} \bigoplus_{s\in S_r}a_{t_s r}^* a_{r^{-1}
h_s r}^* K \sto{a_{t_s r}^*(b_{r^{-1}h_s r})} \bigoplus_{r\in R}
\bigoplus_{s\in S_r} a_{t_s r}^* K.
\]
Donc l'isomorphisme \eqref{eq.prop.cb} est $H'$-équivariant.
L'image de $L_r$ ne dépend que de la double classe $H' r G$.
\end{proof}

\begin{prop}\label{prop.quot} Soit $(f,\al): (X,G) \to (Y,H)$ un morphisme dans $\Eq/S$.

(a) Supposons que $\al : G \to H$ soit surjectif et que $f$ fasse de
$Y$ un quotient de $X$ par $\Ker \al$. Alors pour tout $K\in
\Dbc(Y,H,\Qlb)$, la flèche d'adjonction $K \to (f,\al)_* (f,\al)^*
K$ est un isomorphisme.

(b) Si de plus $\Ker \al$ opère librement sur $X$, de sorte que $X$
est un $(\Ker \al)$-torseur sur~$Y$, alors pour tout $L\in
\Dbc(X,G,\Qlb)$, la flèche d'adjonction $(f,\al)^* (f,\al)_*L \to
L$ est un isomorphisme.
\end{prop}

On renvoie à \ref{ss.omi} pour une interprétation en termes de
champs dans le cas (b).

\begin{proof}
(a) Pour tout point géométrique $\yb$ de $Y$ à valeur dans un
corps algébriquement clos, $\yb$ est le quotient de $X_\yb^{\red}$
par $\Ker \al$ (ceci découle de \cite[V 1.1 (iii)]{SGA1}). Quitte
à changer $S$, \ops $H=\{e\}$, $Y=y=\Spec k$ avec $k$ un corps
\algt clos et $X$ réduit. Alors il existe un sous-groupe $G_1$
de~$G$ \tq $X=\coprod_{h\in G_1\backslash G} x_h$, où $x_h=\Spec k$
et l'action de $G$ est donnée par $a_g(x_h) = x_{hg}$. Soit $K\in
\Dbc(Y,\Qlb)$. Alors $((f,\al)^* K)_{x_h} = K$, et $g\in G$ agit par
\[a_g^* (f,\al)^* K \simeq (f,\al)^* K \sto{\Id} (f,\al)^* K.\]
Donc $(f,\Id)_* (f,\al)^* K \simeq \bigoplus_{h\in G/G_1} K_h$, où
$K_h =K$, et $g\in G$ agit par
\[\bigoplus_{h\in G/G_1} (K_h
\sto{\Id} K_{gh}).
\]
Donc la diagonale induit un isomorphisme $K \to
(\bigoplus_{h\in G/G_1} K_h)^G \simeq (f,\al)_* (f,\al)^* K$.

(b) Si $\Ker \al$ opère librement sur $X$, alors pour tout point
géométrique $\yb \to Y$, $X_\yb$ est un $\Ker \al$-torseur sur
$\yb$. On fait le même dévissage qu'en (a). Dans notre cas $G_1
=\{1\}$. Soit $(g,\beta)$ une section de $(f,\al)$ (par exemple
$(y,1) \mapsto (x_1,1)$). Alors $(g,\beta)$ est cocartésien, donc
$(g,\beta)_*$ est une équivalence de catégories. Il en résulte
que $(f,\al)_*$ l'est aussi, d'où le résultat.
\end{proof}

\begin{remqe}

Soient $k/k_0$ une extension finie galoisienne, $k'$ une extension
finie de $k$ \tq $k'$ soit galoisien sur $k_0$. Alors $(\Spec
k',\Gal(k'/k_0)) \to (\Spec k, \Gal(k/k_0))$ vérifie les
hypothèses de (b).

\end{remqe}

\subsection{Traces locales}

Soit $G$ un groupe fini agissant sur $x=\Spec \kappa(x)$. Le corps
$\kappa(x)$ est muni d'une action de~$G$ à gauche. Soit $\xb$ un
point \geoq algébrique au-dessus de $x$, \ie le spectre d'un corps
$\kappa(\xb)$, clôture séparable de $\kappa(x)$. On pose
\[G_\xb
= G \times_{\Gal(\kappa(x)/\kappa(x)^{G})}
\Gal(\kappa(\xb)/\kappa(x)^{G}).\]

Soit $l$ un nombre premier inversible sur~$x$. On a une équivalence
de \cats
\begin{align*}
\Mod_c(x,G,\Qlb) &\to \Rep(G_\xb,\Qlb) \\
L &\mapsto L_\xb,
\end{align*}
où $\Rep(G_\xb,\Qlb)$ est la \cat des $\Qlb$-\repns continues de
$G_\xb$. Pour $L\in \Mod_c(x,G,\Qlb)$, $(DL)_\xb$ s'identifie à la
représentation contragrédiente de $L_\xb$.

Soient $H$ un groupe fini agissant sur  $y= \Spec \kappa(y)$ et
$(f,\al): (x,G) \to (y, H)$ un morphisme avec $f$ fini. Prenons
$\yb$ un point géométrique algébrique au-dessus de $y$
s'insérant dans un carré commutatif
\[
\xymatrix{
  \xb \ar[d] \ar[r] & x \ar[d]^f \\
  \yb \ar[r] & y   }
\]
À l'aide de $f$, on identifie $\kappa(y)$ à un sous-corps de
$\kappa(x)$. On a $\kappa(y)^{H} \subset \kappa(y)^{G} \subset
\kappa(x)^{G}$, d'où un \diagcomm de groupes
\[\xymatrix{
G \ar[r]\ar[d]_\al & \Gal(\kappa(x)/\kappa(x)^{G})\ar[d] &
\ar[l]\ar[d] \Gal(\kappa(\xb)/\kappa(x)^{G})\\
H \ar[r] & \Gal(\kappa(y)/\kappa(y)^{H}) & \ar[l]
\Gal(\kappa(\yb)/\kappa(y)^{H})}
\]
qui induit un homomorphisme $\al_\xb = (f,\al)_\xb : G_\xb \to
H_\yb$. Soit $L\in \Mod_c(X,G,\Qlb)$. On a $((f,\al)_* L)_\yb \simeq
(\al_\xb)_* L_\xb$ comme $H_\yb$-représentations. Ici
$(\al_{\xb})_*$ est la co-induction. (Si on note $(\al_{\xb})_!$
l'induction, on a un isomorphisme canonique $(\al_{\xb})_* \simto
(\al_{\xb})_!$. \Cf \ref{ss.6op} (ii) et (iii).) Le caractère de la
représentation $(\al_{\xb})_*$ peut se calculer en utilisant le
lemme suivant.

\begin{lemm}\label{lemm.traces}
Soient $F$ un corps, $\al : G_1 \to G_2$ un homomorphisme de groupes
(abstraits)
 avec $N=\Ker \al$ d'ordre fini et $I=\Img \al$
d'indice fini dans $G_2$, $L$ un $F$-\ev \df, $\rho: G_1 \to
\GL_F(L)$ une représentation \lin. On note $\al_* L$ le module
co-induit, \ie $\al_* L = \Hom_{F[I]}(F[G_2],L^N)$ avec l'action de
$G_2$ induite par multiplication à droite sur $F[G_2]$. Alors pour
tout $g\in G_2$,
\beq\label{eq.traces}
\sharp N \cdot \Tr(g,\al_* L) = \sum_{s} \sum_{\substack{ t\in G_1 \\
\al(t)= s^{-1}gs}} \Tr(t,L),
\eeq
où $s$ parcourt un système de représentants de $G_2/I$.
\end{lemm}


\begin{proof} Il suffit de traiter les deux cas spéciaux suivants.

(i) Cas $N=\{1\}$. La formule est immédiate (voir, par exemple,
\cite[3.3]{Serre}).

(ii) Cas $I=G_2$. Soit $\gb = \sum_{\substack{t\in G_1 \\
\al(t)=g}} \rho(t) \in \End_F(L)$. Alors $\gb(L)\subset L^N$. Le
membre de droite de \eqref{eq.traces} est
\[\Tr(\gb,L) = \Tr(\gb, L^N) + \Tr(\gb, L/L^N) = \sharp N \cdot \Tr(g,L^N).\]
\end{proof}

\begin{ssect}\label{ss.traces}


Soit $(X,G)\in \Eq/S$. 
Pour tout point $x\in X$, $\kappa(x)$ est muni d'une action à
gauche du groupe de décomposition
\[G_d(x) = \ensdr{g\in G}{a_g(x) = x}.
\]
On désigne le \sch $\Spec \kappa(x)$ encore par $x$ et on note
$i_x$ le morphisme $(x,G_d(x))\to (X,G)$. Soit $\xb$ un point \geoq
algébrique au-dessus de~$x$. On pose
\beq\label{eq.Gdxxb}
G_d(x,\xb) = G_d(x)_\xb = G_d(x)
\times_{\Gal(\kappa(x)/\kappa(x)^{G_d(x)})}
\Gal(\kappa(\xb)/\kappa(x)^{G_d(x)}).
\eeq
En d'autres termes, un éléments de $G_d(x,\xb)$ est un couple
$(g,\phi)$ où $g\in G_d(x)$ et $\phi\in \Aut(\kappa(\xb))$
induisent le même automorphisme sur $\kappa(x)$.

Soient $(f,\al): (X,G) \to (Y, H)$ un morphisme dans $\Eq/S$, $x\in
X$, $y=f(x)$, $\xb \to x$ au-dessus de $\yb \to y$. Alors
$(f,\al)$ induit $(f_x, \al_x) : (x,G_d(x)) \to (y,H_d(y))$ et
$\al_\xb = (f,\al)_\xb : G_d(x,\xb) \to H_d(y,\yb)$.

Soit $L\in \Dbc(X,G,\Qlb)$. On pose $L_x = i_x^* L \in
\Dbc(x,G_d(x),\Qlb)$.
Pour $g\in G$, $b_g$ induit $L_{a_g x} \simeq (T_{g,x})_* L_x$, où
$T_g=(a_g,h\mapsto g^{-1}hg)$ comme dans \ref{ss.cat}.

Soit $f: (X,G) \to (Y,H)$ un morphisme fini. Soit $y\in Y$. Notons
$O_y$ l'orbite de~$y$ sous~$H$. On choisit un système~$S$ de
représentants des orbites dans $X$ sous~$G$ au-dessus de $O_y$.
Pour $x\in S$, notons $O_x$ l'orbite de $x$ sous~$G$, $f_{O_x}: O_x
\to O_y$ la restriction de $f$ à $O_x$, et $z=f(x)$. Choisissons
$h_x\in H$ \tq $a_{h_x} (y) = z$. Alors $((f,\al)_*L)|O_y \simeq
\bigoplus_{s\in S} (f_{O_x},\al)_* (L|O_x)$. Donc
\[
((f,\al)_* L)_y \simeq \bigoplus_{x\in S} ((f_{O_x},\al)_*(L|O_x))_y
\simeq \bigoplus_{x\in S} T_{h_x,y}^*((f_{O_x},\al)_*(L|O_x))_{z}
\simeq \bigoplus_{x\in S} T_{h_x,y}^* (f_x,\al_x)_* L_x.
\]
Prenons $\xb \to x$ et $\zb\to z$ au-dessus de $\yb \to y$. Alors
pour tout $h\in H_d(y,\yb)$,
\begin{align}
\notag \Tr(h,((f,\al)_* L)_\yb) &= \sum_{x\in S} \Tr(T_{h_x,\yb}(h),
(\al_\xb)_* L_\xb) \\
\label{eq.ftraces} &= \sum_{x\in S} \frac{1}{\sharp N_x} \sum_{s\in
H_d(z,\zb)/I_x} \sum_{\substack{t\in G_d(x,\xb)\\\al_\xb(t) =
s^{-1}T_{h_x,\yb}(h) s}} \Tr(t,L_\xb), \quad\text{d'après
\eqref{eq.traces}}
\end{align}
où $N_x = \Ker \al_\xb$, $I_x = \Img \al_\xb$.

\end{ssect}

\begin{prop}\label{prop.quot2}
On utilise les notations de \ref{prop.quot}. Soient $x\in X$,
$y=f(x)$, $\yb \to y$ au-dessus de $\xb\to x$. Dans le cas (a)
(\resp (b)) de \ref{prop.quot}, l'homomorphisme $\al_\xb :
G_d(x,\xb) \to H_d(y,\yb)$ est surjectif (\resp un isomorphisme).
\end{prop}

\begin{proof}
D'après \cite[V 1.1 (iii)]{SGA1}, $\kappa(x)$ est une extension
quasi-galoisienne de $\kappa(y)$ et l'homomorphisme $\Ker \al_x \to
\Gal(\kappa(x)/\kappa(y))$ est surjectif. Donc $\kappa(x)^{G_d(x)}$
est une extension radicielle de $\kappa(y)^{H_d(y)}$, d'où un
carré cartésien de groupes
\[\xymatrix{
  G_d(x,\xb) \ar[d] \ar[r]^{\al_{\xb}} & H_d(y,\yb) \ar[d] \\
  G_d(x) \ar[r] & H_d(y) \times_{\Gal(\kappa(y)/\kappa(y)^{H_d(y)})} \Gal(\kappa(x)/\kappa(x)^{G_d(x)})  } \]
Donc $\al_\xb$ est surjectif.

Dans le cas (b), l'homomorphisme $\Ker \al_x \to
\Gal(\kappa(x)/\kappa(y))$ est un isomorphisme. Donc $\al_\xb$ est
un isomorphisme.
\end{proof}

\begin{coro}\label{prop.DL1}
  Soient $K$ un corps, $K'$ une extension finie galoisienne, $\eta= \Spec
  K$, $\eta'= \Spec K'$, $(X,G)\in \Eq/\eta$. Considérons $(\pr_1,\Id) : (X_{\eta'},G) \to
  (X,G)$.  Soient $x\in X$, $y\in X_{\eta'}$ au-dessus de $x$, $\yb \to y$ au-dessus de $\xb\to x$.
  Alors l'image de $\Id_\yb : G_d(y,\yb) \to G_d(x,\xb)$ est
  \[G'=G_d(x)
  \times_{\Gal(\kappa(x)/\kappa(x)^{G_d(x)})} \Gal(\kappa(\xb)/\kappa(x)^{G_d(x)} \cdot
  K').\]
\end{coro}

\begin{proof}
  Il est clair que $\Img \Id_\yb \subset G'$. Soit $(g,\phi)\in G'$.
  Notons $H=\Gal(K'/K)$. Formons le diagramme commutatif à carré cartésien
  \[\xymatrix{(X_{\eta'},G) \ar[r]\ar@/^1pc/[rr]^{(\pr_1,\Id)} & (X_{\eta'},G\times H) \ar[r]^{\pr_1} \ar[d]_{\pr_2} & (X,G)\ar[d]\\
  & (\eta',H) \ar[r] & (\eta,\{1\})}
  \]
  D'après \ref{prop.quot2} (b), l'homomorphisme $(\pr_1)_\yb : (G\times H)_d(y,\yb) \to
  G_d(x,\xb)$ est un isomorphisme. Soit $(g,h,\psi)$ l'image inverse
  de $(g,\phi)$. Alors $\psi|K'= \Id_{K'}$. En appliquant
  $(\pr_2)_\yb$, on obtient $(h,\psi) \in H_d(\eta',\etab)$, donc
  $h=1$, \ie $(g,\phi) = \Id_\yb (g,\psi)$.
\end{proof}

\subsection{$E$-\compte}\label{ss.compte}

Soient $p$ un nombre premier, $E$ un corps de caractéristique $0$,
$I$ un ensemble, $\gamma$ une application
\[I\to \ensdr{(l,\iota)}{l\text{ un nombre premier $\neq p$,
$\iota: E\hra \Qlb$ un plongement de corps}},
\]
où, pour chaque $l$, $\Qlb$ est une clôture algébrique de
$\Q_l$. Pour $\lbd \in I$, on écrit $\gamma(\lbd) = (l_\lbd,
\iota_\lbd)$.

Soient $K= \F_{p^f}$ un corps fini (\resp un corps local de corps
résiduel $k= \F_{p^f}$, \ie le corps des fractions d'un anneau de
valuation discrète hensélien excellent $R$ de corps résiduel
$k$), $\eta = \Spec K$. Soit $(X,G)\in \Eq/\eta$. On désigne par
$|X|$ l'\ens des points fermés de $X$. Pour $x\in |X|$, $\kappa(x)$
est une extension finie de~$K$. Soit $\xb$ un point \geoq
algébrique au-dessus de $x$. On note $F_0\in \Aut(\kappa(\xb))$ le
Frobenius \geoq (absolu) qui envoie $a$ sur $a^{1/p}$ (\resp $R_x$
l'anneau des entiers de $\kappa(x)$, $R_{\xb}$ le normalisé de
$R_x$ dans $\kappa(\xb)$, $x_0$ le point fermé de $\Spec R_x$,
$\xb_0$ le point fermé de $\Spec R_\xb$, $F_0\in
\Aut(\kappa(\xb_0))$ le Frobenius \geoq). On note $\rho$ l'inclusion
$\Gal(\kappa(\xb)/\kappa(x)^{G_d(x)}) \hra \Aut(\kappa(\xb))$ (\resp
l'homomorphisme composé $\Gal(\kappa(\xb)/\kappa(x)^{G_d(x)}) \to
\Gal(\kappa(\xb_0)/\kappa(x_0)^{G_d(x)}) \hra \Aut(\kappa(\xb_0))$).

\begin{defn}\label{def.Ecomp}
(i) On dit qu'un système $(t_\lbd)_{\lbd\in I}\in \prod_{\lbd\in I}
\lbar{\Q_{l_\lbd}}$ est $(E,I,\gamma)$-\emph{compatible} (ou
$E$-compatible s'il n'y a pas de confusion à craindre) s'il existe
$c\in E$ \tq $t_\lbd = \iota_\lbd (c)$ pour tout $\lbd\in I$.

(ii) On dit qu'un système de complexes $(L_\lbd)_{\lbd\in I} \in
\prod_{\lbd\in I} \Dbc(X,G,\lbar{\Q_{l_\lbd}})$ est
$(E,I,\gamma)$-\emph{compatible} (ou \Ecomp s'il n'y a pas de
confusion à craindre) si pour tout $x\in |X|$, tout $\xb\to x$ et
tout $(g,\phi)\in G_d(x,\xb)$ \eqref{eq.Gdxxb} avec $\rho(\phi)$ une
puissance entière de $F_0^f$, le système
$(\Tr((g,\phi),(L_\lbd)_\xb))_{\lbd\in I}$ est
$(E,I,\gamma)$-compatible.
\end{defn}

Les systèmes \Ecomps sur $(X,G)$ forment une sous-\cat triangulée
de la catégorie produit $\prod_{\lbd\in I}
\Dbc(X,G,\lbar{\Q_{l_\lbd}})$, notée $\Dbc(X,G,E)$. Lorsque
$G=\{1\}$, on écrit simplement $\Dbc(X,E)$. La \Ecompte de
$(L_\lbd)_{\lbd\in I}$ ne dépend que de $([L_\lbd])_{\lbd\in I}\in
\prod_{\lbd\in I} K(X,G,\lbar{\Q_{l_\lbd}})$ (\ref{ss.Dbc}).

\begin{prop}\label{prop.pos}
Soient $(L_\lbd)_{\ldI} \in \prod_\ldI
\Dbc(X,G,\lbar{\Q_{l_\lbd}})$, $x\in |X|$, $\xb\to x$, $N\in \Z$.
Supposons que pour tout $(g,\phi)\in G_d(x,\xb)$ avec $\rho(\phi) =
F_0^{fn}$, $n \ge N$, le système
$(\Tr((g,\phi),(L_\lbd)_\xb))_{\lbd\in I}$ est \Ecomp. Alors
$((L_\lbd)_x)_\ldI$ est \Ecomp.
\end{prop}

\begin{proof}
\Ops $X=x$ et $\sharp I=2$. \Ops $[L_\lbd] = [\cF_\lbd] -
[\cG_\lbd]$ avec $\cF_\lbd, \cG_\lbd \in
\Mod_c(x,G,\lbar{\Q_{l_\lbd}})$ semi-simples, $\lbd \in I$. Il
suffit de montrer que le système
$(\Tr((g,\phi),(L_\lbd)_\xb))_{\lbd\in I}$ est \Ecomp pour tout
$(g,\phi)\in G_d(x,\xb)$ avec $\rho(\phi) = F_0^{f(N-1)}$. On pose
$Z=\Aut(\kappa(\xb))$ (\resp $Z=\Aut(\kappa(\lbar{x_0}))$). On note
$\tau$ l'homomorphisme composé
\[G_d(x,\xb) \sto{\pr_2} \Gal(\kappa(x)/\kappa(x)^{G_d(x)}) \sto{\rho} Z \sto{r} \Zch,\]
où $r$ est l'isomorphisme qui envoie $F_0$ sur $1$. Dans le cas
fini, on pose $W=\tau^{-1}(\Z)$ ; dans le cas local, il existe un
sous-groupe ouvert $I_1$ de $\tau^{-1}(0)$ agissant trivialement sur
les $(\cF_\lbd)_\xb$ et $(\cG_\lbd)_\xb$, $\ldI$ en vertu du
théorème de monodromie locale de Grothendieck, et on pose
$W=\tau^{-1}(\Z)/I_1$. Dans les deux cas, $\tau$ induit un
homomorphisme surjectif $\nu : W\to f'\Z$ de noyau fini avec $f\mid
f'$. Le groupe infini $W$ agit par conjugaison sur l'ensemble fini
$\nu^{-1}(f')$. Le centre $Z(W)$ de $W$, qui est le fixateur de
$\nu^{-1}(f')$, est d'indice fini. Soit $h\in \tau^{-1}(\Z)$ \tq
$\tau(h)\ge f$ et que l'image de $h$ dans $W$ appartienne à $Z(W)$.
Alors pour tout $(g,\phi)\in G_d(x,\xb)$ avec $\tau(g) = f(N-1)$,
$gh$ et $hg$ définissent le même automorphisme sur
$(\cF_\lbd)_\xb$. Il en est de même sur $(\cG_\lbd)_\xb$. On
conclut en appliquant \cite[8.1]{IllMisc}.
\end{proof}





Le théorème suivant est le résultat principal de cet article.

\begin{theo}\label{theo.6op}
Sur un corps fini ou local comme dans \ref{ss.compte}, les six
opérations préservent la \Ecompte. Plus précisément, pour
$(f,\al) : (X,G) \to (Y,H)$ dans $\Eq/\eta$, $(L_\lbd)_\ldI ,
(M_\lbd)_\ldI \in \Dbc(X,G,E)$, $(N_\lbd)_\ldI \in \Dbc(Y,H,E)$, on
a
\begin{gather*}
(L_\lbd \otimes M_\lbd)_\ldI, (R\cHom(L_\lbd,M_\lbd))_\ldI,
 ((f,\al)^*N_\lbd)_\ldI
 \in \Dbc(X,G,E),\\
(R(f,\al)_*L_\lbd)_\ldI \in \Dbc(Y,H,E),
\end{gather*}
et, lorsque $f$ est séparé,
\[(R(f,\al)^!N_\lbd)_\ldI
 \in \Dbc(X,G,E), \quad (R(f,\al)_!L_\lbd)_\ldI \in \Dbc(Y,H,E).\]
\end{theo}

\begin{coro}\label{coro.D}
La dualité préserve la \Ecompte. Plus précisément, pour
$(X,G)\in \Eq/S$, $(L_\lbd)_\ldI \in \Dbc(X,G,E)$, on a
$(D_{(X,G)}L_\lbd)_\ldI \in \Dbc(X,G,E)$.
\end{coro}

\begin{proof}
Pour montrer que $(D_X L_\lbd)_\ldI$ est \Ecomp, on s'intéresse aux
traces en un point $x\in |X|$. Soit $V$ un voisinage ouvert séparé
de $x$ dans~$X$. Quitte à remplacer $(X,G)$ par $(\bigcap_{g\in
G_d(x)} a_g(V), G_d(x))$, \ops $X$ séparé. Ce cas découle de
\ref{theo.6op}.
\end{proof}

En fait, \ref{theo.6op} et \ref{coro.D} seront démontrés en même
temps (voir \ref{ss.demo.fini} et \ref{ss.demo.sousAB}).


\begin{remq}\label{rq.th}
(i) On retrouve le théorème de Gabber \cite[Th.~2]{Fujiwara} pour
$G=\{1\}$ et $K$ fini.

(i') Dans le cas $K$ local, soient $\eta_1$ une extension finie de
$\eta$, $\lbar{\eta_1} \to \eta_1$, $X\to \eta_1$ un morphisme séparé de
type fini muni d'une action d'un groupe fini $G$. Si $(g,\phi)\in
G_d(\eta_1,\lbar{\eta_1})$ avec $\rho(\phi) = F_0^{fn}$, $n\in \Z$,
alors \ref{theo.6op} combiné avec \cite[(2.4.1) et
(2.4.4)]{Zheng-inte} implique que
\[ \sum_{i\in \Z} (-1)^i \Tr((g,\phi), H_c^i(X_{\lbar{\eta_1}}, \Q_l)) \in p^{fd\min\{0,n\} } \Z \]
est indépendant de~$l$, où $d= \dim X$. On retrouve ainsi
\cite[Th.~B]{Ochiai}, le cas $s= \Spec \lbar{\F_p}$ de
\cite[4.2]{Vidal} et le cas (très particulier) de
\cite[0.1]{Saito-weight} où $\Gamma$ est la classe du graphe de
$a_g$.

(ii) \ref{theo.6op} équivaut à l'énoncé pour $\sharp I = 2$.

(iii) On peut remplacer $E$ par une extension algébrique quelconque
$E'$. En effet, pour tout $\lbd\in I$, notons $J_\lbd$ l'ensemble
des plongements $\mu : E'\to \Qllb$ qui prolongent $\iota_\lbd$.
Soit $I'$ une partie de $\coprod_{\lbd\in I} J_\lbd$, qui contient
au moins un élément de chaque $J_\lbd$, $\ldI$, et au moins un
$J_\lbd$, $\ldI$, au cas où $I$ est non vide. On note $\gamma' : I'
\to \left\{ (l,\iota) \right\}$ l'application qui envoie $\mu \in
J_\lbd$ sur $(l_\lbd,\mu)$. Alors un système $(t_\lbd)_\ldI \in
\prod_\ldI \Qllb$ est $(E,I,\gamma)$-compatible \ssi
$(t'_\mu)_{\mu\in I'} \in \prod_{\mu\in I'} \lbar{\Q_{l_\mu}}$ est
$(E',I',\gamma')$-compatible, où $t'_\mu = t_\lbd$ pour $\mu \in
J_\lbd$. Le foncteur
\begin{align*}
M: \prod_\ldI \Dbc(X,G,\Qllb) &\to \prod_{\mu\in I'} \Dbc(X,G,\lbar{\Q_{l_\mu}})\\
(L_\lbd)_\ldI &\mapsto (L'_\mu)_{\mu\in I'}, \quad \text{où $L'_\mu
= L_\lbd$ pour $\mu\in J_\lbd$,}
\end{align*}
commute aux six opérations et à la dualité. Un système
$(L_\lbd)_\ldI \in \prod_\ldI\Dbc(X,G,\Qllb)$ est $E$-compatible
\ssi $M((L_\lbd)_\ldI)$ est $E'$-compatible.

(iv) Les énoncés pour $(f,\al)^*$ et $\ot$ sont triviaux.
L'énoncé pour $(f,\al)_*$ avec $f$ fini résulte de
\eqref{eq.ftraces}. L'énoncé pour $D$ dans le cas où $X$ est
régulier et où les $L_\lbd$ sont à faisceaux de cohomologie
lisses découle du théorème de pureté : pour $x\in |X|$, $(D_X
L_\lbd)_x \simeq (D_x(L_\lbd)_x)(d_x)[2d_x]$, où $d_x = \dim_x X$,
et évidemment $D_x$ préserve la \Ecompte.

(v) On verra en \ref{theo.RP} que $R\Psi$ préserve aussi la
\Ecompte.
\end{remq}

\section{Descente galoisienne et indépendance de~$l$ sur un corps fini}\label{s2}

\begin{ssect}
Une technique de Deligne-Lusztig \cite[démonstration de
3.3]{Deligne-Lusztig} permet de se débarrasser de l'action de
groupe dans certaines circonstances.

Soit $S$ un schéma noethérien. On désigne par $\Eq^\adm/S$ la
sous-catégorie pleine de $\Eq/S$ formée des couples $(X,G)$ où
$G$ agit sur $X$ de façon admissible (\ref{ss.scheq}).
\end{ssect}
\begin{lemm}
Soient $G$ un groupe fini, $f : S'\to S$ un $G$-torseur. Le foncteur
\begin{align*}
\eps_f : \bC_{S} &\simto (\Eq^\adm/S)_G/(S',G)\\
X &\mapsto (X\times_{S}S',G),
\end{align*}
où $G$ agit sur $X\times_{S}S'$ via son action sur $S'$, est une
équivalence de catégories.
\end{lemm}

\begin{proof}
C'est un cas particulier de la descente galoisienne (voir
\cite[VIII~7.6]{SGA1} et \cite[6.2.B]{BLR}). Un quasi-inverse de
$\eps_f$ est donné par $(Y,G) \mapsto Y/G$.
\end{proof}

Notons $e_{f,X} : (X\times_{S}S', G) \to (X,\{1\})$ la projection.

\begin{ssect}\label{ss.Cmg}
Soit
\[\left((S_m,\Z/m\Z)_{m\ge 1}, \left((f_{m,n},\al_{m,n}) : (S_{mn},\Z/{mn}\Z) \to (S_m,\Z/m\Z)\right)_{m,n\ge 1}\right)\]
un système projectif, où $S_1=S$, $\al_{m,n} : \Z/mn\Z \to \Z/m\Z$
désigne l'homomorphisme qui envoie $\bar{1}$ sur $\bar{1}$, \tq
$f_{m,n}$ soit un  $\Ker \al_{m,n}$-torseur, $m,n\ge 1$. Soit
$(X,G)\in \Eq^\adm/S$. Soient $m\ge 1$, $g\in G$. Notons $n_g$
l'ordre de $g$. On considère l'action diagonale de $\Z/n_g\Z$ sur
$X_{S_m} \times_{S_m} S_{n_g m}$, celle sur $X_{S_m}$ étant donnée
par l'homomorphisme $i_g : \Z/n_g\Z \to G$ qui envoie $\bar{1}$ sur
$g$ et celle sur $S_{n_g m}$ donnée par l'homomorphisme $\Z/n_g \Z
\to \Z/n_g m \Z$ qui envoie $\bar{1}$ sur $\bar{m}$. On définit un
foncteur
\begin{align*}
(\Eq^\adm/S)_G &\to \bC_{S_m}\\
(X,G) &\mapsto X^{(m,g)} = \eps_{f_{m,n_g}}^{-1}(X_{S_m}
\times_{S_m} S_{n_g m}, \Z/n_g\Z).
\end{align*}
Posons
\[d_{(X,G),m,g}=(\pr_1,i_g) : (X_{S_m} \times_{S_m}
S_{n_g m}, \Z/n_g\Z) \to (X,G).
\]
Rappelons
\[e_{f_{m,n_g},X^{(m,g)}} : (X_{S_m} \times_{S_m}
S_{n_g m}, \Z/n_g\Z) \to (X^{(m,g)},\{1\}). \]

Supposons que $S$ soit régulier de dimension $\le 1$. Posons
\[C_{(X,G),m,g} = e_{f_{m,n_g},X^{(m,g)}*} d_{(X,G),m,g}^* :
\Dbc(X,G,\Qlb) \to \Dbc(X^{(m,g)},\Qlb).\]

La définition de $C_{m,g}$ commute aux six opérations et à la
dualité. Plus précisément, on a
\begin{gather*}
C_{(X,G),m,g} R\cHom_{(X,G)}(-,-) \simeq
R\cHom_{X^{(m,g)}}(C_{(X,G),m,g}-, C_{(X,G),m,g}-),\\
C_{(X,G),m,g}(-\otimes-) \simeq C_{(X,G),m,g}- \otimes
C_{(X,G),m,g}-,
C_{(X,G),m,g}D_{(X,G)} \simeq D_{X^{(m,g)}}
C_{(X,G),m,g},
\end{gather*}
et, pour $(f,\Id): (X,G) \to (Y,G)$ dans $\Eq^\adm/S$,
\[C_{(X,G),m,g}(f,\Id)^* \simeq (f^{(m,g)})^* C_{(Y,G),m,g},
C_{(Y,G),m,g}R(f,\al)_* \simeq R(f^{(m,g)})_* C_{(X,G),m,g},
\]
et, lorsque $f$ est séparé,
\[C_{(X,G),m,g} R(f,\Id)^! \simeq R(f^{(m,g)})^! C_{(Y,G),m,g},
C_{(Y,G),m,g} R(f,\al)_! \simeq R(f^{(m,g)})_! C_{(X,G),m,g}.
\]
\end{ssect}

Reprenons les notations de \ref{ss.compte}. En particulier, $K=
\F_{p^f}$ (\resp $K$ est un corps local de corps résiduel
$\F_{p^f}$). Pour $m\ge 1$, soient $K_m = \F_{p^{fm}}$ (\resp $K_m$
une extension \nre de $K$ de corps résiduel $\F_{p^{fm}}$), $\eta_m
= \Spec K_m$. Pour $m \ge 1$, $\Z/m\Z$ agit sur $\eta_m$ via
l'isomorphisme de groupes
\beq\label{eq.Z/mZ}
\Z/m\Z \to \Gal(K_{m}/K)
\eeq
où l'image $\phi$ de $\bar{1}$ est donné par $\phi =
F_0^{f}|K_{m}$ (\resp $\bar{\phi}= F_0^f|\F_{p^{fm}}$, où
$\bar{\phi}$ est la réduction de~$\phi$). Appliquons \ref{ss.Cmg}
au système projectif $(\eta_m)_{m\ge 1}$.

\begin{prop}\label{prop.DL} 
Soient $(X,G)\in \Eq^\adm/\eta$, $(L_\lbd)_\ldI \in \prod_\ldI
\Dbc(X,G,\Qllb)$. Pour que $(L_\lbd)_\ldI$ soit \Ecomp, il faut et
il suffit que pour tout $m\ge 1$ et tout $g\in G$,
\[(C_{(X,G),m,g}
L_\lbd)_\ldI \in \prod_\ldI \Dbc(X^{(m,g)},\Qllb)
\]
soit \Ecomp.
\end{prop}

\begin{proof}
La nécessité est claire. Soit $(L_\lbd)_\ldI \in
\prod_{\ldI}\Dbc(X,G,\lbar{\Q_{l_\lbd}})$ \tq pour tout $m\ge 1$ et
tout $g\in G$, $(C_{m,g} L_\lbd)_\ldI$ soit $E$-compatible. Comme
$e_{f_{m,n_g}}^* e_{f_{m,n_g}*} \simeq \Id$ (\ref{prop.quot} (b)),
on a $(d_{m,g}^* L_\lbd)_\ldI \in \Dbc(X_{\eta_m} \times_{\eta_m}
\eta_{n_g m},\Z/n_g\Z, E)$. D'après \ref{prop.pos}, il suffit de
montrer que pour tout $x\in |X|$, tout $\xb\to x$ et tout
$(g,\phi)\in G_d(x,\xb)$ avec $\rho(\phi)=F_0^{fm}$, $m\ge 1$,
$(\Tr((g,\phi),(L_\lbd)_\xb))_\ldI$ est \Ecomp. Soient $y\in
|X_{\eta_m}|$ au-dessus de $x$, $\yb\to y$ au-dessus de $\xb \to x$.
D'après \ref{prop.DL1},
 on a $(g,\phi)\in
G_d(y,\yb)$, correspondant à $(\bar{1},\phi)\in
(\Z/n_g\Z)_d(y,\yb)$. On prend le diagramme commutatif à carré
cartésien
\[\xymatrix{(X_{\eta_m}\times_{\eta_m}\eta_{n_g m}, \Z/n_g\Z) \ar[r]^-{(\Id,\Delta)} &
(X_{\eta_m}\times_{\eta_m}\eta_{n_g m}, \Z/n_g\Z\times {\Z/n_g\Z})
\ar[r] \ar[d] &
(X_{\eta_m},\Z/n_g\Z)\ar[d]\\
& (\eta_{n_g m},\Z/n_g\Z) \ar[r] & (\eta_m,\{1\})}
\]
Soient $z\in |X_{\eta_m}\times_{\eta_m} \eta_{n_g m}|$ au-dessus de
$y$, $\zb\to z$ au-dessus de $\yb\to y$. D'après \ref{prop.quot2}
(b), l'homomorphisme
\[(\Z/n_g\Z\times \Z/n_g\Z)_d(z,\zb) \to (\Z/n_g \Z)_d(y,\yb)\]
est un isomorphisme. Soit $(\bar{1},b,\phi) \in (\Z/n_g\Z\times
\Z/n_g\Z)_d(z,\zb)$ l'image inverse de $(\bar{1},\phi)$. Alors
$(b,\phi)\in (\Z/n_g\Z)_d(\eta_{n_g m},\etab)$, donc $b=\bar{1}$.
Bref, $(d_{m,g})_\zb : (\Z/n_g \Z)_d(z,\zb) \to G_d(x,\xb)$ envoie
$(\bar{1},\phi)$ sur $(g,\phi)$. Donc les
\[\Tr((g,\phi),(L_\lbd)_\xb) = \Tr((\bar{1},\phi),(d_{m,g}^*
L_\lbd)_\zb)
\]
forment un système \Ecomp.
\end{proof}

\begin{escho}
Avec les notations de la démonstration, si on note $w$ l'image de
$z$ par $e_{f_{m,n_g}}$, on choisit $\wb \to w$ au-dessous de
$\zb\to z$ et on désigne encore par $\phi$ l'image de
$(\lbar{1},\phi)$ par
\[(e_{f_{m,n_g}})_\zb : (\Z/n_g\Z)_d (z,\zb) \to \{1\}_d(w,\wb) =
\Gal(\kappa(\wb)/\kappa(w)),
\]
alors
\[ \Tr((g,\phi),(L_\lbd)_\xb) = \Tr(\phi, (C_{(X,G),m,g} L_\lbd)_\wb). \]
\end{escho}

\begin{ssect}\label{ss.demo.fini}
\begin{proof}[Démonstration de \ref{theo.6op} dans le cas d'un corps fini]
On sait que $(f,\al)^*$ et $\otimes$ préservent la \Ecompte. Pour
montrer \ref{theo.6op}, il suffit d'établir la stabilité par $D$
et $R(f,\al)_*$. La stabilité par $R(f,\al)_!$, $R(f,\al)^!$ et
$R\cHom$ s'obtiendront par dualité. Comme $R(f,\al)_* = (\Id,\al)_*
R(f,\Id)_*$, \ops $\al=\Id$. D'après \ref{prop.DL} et le théorème
de Gabber, $D$ et $R(f,\Id)_*$ préservent la \Ecompte lorsque les
actions de groupe sont admissibles.

Soit $(L_\lbd)_\ldI\in \Dbc(X,G,E)$. Pour montrer que $(D_X
L_\lbd)_\ldI$ est \Ecomp, on s'intéresse aux traces en un point
$x\in |X|$. Soient $V$ un voisinage ouvert séparé de $x$ dans~$X$,
$U$ un voisinage ouvert affine de $x$ dans $\bigcap_{g\in G_d(x)}
a_g(V)$. Quitte à remplacer $(X,G)$ par $(\bigcap_{g\in G_d(x)}
a_g(U), G_d(x))$, \ops $X$ affine. Ce cas est déjà connu.

Soit $f : X\to Y$ un morphisme $G$-\eqrt. Pour montrer que $Rf_*$
préserve la \Ecompte, on fait une récurrence sur $d=\dim X$. Le
cas $d< 0$ est trivial. Pour $d\ge 0$, on prend un ouvert dense
affine $G$-stable $U$ de $X$ et on note $j: U\hra X$. Pour
$(L_\lbd)_\ldI\in \Dbc(X,G,E)$, comme $[L_\lbd]-[Rj_* j^* L_\lbd]$
est à support de dimension $\le d-1$, il suffit de montrer que
$(Rj_*j^*L_\lbd)_\ldI$ et $(R(fj)_* j^* L_\lbd)_\ldI$ sont \Ecomps.
Donc \ops $X$ affine. Comme le problème est local sur $Y$, un
argument similaire à celui dans l'alinéa précédent permet de
supposer en plus $Y$ affine. Ce cas est déjà connu.
\end{proof}
\end{ssect}

\begin{ssect}\label{ss.j!*}
Rappelons le théorème de Gabber de stabilité de la \Ecompte par
extension intermédiaire \cite[Th.~3]{Fujiwara} : pour $j: U \hra X$
une immersion ouverte de schémas séparés de type fini sur un
corps fini $k = \F_{p^f}$, $(L_\lbd)_\ldI \in \Dbc(U,E)$ un système
$E$-compatible de faisceaux pervers purs de poids $w\in \Z$,
$(j_{!*} L_\lbd)_\ldI$ est $E$-compatible.
\end{ssect}

Rappelons que $K(X,\Qlb)$ est un groupe abélien libre engendré par
les classes d'isomorphie de faisceaux pervers simples. Pour $a\in
\Z$, notons $K_m(X,\Qlb)$ (\resp $K_{\le a}(X,\Qlb)$, \resp $K_{\ge
a} (X,\Qlb)$, \resp $K_a(X,\Qlb)$) le groupe de Grothendieck de la
catégorie des faisceaux pervers mixtes (\resp de poids $\le a$,
\resp de poids $\ge a$, \resp purs de poids $a$).

\begin{prop}\label{prop.poids}
Les projections canoniques
\[\text{$p_{\le a} : K_m(X,\Qlb) \to
K_{\le a} (X,\Qlb)$ et $p_{\ge a} : K_m(X,\Qlb) \to K_{\ge a}
(X,\Qlb)$}
\]
préservent la \Ecompte. Plus précisément, pour $(L_\lbd)_\ldI \in
\prod_\ldI K_m(X,\Qlb)$ un système \Ecomp, les systèmes $(p_{\le
a} L_\lbd)_\ldI$ et $(p_{\ge a} L_\lbd)_\ldI$ sont \Ecomps.
\end{prop}

\begin{proof}
Il suffit de montrer que la projection canonique $p_w : K_m(X,\Qlb)
\to K_w(X,\Qlb)$ préserve la $E$-compatibilité pour tout $w\in
\Z$, car $p_{\le a} = \sum_{w\le a} p_w$, $p_{\ge a} = \sum_{w\ge a}
p_w$. On peut supposer $\sharp I = 2$. Soient $L_\lbd = \sum_{w\in
\Z} L_{\lbd,w}$, $L_{\lbd,w} = [M_{\lbd,w,0}] - [M_{\lbd,w,1}]$, où
$M_{\lbd,w,0}$ et $M_{\lbd,w,1}$ sont des faisceaux pervers purs de
poids $w$.

(a) \emph{Cas où $X$ est lisse sur $k$ et $\fH^e(M_{\lbd,w,\al})$
lisse sur $X$ pour $\ldI$, $w\in \Z$, $\al=0,1$, $e\in \Z$}. Alors pour tout
$x\in |X|$, le facteur $L$ local $L_x(L_{\lbd,w},t)$ peut être
extrait de $L_x(L_\lbd,t)$ comme la partie de poids $w-\dim_x X$.
Donc les $L_{\lbd,w}$ forment un système $E$-compatible.

(b) \emph{Cas général.} On fait une récurrence sur $d=\dim X$,
simultanément pour tout $w\in \Z$. Le cas $d<0$ est trivial. Pour
$d\ge 0$, \ops $X$ réduit. Prenons un ouvert dense $j: U\hra X$
lisse sur $k$ tel que $\fH^e(j^* M_{\lbd,w,\al})$ soit lisse sur
$U$, pour $\ldI$, $w\in \Z$, $\al=0,1$, $e\in \Z$. Notons $i: Y\to
X$ le fermé complémentaire. D'après (a), $(j^* L_{\lbd,w})_\ldI$
est $E$-compatible, donc $(j_{!*} j^*L_{\lbd,w})_\ldI$ l'est aussi
en vertu de la démonstration du théorème de Gabber rappelé en \ref{ss.j!*}. D'après
\cite[5.3.11]{BBD}, pour tout $w\in \Z$, il existe $L'_{\lbd,w} \in
K_w(Y,\Qlb)$ tel que
\beq\label{eq.j!*}
L_{\lbd,w} = j_{!*} j^* L_{\lbd,w} + i_* L'_{\lbd,w}.
\eeq
En sommant pour tout $w\in \Z$, on obtient
\[ L_\lbd = \sum_{w\in \Z} j_{!*} j^* L_{\lbd,w} + i_* \sum_{w\in \Z} L'_{\lbd,w}, \]
d'où la $E$-compatibilité de $(\sum_{w\in \Z} L'_{\lbd,w})_\ldI$.
Par conséquent, pour tout $w\in \Z$, $(L'_{\lbd,w})_\ldI$ est
$E$-compatible en vertu de l'hypothèse de récurrence. La
$E$-compatibilité de $(L_{\lbd,w})_\ldI$ résulte alors de
\eqref{eq.j!*}.
\end{proof}

Notons $\Dbm(X,\Qlb)$ la sous-catégorie triangulée de $X$ formée
des complexes mixtes. Pour $a\in \Z$, notons ${}^wD^{\le a}
(X,\Qlb)$ (\resp ${}^wD^{\ge a} (X,\Qlb)$) la sous-catégorie
triangulée de $\Dbc(X,\Qlb)$ formée des complexes mixtes $K$ tels
que pour tout $i\in \Z$, ${}^{p}\fH ^{i} K$ soit de poids $\le a$
(\resp $\ge a$). S.~Morel 
\cite[3.1.1]{Morel-pondere} a défini un foncteur exact troncature
par le poids
\[\text{$w_{\le a} :
\Dbm(X,\Qlb) \to {}^wD^{\le a}(X,\Qlb)$ (\resp $w_{\ge a} :
\Dbm(X,\Qlb) \to {}^wD^{\ge a}(X,\Qlb)$),}
\]
adjoint à gauche (\resp à droite) de l'inclusion
\[\text{${}^wD^{\le a}(X,\Qlb)\subset
\Dbm(X,\Qlb)$  (\resp ${}^wD^{\ge a}(X,\Qlb) \subset
\Dbm(X,\Qlb)$),}
\]
rendant commutatif le diagramme
\[\xymatrix{
\Ob(^wD^{\le a}(X,\Qlb))\ar[d] & \ar[l]_{\Ob(w_{\le a})}
\Ob(\Dbm(X,\Qlb))
\ar[r]^{\Ob(w_{\ge a})}\ar[d] & \Ob(^wD^{\ge a}(X,\Qlb))\ar[d]\\
K_{\le a}(X,\Qlb) & \ar[l]_{p_{\le a}} K_m(X,\Qlb)\ar[r]^{p_{\ge a}}
& K_{\ge a}(X,\Qlb)}
\]

Pour $j: U\hra X$ un ouvert et $L$ un faisceau pervers pur de poids
$a\in \Z$, la flèche canonique $w_{\ge a} j_{!} L \to j_{!*} L$
est un isomorphisme 
\cite[3.1.4]{Morel-pondere}. En vue de cela, la stabilité par la
troncature par le poids (\ref{prop.poids}) généralise le
théorème de Gabber (\ref{ss.j!*}).

\newcommand{\anciencasfini}{%
On suppose que $K$ est un corps fini en
\ref{ss.Rf!} et \ref{prop.varab}.

\begin{ssect}\label{ss.Rf!} Soit $(f,\Id) : (X,G) \to (Y,G)$ un
morphisme dans $\Eq/\eta$ avec $f$ \sep. Alors $R(f,\Id)_!$
préserve la \Ecompte.

En effet, d'après le théorème de changement de base, \ops $Y=y =
\Spec K_1$, $K_1$ une extension finie de $K$. Soit $\Fr_X : X\to X$
le Frobenius absolu. Pour $(g,F_0^n) \in G_d(y,\yb)$ avec $n\ge 1$,
$a_g^{(n)} = a_g \circ \Fr_{X}^n = \Fr_X^n \circ a_g$ est un
$y$-morphisme. On a $(a_g^{(n)})^{\ord(a_g)} = \Fr_X^{n \cdot
\ord(a_g)}$, donc $\Fix(a_g^{(n)}) \subset \Fix(\Fr_X^{n \cdot
\ord(a_g)})$. Pour $L\in \Dbc(X,G,\Qlb)$, on considère la
correspondance $u_{g,L}^{(n)}$ sous-jacente à l'endomorphisme
 $a_g^{(n)}$ donnée par
\[a_g^{(n)}L \simeq a_g^* (\Fr_X^n)^* L \sto{a_g^* (\Fr^n_{L/X})} a_g^* L \sto{b_g} L,\]
où $\Fr^n_{L/X} : (\Fr_X^n)^* L \to L$ est la correspondance de
Frobenius itérée \cite[XV 2.1]{SGA5}. Pour $x\in \Fix(a_g^{(n)})$,
$(g,F_0^n)$ appartient à $G_d(x,\xb)$ et le terme local naïf de
$(u_{g,L}^{(n)})_\yb$ en $\xb$ est $\Tr((g,F_0^n),L_\xb)$. D'après
la formule des traces \cite[2.3.2]{Varshavsky}, il existe $N \ge 1$
qui ne dépend que de l'action de $G$ sur $X$, \tq dès que $n\ge
N$, on a
\[\Tr((g,F_0^n),(Rf_!L)_\yb) = \sum_{z\in Z(\yb)} \Tr((g,F_0^n),L_z),\]
où $Z=\Fix(a_g^{(n)})$, $Z(\yb)$ est l'ensemble des sections du
$\yb$-\sch $Z\times_y \yb$. Donc $Rf_!$ préserve la \Ecompte.

\begin{remqe}
Si $X$ est propre sur $y$, la formule précédente vaut pour $N=1$.
En effet, on a
\[
\xymatrix{\Tr((g,F_0^n),(Ra_{X!}L)_\yb) \ar@{=}[r] \ar@{=}[d]
& \sum_{z\in |Z \times_\eta \yb|} \Tr((g,F_0^n),L_z) \ar@{=}[d]\\
d\Tr((g,F_0^n),(Rf_!L)_\yb) & d\sum_{z\in Z(\yb)}
\Tr((g,F_0^n),L_z)}
\]
où $a_X : X\to \eta$, $d=[K_1 :K]$.
\end{remqe}
\end{ssect}

\begin{prop}\label{prop.varab}
Soient $K_1$ un corps extension finie de $K$, $\eta_1 = \Spec K_1$,
$a_A : A \to \eta_1$ une \varab, $G$ un groupe fini agissant sur
$a_A$ par des $\eta$-automorphismes. Supposons que $G_d(0_A)=G$.
Alors la dualité $D_A$ préserve la \Ecompte.
\end{prop}

La démonstration est calquée sur celle de \cite[Th.~2]{Fujiwara}.

\begin{proof}
Soit $(M_\lbd)_\ldI$ un système \Ecomp sur $A$. Soient $w \in |A|$,
$(g,F_0^n)\in G_d(w,\bar{w})$ avec $n\ge 1$ un multiple de $r$. Il
suffit de montrer que $(\Tr((g,F_0^n),(DM_\lbd)_{\bar{w}}))_\ldI$
est $E$-compatible.

On a $(g,F_0^n)\in G_d(\eta_1,{\etab})$. Donc $a_g^{(n)} : A\to A$
est un $\eta_1$-morphisme. Comme $G_d(0_A) = G$, $a_g^{(n)}$ fixe
$0_A$, donc est un homomorphisme. Pour tout entier $m\ge 1$, on
considère l'analogue de l'isogénie de Lang \cite[Sommes trig.\
1.4--6]{SGA4d}
\begin{align*}
L_m : A&\to A\\
x&\mapsto (a_g^{(n)})^m x \cdot x^{-1}.
\end{align*}
C'est une isogénie de noyau $Z_m\subset \Fix(\Fr_A^{mn\cdot
\ord(a_g^m)})$. Elle définit un $Z_m$-torseur $\fL_m$ sur~$A$. Pour
$h\in \langle g \rangle$, le carré cartésien
\[\xymatrix{A\ar[r]^{a_h} \ar[d]^{L_m} & A\ar[d]^{L_m}\\
A\ar[r]^{a_h} & A}
\]
donne un isomorphisme $\fL_m \to a_h^* \fL_m$. Son inverse définit
une action de $\langle g \rangle$ sur le $Z_m$-torseur $\fL_m$ sur
$A$. Soit $z\in Z_m$. Sur $L_m^{-1}(Z_m)$, $a_g^{(n)}x = a_g^{(n)}x
\cdot x^{-1} \cdot x= z\cdot x$. Donc sur $\fL_{m,z} = L_m^{-1}(z)$,
$b_g$ est donné par $x\mapsto z^{-1}x$. On a un \diagcomm
\[\xymatrix{0\ar[r] & Z_m\ar[d]^{T_m} \ar[r] & A \ar[d]^{P_m}\ar[r]^{L_m}
& A\ar@{=}[d] \ar[r] & 0\\
0\ar[r]& Z_1 \ar[r] & A\ar[r]^{L_1} & A\ar[r] &0}
\]
où $P_m$, $T_m$ sont donnés par $x\mapsto \prod_{i=0}^{m-1}
(a_g^{(n)})^i x$. Ceci fournit un isomorphisme de $Z_1$-torseur $T_m
\fL_m \simto \fL_1$. Pour $z\in Z_m$, $g$ agit sur $\fL_{1,z}$ par
$x\mapsto T_m(z^{-1})x$.

On définit les fonctions $f_{m,\lbd}, g_{m,\lbd} : Z_1(\etab) \to
\Qllb$ par
\begin{align*}
f_{m,\lbd}(x) &= \sum_{\substack{y\in Z_m(\etab)\\ T_m(y)=x}}
\Tr((g^m, F_0^{mn}),(M_\lbd)_y),\\
g_{m,\lbd}(x) &= \sum_{\substack{y\in Z_m(\etab)\\ T_m(y)=-x}}
\Tr((g^m, F_0^{mn}),(DM_\lbd)_y).
\end{align*}
Comme $w\in Z_1$, il suffit de montrer que pour tout $x\in
Z_1(\etab)$, $(g_{1,\lbd}(x))_\ldI$ est $E$-compatible. Par
hypothèse, les $(f_{m,\lbd}(x))_\ldI$ sont $E$-compatible. Donc il
suffit de montrer que pour tout $\ldI$ et tout $x\in
A_{g,n}(\etab)$, il existe une suite récurrente linéaire
$(s_{a,\lbd}(x))_{a\in \Z}$ à valeurs dans $\Qllb$ \teq
$f_{m,\lbd}(x)= s_{m,\lbd}(x)$, $g_{m,\lbd}(x)= s_{-m,\lbd}(x)$ pour
tout $m\ge 1$. Dans la suite, on fixe un $\lbd\in I$ et on l'enlève
de la notation.

Soit $K_2$ une extension finie de $K_1$ \teq $(Z_1)_{\eta_2}$ soit
constant, où $\eta_2 = \Spec K_2$. Pour un caractère $\rho :
Z_1(\etab) \to \Qlb^{\times}$, on associe le $\Qlb$-\faisc $\fL_\rho
= \rho^{-1} ((\fL_1)_{\eta_2})$ lisse de rang~$1$ sur $A_{\eta_2}$.
Pour une fonction $f : Z_1(\etab) \to \Qlb$, on pose $(\cF f)(\rho)=
\sum_{x\in Z_1(\etab)} f(x) \rho(x)$, la transformée de Fourier sur
$Z_1(\etab)$. On note $q : A_{\eta_2} \to \eta_2$. D'après la
formule des traces,
\begin{align*}
\Tr((g^m, F_0^{mn}), (Rq_! (M_{\eta_2}\ot \fL_\rho))_\etab) &=
\sum_{x\in Z_m(\etab)} \Tr((g^m,F_0^{mn}),M_x) \rho (T_m(x)) \\
&= (\cF
f_m)(\rho),\\
\Tr((g^m, F_0^{mn}), (Rq_! ((D_A M)_{\eta_2}\ot \fL_\rho))_\etab) &=
\sum_{x\in Z_m(\etab)} \Tr((g^m,F_0^{mn}),(D_A M)_x) \rho (-T_m(x))\\
&= (\cF g_m)(\rho).
\end{align*}
Puisque $Rq_! (M_{\eta_2}\ot \fL_\rho)$ et $Rq_! ((D_A
M)_{\eta_2}\ot \fL_{\rho^{-1}})$ sont duaux, il existe, pour tout
$\rho$, une suite récurrente linéaire $(S_a(\rho))_{a\in \Z}$ à
valeurs dans $\Qlb$ \teq $(\cF f_m)(\rho) = S_m(\rho)$, $(\cF
g_m)(\rho) = S_{-m}(\rho)$ pour tout $m\ge 1$. Il suffit alors de
prendre $s_a = \cF^{-1}S_a$.
\end{proof}
}

\section{Altérations galoisiennes et réduction au cas des
courbes}\label{s3}

On reprend les notations de \ref{ss.compte}. En particulier, $K$ est
un corps fini ou local.

Soit $X$ un \sch \tf sur $\eta=\Spec K$. On considère la condition
suivante portant sur $L\in \Dbc(X,\Qlb)$ :
\begin{numcond}\label{cond.mr}
Pour tout $i\in \Z$, il existe une extension finie $F$ de $\Q_l$ et
un $\cO$-faisceau $\cF_{\cO}$, \tsq $\fH^i L \simeq (\cF_\cO
\ot_{\cO} F) \ot_F \Qlb$ et que $\cF_\cO \ot_{\cO} (\cO/\fm)$ soit
constant sur toute composante connexe de $X$, où $\cO$ est l'anneau
des entiers de~$F$, $\fm$ est l'idéal maximal de~$\cO$.
\end{numcond}
Cette condition est stable par image inverse et $\ot$. Si $X\hra Y$
est un ouvert de complémentaire un \dcn, alors \ref{cond.mr}
implique que $L$ est à faisceaux de cohomologie lisses sur $X$,
modérément ramifiés sur $Y$.

Le résultat principal de ce \S~est la proposition suivante.

\begin{prop}\label{prop.sousAB}
On fait les hypothèses suivantes :

\tm{(A)} 
Soient $K_1$ une extension finie de $K$, $\eta_1= \Spec K_1$, $a_X :
X\to \eta_1$ une courbe lisse affine,
$(L_\lbd)_\ldI \in \Dbc(X,E)$ avec tous les $L_\lbd$ vérifiant
\ref{cond.mr}. Alors $(R a_{X!} L_\lbd)_\ldI$ est $E$-compatible.

\tm{(B)} Soient 
$X$ une courbe sur $\eta$, normale, 
$j: U\hra X$ un ouvert dense, $(L_\lbd)_\ldI\in \Dbc(U,E)$ avec tous
les $L_\lbd$ \ver ant \ref{cond.mr}. Alors $(Rj_* L_\lbd)_\ldI$ est
\Ecomp.

Alors la conclusion de \ref{theo.6op} est vraie.
\end{prop}

Dans le cas d'un corps local, (A) et (B) seront démontrés en
\ref{ss.demoAB}.

\begin{ssect}\label{ss.demoABfini}
Dans le cas d'un corps fini, 
(A) est une conséquence de la formule des traces de Grothendieck,
tandis que (B) est une variante de \cite[9.8]{constantes}. Prouvons
(B) dans une plus grande généralité, avec la condition \ref{cond.mr} remplacée par ce que les $L_\lambda$ sont à faisceaux de cohomologie lisses. \Ops $X$ irréductible. Quitte à remplacer $\eta$ par le corps de
définition de $X$, \ops $X$ géométriquement irréductible. \Ops
$X$ projective. D'après \ref{rq.th} (iii), \ops que $E$ contient
les racines $p$-ièmes de l'unité. D'après \ref{rq.th} (ii), \ops
$\sharp I = 2$. Notons $D=X-U$. Pour $x\in D$, posons $n =
\max_{\ldI, i\in \Z, y\in D-\{x\}} \Swan_y(\fH^i L_\lbd)$. D'après
\cite[Construction, p.~217]{Katz}, il existe
 un système  $E$-compatible $(\cF_\lbd)_\ldI$ de faisceaux lisses de
rang $1$ sur $X-(D-\{x\})$ \tq pour tout $y\in D-\{x\}$,
$\Swan_y(\cF_\lbd) > n$. Alors $Rj_* (L_\lbd \ot \cF_\lbd)_y = 0$,
pour $y\in D-\{x\}$ \ibid[Prop., p.~216], et $Rj_*(L_\lbd \ot
\cF_\lbd)_x \simeq (Rj_* L_\lbd)_x \ot (\cF_\lbd)_x$. Donc, d'après
la formule des traces, pour $m\ge 1$,
\begin{multline*}
\Tr(F_0^{rm}, (Ra_{U*} L_\lbd)_\etab) \\
= \sum_{y\in U(\F_{p^{rm}})}
\Tr(F_0^{rm},(L_\lbd)_\yb)\Tr(F_0^{rm},(\cF_\lbd)_\yb) + \frac{r}{f}
\Tr(F_0^{rm},(Rj_* L_\lbd)_\xb)\Tr(F_0^{rm},(\cF_\lbd)_\xb),
\end{multline*}
où $a_U: U\to \eta$, $r$ et $f$ sont \tsq $\kappa(x) = \F_{p^r}$,
$\kappa(\eta) = \F_{p^f}$. D'après \ref{rq.th} (iv) et la formule
des traces, $(Ra_{U*} L_\lbd)_\ldI \simeq (D_\eta Ra_{U!} D_U
L_\lbd)_\ldI$ est \Ecomp, donc $((Rj_* L_\lbd)_x)_\ldI$ est \Ecomp.
D'où (B).

Cette démonstration, jointe à \ref{prop.sousAB}, donne une
démonstration de \ref{theo.6op} dans le cas d'un corps fini
indépendante du théorème de Gabber.
\end{ssect}

\begin{ssect}
Si $K$ est un corps local, $E$ est un corps de nombres, pour
$\lbd_1$, $\lbd_2$ deux places finies de~$E$ ne divisant pas $p$,
$\rho_1$, $\rho_2$ des représentations $\lbd_1$- et
$\lbd_2$-adiques de $\Gal(\Kb/K)$, Deligne a défini une notion de
compatibilité entre $\rho_1$ et $\rho_2$ \cite[8.8]{constantes}.
Prenons pour $\cF_i$ le $E_{\lbd_i}$-faisceau sur $\eta=\Spec K$
correspondant à $\rho_i$, choisissons un plongement $E_{\lambda_i}
\hra \overline{\Q_{l_i}}$, $i\in I=\{1,2\}$, et posons
\begin{align*}
\gamma : I &\to \left\{(l,\iota)\right\}\\
i &\mapsto (l_i, E\hra E_{\lbd_i}\hra \lbar{\Q_{l_i}}).
\end{align*}
Alors $(\cF_i)_{i\in I}$ est $(E,I,\gamma)$-compatible au sens de
\ref{def.Ecomp} \ssi les semi-simplifiées de $\rho_1$ et $\rho_2$
sont compatibles au sens de Deligne.

En vue de cette interprétation, la proposition suivante
généralise \cite[9.8]{constantes}.
\end{ssect}

\begin{prop}\label{prop.saeta}
Soient $S$ une courbe lisse sur $k=\F_{p^f}$, $(X,G)\in \Eq/S$,
$(L_\lbd)_\ldI \in \Dbc(X,G,E)$. Soit $s$ un point fermé de $S$.
Notons $S_{(s)}$ le hensélisé de $S$ en $s$, $\eta_s$ son point
générique. Alors le système $(L_\lbd|X_{\eta_s})_\ldI$ est \Ecomp
(au sens de \ref{def.Ecomp} pour les schémas de type fini sur un
corps local).
\end{prop}

\begin{proof}
Quitte à remplacer $k$ par une extension finie, \ops $\kappa(s)=k$.
Le problème étant local sur $X$, \ops $X$ affine. Pour $m\ge 1$,
$g\in G$, on applique \ref{ss.Cmg} à $(X,G)$ sur $\F_{p^f}$ et à
$(X_{\eta_s},G)$ sur $\eta_s$. On obtient
\[\xymatrix{&X^{(m,g)} \ar[d] \\ s_m = s_{\F_{p^{fm}}}\ar[r] & S_{\F_{p^{fm}}}}\]
Notons $\eta_{s_m}$ le point générique du hensélisé de
$S_{\F_{p^{fm}}}$ en $s_m$. On a $(X^{(m,g)})_{\eta_{s_m}} \simeq
(X_{\eta_s})^{(m,g)}$ et $(C_{m,g}L_\lbd)|(X_{\eta_s})^{(m,g)}
\simeq C_{m,g}(L_\lbd|X_{\eta_s})$. Quitte à remplacer $X$ par
$X^{(m,g)}$ (\ref{prop.DL}), \ops $G=\{1\}$.

Soit $y$ un point fermé de $X_{\eta_s}$. Alors il existe un
morphisme étale de type fini $S'\to S$, une section $s\to S'$, et
un sous-schéma fermé $Y$ de $X_{S'-\{s\}}$, fini sur $S'-\{s\}$,
\tsq $Y_{\eta_s} \simeq y$, d'où le \diagcomm à carrés
cartésiens
\[\xymatrix{y\ar[r]\ar[d] &Y\ar[d]\\
X_{\eta_s}\ar[r]\ar[d] & X_{S'-\{s\}}\ar[r]\ar[d] & X\ar[d]\\
\eta_s \ar[r] & S'-\{s\}\ar[r] & S}
\]
Quitte à remplacer $S$ par $S'$, $X$ par $Y$, \ops que $X$ est
quasi-fini séparé sur~$S$. Prenons un $S$-plongement $j: X\to
\lbar{X}$ où $\lbar{X}$ est un $S$-\sch fini. Quitte à remplacer
$X$ par $\lbar{X}$, $L_\lbd$ par $j_! L_\lbd$, \ops que $X$ est fini
sur~$S$. Quitte à remplacer $X$ par une composante du normalisé de
$X^\red$, \ops de plus $X$ lisse sur $k$. Quitte à remplacer $S$
par $X$, \ops $X=S$.

Ce cas de la proposition est essentiellement \cite[9.8]{constantes}.
Rappelons comment le déduire du cas d'un corps fini de
\ref{theo.6op}. \Ops $\sharp I =2$. On pose $[L_{\lbd}] =
[\cF_\lbd] - [\cG_\lbd]$, où $\cF_\lbd$ et $\cG_\lbd$ sont des
faisceaux sur $S$, semi-simples sur $\eta_s$, $\ldI$. Soit $I_1$ un
sous-groupe ouvert de $I(\Kb/K)$ agissant trivialement sur les
$(\cF_\lbd)_{\lbar{\eta_s}}$ et $(\cG_\lbd)_{\lbar{\eta_s}}$, où
$K=\kappa(\eta_s)$, $\Kb$ est une clôture séparable de $K$. Soit
$\phi \in G_K= \Gal(\Kb/K)$ avec $\rho(\phi) = F_0^{fm}$, $m\ge 1$.
Notons $K_1$ l'extension finie séparable de $K$ correspondant au
sous-groupe fermé de $G_K$ topologiquement engendré par $I_1$
et~$\phi$. Alors $I(\Kb/K_1) = I_1$. Prenons $S_1$ une courbe lisse
sur $k$, $f: S_1 \to S$ un morphisme quasi-fini et $s_1$ un point de
$S_1$ au-dessus de $s$ \tsq le point générique $\eta_{s_1}$ du
hensélisé de $S_1$ en $s_1$ vérifie $\eta_{s_1} = \Spec K_1$.
Formons le carré cartésien
\[\xymatrix{f^{-1}(S-\{s\})\ar@{^{(}->}[r]^-{j}\ar[d]_{f'} & S_1\ar[d]^{f}\\
S-\{s\} \ar@{^{(}->}[r] & S}
\]
Alors ${f'}^* \cF_\lbd$ et ${f'}^* \cG_\lbd$ se prolongent en des
faisceaux $\cF_\lbd'$ et $\cG_\lbd'$ sur $S_1$, lisses dans un
voisinage de $s_1$. Posons $L_\lbd' = \cF'_\lbd \oplus \cG'_\lbd[1]$.
Alors $[(Rj_* {f'}^* L_\lbd)_{s_1}] \simeq [(L_\lbd')_{s_1}] -
[(L_\lbd')_{s_1}(-1)]$, d'où
\[\Tr(\phi, (L_\lbd)_{\lbar{\eta_s}}) = \Tr(F_0^{fm}, (L_\lbd')_{\lbar{s_1}})
= \frac{1}{1-p^{fm}} \Tr(F_0^{fm}, (Rj_* (f')^*
L_\lbd)_{\lbar{s_1}}),
\]
où $\lbar{\eta_s} \to \eta_s$, $\lbar{s_1} \to s_1$. Ces traces
sont \Ecomps, car $(Rj_* (f')^* L_\lbd)_\ldI$ l'est d'après le cas
d'un corps fini de \ref{theo.6op}. (En fait la démonstration de (B)
donnée en \ref{ss.demoABfini} suffit : on n'a pas besoin d'invoquer
\ref{prop.sousAB}.)
\end{proof}

\newcommand{\anciennedemoABfini}{%
Pour $K$ fini, on a démontré (A) en \ref{ss.Rf!}. On montre
maintenant (B) pour $K$ fini. On s'intéresse aux traces situées en
un point $x\in X-U$. Quitte à remplacer $K_1$ par une extension
fini $K_2$, $X$ par $X_{K_2}$ et $G$ par
$G\times_{\Gal(K_1/K_1^G)}\Gal(K_2/K_1^G)$ (\ref{prop.quot}), \ops
(a) que $x$ soit un point rationnel. Quitte à remplacer $G$ par
$G_d(x)$ et $X$ par un voisinage affine $G_d(x)$-stable de $x$, \ops
(a), (b) $G=G_d(x)$ et (c) que $X$ soit affine. Quitte à remplacer
$X$ par sa compactification lisse $G$-\eqrte (\ref{lemm.compact}),
\ops (a), (b) et (d) que $X$ soit propre. La courbe pointée $(X,x)$
sur $\eta_1$ se plonge dans sa jacobienne $J$, de façon
$G$-\eqrte. On désigne ce plongement par $i$. Alors $i_* Rj_*
L_\lbd = D_J(i_* j_! L_\lbd)$. Il suffit donc d'appliquer
\ref{prop.varab}.
}%



Les lemmes \ref{lemm.fact}, \ref{lemm.compact} et \ref{lemm.deJong}
seront utilisés dans la démonstration de \ref{prop.sousAB}.

\begin{lemm}\label{lemm.fact}
Soient $X$ un \sch noethérien intègre muni d'une action d'un
groupe fini~$G$, $Y\to X$ un morphisme étale de schémas. Alors il
existe un ouvert affine dense $G$-stable $U$ de~$X$, un morphisme
$(f,\al):(Z,G') \to (U, G)$ avec $\al$ surjectif faisant de $Z$ un
$(\Ker \al)$-torseur sur $U$, et une factorisation
\[
\xymatrix{Z\ar@{-->}[r]\ar[rd]^{f} & Y\times_X
U\ar@{^{(}->}[r]\ar[d]
&Y\ar[d]\\
&U \ar@{^{(}->}[r]& X}
\]
\end{lemm}

\begin{proof}
Quitte à remplacer $X$ par un ouvert affine dense $G$-stable de
$X$, \ops $X$ affine. Notons $\eta$ le point générique de $X$,
$G_0 = G_i(\eta)$. Alors $G_0$ agit trivialement sur $X$, l'action
de $G$ sur $X$ se factorise donc par $Q=G/G_0$. Comme $X\to X/G$ est
étale en~$\eta$, il existe un ouvert $U$ affine dense $G$-stable de
$X$ \tq $Y\times_X U \to U$ et $U\to U/G$ soient des revêtements
étales. \Ops que $G_0=\{1\}$. En effet, une fois ce cas établi, on
peut l'appliquer, dans le cas général, à l'action de $Q$ sur $U$.
On obtient $(f,\beta) : (Z,Q') \to (U,Q)$. Il suffit alors de
prendre $G'=G\times_Q Q'$, qui agit sur $Z$ à travers sa projection
sur $Q'$, et prendre pour $\al$ la projection $G'\to G$.

Soit $x$ un point \geoq de $U/G$. On a une équivalence de \cats
\begin{align*}
F: \{\text{revêtements étales de $U/G$}\} &\to
\{\text{$\pi_1$-ensembles finis à gauche}\}\\
T &\mapsto T_x,
\end{align*}
où $\pi_1 = \pi_1(U/G,x)$. Alors $F(U)\simeq\pi_1/A$, où $A$ est
un sous-groupe ouvert \dist \tq $\pi_1/A \simeq G^\op$, et
$F(Y\times_X U) \simeq \coprod_{i=1}^n \pi_1/B_i$, où les $B_i$,
$1\le i \le n$, sont des sous-groupes ouverts de $A$. Soit $C\subset
\bigcap_{i=1}^n B_i$ un sous-groupe ouvert \dist de $\pi_1$. On
prend $G'= (\pi_1/C)^\op \times (\Z/n\Z)$ et $Z$ un revêtement
étale de $U/G$ \tq $F(Z) = (\pi_1/C)\times(\Z/n\Z)$, l'action de
$G'$ sur $Z$ étant donnée par la translation à droite de $G'$ sur
$F(Z)$. On prend $g: Z\to Y\times_X U$ \tq $F(g)$ soit donné par
\begin{align*}
(\pi_1/C)\times(\Z/n\Z) &\to \coprod_{i=1}^n \pi_1/B_i\\
(\sigma, \lbar{i}) &\mapsto \lbar{\sigma}\in \pi_1/B_i.
\end{align*}
Notons $f$ le composé de $g$ et de la projection $Y\times_X U \to
U$, $\al$ le composé $G'\sto{\pr_1} (\pi/C)^\op \to (\pi/A)^\op
\simeq G$. Alors $F(f)$ s'identifie à $\al^\op$, donc $f$ est
$G'$-équivariant et fait de $Z$ un $(\Ker \al)$-torseur sur $U$.
\end{proof}

\begin{ssect}\label{ss.equiv}
Soient $u : X\to Z$ un morphisme de schémas muni d'une action d'un
groupe fini~$G$, $X\sto{v} Y \sto{w} Z$ une factorisation de $u$
dans la catégorie des schémas. Rappelons la construction d'une
factorisation équivariante \cite[7.6]{deJong1}. Pour $g\in G$,
notons $Y_g$ le $Z$-\sch défini par la composition $Y\sto{w} Z
\sto{a_g} Z$. Posons $G=\{g_1,\cdots,g_n\}$,
\[(Y/Z)^G = Y_{g_1} \times_Z \cdots \times_Z Y_{g_n}\]
sur lequel $G$ opère par $(y_{g_1},\cdots,y_{g_n})g =
(y_{gg_1},\cdots,y_{gg_n})$. Alors $f$ a la factorisation
$G$-équivariante $X \sto{v'} (Y/Z)^G \to Z$, où $v'$ envoie $x$
sur $(v(xg_1),\cdots, v(xg_n))$.
\end{ssect}

\begin{lemm}\label{lemm.compact}
Soient $f: X\to Z$ un morphisme séparé de type fini de schémas
noethériens, $G$ un groupe fini agissant sur $f$. Alors on peut
factoriser $f$ en $X\stackrel{j}{\hra} Y \sto{g} Z$ de façon
$G$-\eqrte, où $g$ est propre et $j$ est une immersion.
\end{lemm}

\begin{proof}
D'après le théorème de compactification de Nagata
\cite[4.1]{Conrad}, $f$ se factorise en $X\stackrel{j_0}{\hra} Y_0
\sto{g_0} Z$, où $g_0$ est propre et $j_0$ est une immersion. Il
suffit alors d'appliquer \ref{ss.equiv}.
\end{proof}

\newcommand{\anciencompact}{
Le lemme suivant est décalqué de \cite[4.4 (c)]{Vidal}.

\begin{lemm}\label{lemm.compact}
Soient $S$ un \sch noethérien \univjap, $f:V\to Z$ un morphisme de
$S$-\schs \tf, $G$ un groupe fini qui agit sur $V$ et $Z$ par des
$S$-automorphismes de telle façon que $f$ soit $G$-\eqrt.
Supposons que $V$ soit normal, que $V\to S$ soit \sep et que
l'action de $G$ sur $V$ soit admissible \cite[1.7]{SGA1}.

\tm{(i)} Alors on peut factoriser $f$ en $V\stackrel{j}{\hra} Y \sto{g}
Z$ de façon $G$-\eqrte, où $g$ est propre, $j$ est une immersion
ouverte dominante, $Y$ est normal, et $g$ est $G$-admissible (\ie
pour tout changement de base $Z' \to Z$ $G$-\eqrt avec $Z'$
$G$-admissible, $Y\times_Z {Z'}$ est $G$-admissible).

\tm{(ii)} De plus, si $f$ est quasi-fini, on peut s'arranger pour que $g$
soit fini.
\end{lemm}

\begin{proof}
(i) Traitons d'abord le cas spécial $Z=S$. \Ops que $G$ agit
transitivement sur $\pi_0(V)$. Soient $V_1 \in \pi_0(V)$, $\al:
G_d(V_1) \to G$. Si on sait factoriser $f|V_1$ en $V_1
\stackrel{j_1}{\hra} Y_1 \sto{g_1} S$ de façon
$G_d(V_1)$-équivariante avec les propriétés prescrites, alors $f$
se factorise en $V\simeq \al_* V_1 \stackrel{\al_*j}{\hra} \al_* Y_1
\to S$. Donc \ops $V$ \irr. Alors $U=V/G$ est un $S$-\sch \sep \tf
\irr, donc est isomorphe à un sous-\sch ouvert d'un $S$-\sch propre
\irr $X$. On prend $p : Y \to X$, où $Y$ est la normalisation de
$X_\red$ dans $R(V)$. Alors $p$ est fini et $p^{-1}(U) \simeq V$,
d'où un diagramme commutatif $G$-\eqrt
\[\xymatrix{V \ar@{^{(}->}[r]^{j}\ar[d] & Y\ar[d]^{p} \\
 U\ar@{^{(}->}[r]\ar[d] &  X \ar[ld]\\
S}
\]
Comme $G$ agit trivialement sur $X$, il agit sur $Y$ de façon
admissible.

Dans le cas général, on applique ce qui précède et obtient une
factorisation $V\hra P \sto{p} S$, avec $p$ propre et
$G$-admissible. Donc on peut factoriser $f$ en $V\hra P\times_S Z
\sto{p_Z} Z$. Il suffit alors de prendre pour $Y$ le normalisé de
l'adhérence schématique de $V$ dans $P\times_S Z$.

(ii) Si $f$ est quasi-fini, on applique la factorisation de Stein à
$g$ et obtient $Y\to Y' \sto{h} Z$, où $h$  est fini et $V\to Y'$
est une immersion ouverte dominante. Cette factorisation est
$G$-\eqrte. Il suffit donc de remplacer $Y$ par le normalisé de
$Y'_\red$.
\end{proof}

\begin{lemme}\label{lemm.compact2}
Soient $S$ un \sch \noeth, $f: X\to Z$ un morphisme de $S$-\schs \tf
avec $X$ réduit et $X\to S$ \sep, $G$  un groupe fini agissant sur
la situation. Supposons vérifiée l'une des deux conditions
suivantes
\begin{cond}
(a) $S$ est universellement japonais et $G$ agit trivialement
là-dessus.

(b) $S$ est affine.
\end{cond}
Alors il existe un \diagcomm $G$-\eqrt
\[\xymatrix{X'\ar@{^{(}->}[r]^j\ar[d]^{\psi} &Y\ar[d]^{g}\\
X\ar[r]^{f} & Z}
\]
où $j$ est une immersion ouverte, $g$ est propre, et $\psi$ est
projectif et birationnel. De plus, on peut s'arranger pour que $Y$
soit réduit et $j$ soit dominant.
\end{lemme}

\begin{proof}
Dans le cas (a), on choisit un système de représentants $C$ des
orbites de $G$ agissant sur l'ensemble des composantes \irrs de $X$.
Pour $W\in C$, on note $G_d(W) = \ensdr{g\in G}{a_g(W) = W}$, $i_W :
G_d(W) \to G$. Alors $w:\coprod_{W\in C} \al_{W*}(W) \to X$ est
finie et birationnel. Donc il suffit de montrer le lemme pour
$f\circ w$. Pour le faire, il suffit de montrer le lemme pour chaque
$\al_{W*}(W) \to Z$ et de prendre la somme directe. Pour cela, il
suffit de montrer le lemme pour $G_d(W)$ agissant sur $W\to Z$ et
d'appliquer $\al_{W*}$. Donc \ops $X$ intègre. En appliquant le
lemme de Chow, on peut supposer en plus $X$ quasi-projectif sur $S$
(voir \cite[7.6]{deJong1}). Quitte à remplacer $X$ par son
normalisé, \ops $X$ normal, quasi-projectif sur $S$. Il suffit
alors d'appliquer \ref{lemm.compact}.

Dans le cas (b), le passage à la limite fournit un \diagcomm
$G$-\eqrt à carrés cartésiens
\[
\xymatrix{X\ar[d]^f \ar[r]^h & X_0\ar[d]^{f_0}\\
Z\ar[r]\ar[d] & Z_0\ar[d]\\
S\ar[r] & S_0}
\]
avec $S_0$ un \sch affine \tf sur $\Spec \Z$, $X_0$, $Z_0$ \tf sur
$S_0$ et $X_0$ \sep et réduit. Quitte à remplacer $X_0$ par
l'image schématique de $h$, \ops $h$ dominant. Il suffit alors
d'appliquer le cas (a) à $f_0$ sur $\Spec \Z$.
\end{proof}
}

Le lemme clef suivant découle de raffinements, dus à Gabber, de
résultats de de Jong (\cite{deJong2}, \cite[4.4]{Vidal}) sur les
altérations équivariantes.

On appelle \emph{altération galoisienne} \cite[4.4.1]{Vidal} un
morphisme $(f,\al):(S',G') \to (S,G)$ où $\al$ est surjectif, $f$
est un morphisme propre dominant \gent fini de \schs \noeths
intègres, \tq $K(S')^{H}$ soit une extension radicielle de $K(S)$,
où $H=\Ker \al$, $K(S')$ (\resp $K(S)$) est le corps des fonctions
de $S'$ (\resp $S$). Alors
il existe un ouvert affine dense $G$-stable $U$ de $S$ \tq
$f^{-1}(U) \to f^{-1}(U)/{H}$ soit un revêtement étale galoisien
de groupe $H/H_0$, où $H_0=H \cap \Ker (G'\to \Aut(K(S')))$, et que
$f^{-1}(U)/{H} \to U$ soit un homéomorphisme universel fini et
plat. La composition de deux altérations galoisiennes est une
altération galoisienne.

Soit $G$ un groupe fini. Soit $X$ un schéma noethérien régulier muni d'une
action de $G$. Rappelons qu'un \dcn $G$-stable $D$ de $X$ est dit
$G$-\emph{strict} s'il est strict et si pour toute composante
irréductible $D_i$ de $D$ et tout $g\in G$ \tq $D_i g \neq D_i$,
$D_i$ et $D_i g$ ne se coupent pas. Soit $S=\Spec R$ un trait de point fermé
$s$. Rappelons qu'un \emph{couple semi-stable} \cite[5.5
(a)]{Zheng-inte} est un couple (X,Z), où $X$ est un $S$-\sch \tf et
$Z$ est une partie fermée de $X$ contenant $X_s$, qui est,
localement pour la topologie étale, de la forme
\[(\Spec R[t_1,\dots,t_n]/(t_1\cdots t_a -\pi), Z), \]
où $\pi$ est une uniformisante de $R$, $Z$ est défini par l'idéal
$(t_1\cdots t_b)$, $1\le a \le b \le n$. Alors $Z$ est un \dcn de
$X$. Si $X\to S$ est muni d'une action de~$G$ sous laquelle $Z$ est
stable, on dit que le couple semi-stable $(X,Z)$ est
\emph{$G$-strict} si $Z$ est un \dcn $G$-strict de~$X$.

Soit $S$ un \sch noethérien. On dit qu'une courbe nodale $G$-équivariante $X\to S$ est
\emph{$G$-scindée} si elle est scindée et si pour tout $s\in S$,
toute composante irréductible $C$ de $X_s$ et tout $g\in G_d(s)$
\tq $Cg \neq C$, $C$ et $Cg$ ne se coupent pas (condition $(*)$ dans
\cite[4.4.1]{Vidal}). On appelle \emph{fibration plurinodale
$G$-scindée} (\cf \cite[5.8]{deJong2}) un système $(X_d\sto{f_d}
\cdots\sto{f_1} X_0, \{\sigma_{ij}\}, Z_0)$, où
\begin{cond}
--- $f_i:X_i \to X_{i-1}$ est une courbe projective, nodale et
$G$-scindée,

--- $\sigma_{ij}:X_{i-1}\to X_i$, $j=1,\cdots,n_i$, sont des sections
disjointes de $f_i$ dans le lieu lisse de $f_i$, permutées par $G$
(\ie pour tout $j$, il existe $j'$ \tq $a_g\circ \sigma_{ij} =
\sigma_{ij'}\circ a_g$),

--- $Z_0\subsetneq X_0$ est un fermé $G$-stable,
\end{cond}
\tq $f_i$ soit lisse sur $X_{i-1}-Z_{i-1}$, où on a posé $Z_i =
\bigcup_{j=1}^{n_i} \sigma_{ij}(X_{i-1}) \cup f_i^{-1}(Z_{i-1})$,
$1\le i \le d$. Ici on dit qu'un morphisme $f: X\to Y$
est \emph{projectif} si $f$ admet une factorisation
$X\sto{i} \bP^n_Y \to Y$, où $i$ est une immersion fermée. Dans la
définition d'une courbe plurinodale $G$-scindée, $X_d \to X_0$ est
projectif et plat.

\begin{lemm}\label{lemm.deJong}
Soient $G$ un groupe fini, $f: X \to S$ un morphisme \sep \tf
$G$-\eqrt de \schs \noeths intègres, $U$ un ouvert dense $G$-stable
de $X$. Supposons que la fibre générique $X_\eta$ de $f$ soit
\geoqt \irr. On fait aussi les hypothèses suivantes :
\begin{cond}
\tm{(i)} $S$ est excellent.

\tm{(ii)} pour toute altération galoisienne $(S_1,G_1) \to (S,G)$ et
toute partie fermée $F_1\subsetneq S_1$, il existe une altération
galoisienne $(S_2,G_2)\to (S_1,G_1)$ \teq $S_2$ soit régulier et
que l'image inverse de $F_1$ dans $S_2$ soit contenue dans un \dcn
$G_2$-strict.
\end{cond}
Alors il existe un groupe fini $G'$, un homomorphisme de groupes
$\al : G'\to G$, un diagramme commutatif $G'$-\eqrt
\[\xymatrix{X'\ar[d]^{v}\ar@{^{(}->}[r]^{j} & X_d\ar[r]^{g} & S'\ar[d]^{u}\\
 X \ar[rr]^{f} & & S}
\]
où $(u,\al)$, $(v,\al)$ sont des \altgals, $j$ est une immersion
ouverte dominante, et une fibration plurinodale $G'$-scindée de $g$
\[(X_d\sto{f_d} \cdots\sto{f_1} X_0=S',
\{\sigma_{ij}\}, Z_0)
\]
\teq $X_i$ soit régulier, $Z_i$ soit un \dcn $G'$-strict de $X_i$
pour $0\le i\le d$ et que $Z_d$ contienne $X_d - v^{-1}(U)$.
\end{lemm}

On appliquera le lemme uniquement quand $S$ est le spectre d'un
corps (\ref{ss.demo.sousAB}) ou un trait excellent (\ref{ss.demoAB},
\ref{ss.demoRP}). Dans ces cas, les hypothèses (i) et (ii) sont
automatiquement satisfaites. Lorsque $S$ est un trait excellent, $u$
est nécessairement un changement de traits fini, et $(X_d,Z_d)$ est
un couple semi-stable $G'$-strict sur $S'$.

\begin{proof}
D'après \ref{lemm.compact}, \ops $f$ propre. Si $\dim X_\eta=0$,
alors $(f,\Id_G)$ est une altération galoisienne, donc il suffit
d'y appliquer l'hypothèse (ii) et de prendre $g=\Id$. Si $\dim
X_\eta \ge 1$, on applique \cite[4.4.3]{Vidal}.
On obtient un groupe fini~$G'$, un
homomorphisme de groupes $\al : G'\to G$, un \diag $G'$-\eqrt
\[
\xymatrix{X_d\ar[r]^{g}\ar[d]_v & S'\ar[d]^{u}\\
 X \ar[r]^{f} & S}
\]
où $(u,\al)$, $(v,\al)$ sont des \altgals, et une fibration
plurinodale $G'$-scindée de $g$
\[(X_d\sto{f_d} \cdots\sto{f_1} X_0=S',
\{\sigma_{ij}\}, Z_0)
\]
\teq $Z_d$ contienne $X_d - v^{-1}(U)$. En appliquant (ii), on se
ramène au cas où $X_0$ est régulier et $Z_0$ est un \dcn
$G'$-strict. Appliquant alors \cite[4.4.4]{Vidal} successivement à
$f_1,\cdots, f_d$, on modifie $X_i$ de sorte que $X_i$ devient
régulier et $Z_i$ devient un \dcn $G'$-strict, $1\le i \le d$.
\end{proof}

\begin{ssect}\label{ss.demo.sousAB}
\begin{proof}[Démonstration de \ref{prop.sousAB}]
On sait que $f^*$ et $\ot$ préservent la \Ecompte. Prouvons d'abord
l'énoncé suivant :
\begin{cond}
$(A')$ Pour $(f,\al) :(X,G) \to (Y,H)$ un morphisme de $\Eq/\eta$
avec $f$ séparé, $R(f,\al)_!$ préserve la \Ecompte.
\end{cond}
Comme $R(f,\al)_! = (\Id,\al)_! R(f,\Id)_!$ et $(\Id,\al)_!$
préserve la $E$-compatibilité (\ref{rq.th} (iv)), \ops $\al=\Id$.
D'après le théorème de changement de base, \ops que $Y=\eta_1$
est le spectre d'une extension finie de $K$. \Ops $\sharp I = 2$.

(a) \emph{Cas $X=\bA^1_{\eta_1}$.} Choisissons $w: W\to X$ étale
dominant \tq les $w^*L_\lbd$ \ver ent \ref{cond.mr}. On applique
\ref{lemm.fact} à $w$ pour trouver un ouvert affine dense
$G$-stable $j : U\hra X$ et un morphisme $(f,\al) :(U',G')\to (U,G)$
avec $\al$ surjectif faisant de $U'$ un $(\Ker\al)$-torseur sur~$U$
\tsq les $(f,\al)^*j^*L_\lbd$ \ver ent \ref{cond.mr}. Soit $i: Z=X-U
\to X$. On a le \trdist
\[Ra_{U!} j^* L_\lbd \to Ra_{X!}L_\lbd \to a_{Z!} i^* L_\lbd \to.\]
Ici les $Ra_{U!} j^* L_\lbd \simeq R(a_{U'},\al)_! (f,\al)^* j^*
L_\lbd$ (\ref{prop.quot} (a), \ref{ss.omi}) forment un système
\Ecomp en vertu de (A) et de \ref{prop.DL}, $(a_{Z!} i^*
L_\lbd)_\ldI$ est \Ecomp car $a_Z$ est fini. Donc
$(Ra_{X!}L_\lbd)_\ldI$ l'est aussi.

(b) \emph{Cas général.} On fait une récurrence sur $d=\dim X$. Le
cas $d\le 0$ est trivial. Si $d\ge 1$, on prend un ouvert affine
dense $G$-stable $U$. Soit $i: Z=X-U \to X$. On a le \trdist
\[Ra_{U!} j^* L_\lbd \to Ra_{X!}L_\lbd \to Ra_{Z!} i^* L_\lbd \to.\]
D'après l'\hyp de \rec, $(Ra_{Z!} i^* L_\lbd)_\ldI$ est \Ecomp.
Donc il suffit de voir que $(Ra_{U!}j^* L_\lbd)_\ldI$ est \Ecomp.
D'après le lemme de normalisation, il existe un morphisme $U/G \to
\bA^1_{\eta}$ à fibres de dimension $\le d-1$. Ceci induit un
$\eta_1$-morphisme $G$-\eqrt $f : U\to \bA^1_{\eta_1}$. On a
$Ra_{U!} \simeq Ra_{\bA^1_{\eta_1}!} Rf_!$. Il suffit donc
d'appliquer l'\hyp de \rec et (a).

Pour prouver \ref{theo.6op}, il suffit de prouver les énoncés
suivants pour tout $d\in \N$ :
\begin{cond}
\textit{$(C_d)$ Pour $(X,G)$ sur $\eta$ avec $\dim X \le d$, $D_X$
préserve la \Ecompte.}

\textit{$(D_d)$ Pour $(f,\al): (X,G) \to (Y,H)$ un morphisme
 de $\Eq/\eta$ (où $f$ est éventuellement non séparé) avec $\dim X \le d$,
$R(f,\al)_*$ la préserve aussi.}
\end{cond}
\Ops $\sharp I =2$. On fait une récurrence sur $d$. 
Le cas $d = 0$ est trivial (\ref{rq.th} (iv)). Pour $d\ge 1$,
supposons $(C_{d-1})$ et $(D_{d-1})$ établis.

Prouvons $(C_d)$. \Ops $X$ réduit. Soit $(L_\lbd)_\ldI$ un système
\Ecomp sur~$X$. On prend un ouvert dense régulier $G$-stable $j:
U\hra X$ \tq les $L_\lbd |U$ soient à \faisx de cohomologie lisses.
On pose $i : Z=X-U \to X$. Le \trdist
\[i_*Ri^! D_X L_\lbd \to D_X L_\lbd \to Rj_* j^* D_X L_\lbd \to\]
se récrit
\[ i_* D_Z(i^* L_\lbd) \to D_X L_\lbd \to Rj_* D_U j^* L_\lbd \to.\]
D'après $(C_{d-1})$, $(D_{Z}(i^*L_\lbd))_\ldI$ est \Ecomp. D'après
\ref{rq.th} (iv), $(D_U j^* L_\lbd)_\ldI$ est \Ecomp. Donc la preuve
de $(C_d)$ se ramène à celle de $(D_d)$.

Prouvons $(D_d)$. On prend $(L_\lbd)_\ldI \in \Dbc(X,G,E)$.

\begin{sssect}\label{sss.modif}
Supposons qu'il existe $(L'_\lbd)_\ldI\in\Dbc(X,G,E)$ \tq pour tout
$\ldI$, $[L_\lbd]-[ L'_\lbd]$ soit à support de dimension $\le
d-1$, \ie il existe $i: Z_\lbd \to X$ sous-\sch fermé $G$-stable de
dimension $\le d-1$ et $L''_\lbd\in K(Z_\lbd,G,\Qllb)$ \tsq
$[L_\lbd] - [L'_\lbd] = i_* L''_\lbd$. En appliquant $(D_{d-1})$, on
voit qu'il suffit de \ver er $(Rf_* L'_\lbd)_\ldI \in \Dbc(Y,G,E)$.
\end{sssect}

\Ops $X$ réduit. Soit $j:U \hra X$ un ouvert dense normal
$G$-stable. On considère $L_\lbd \to Rj_* j^* L_\lbd$.
D'après \ref{sss.modif}, il suffit de montrer que $(Rj_* j^*
L_\lbd)_\ldI$ et $(R(fj)_* j^* L_\lbd)_\ldI$ sont \Ecomps. Donc \ops
$X$ normal. \Ops de plus que $G$ agit \tranvt sur $\pi_0(X)$. Si $C$
est une composante connexe de $X$, alors $(C,G_d(C)) \to (X,G)$ est
cocartésien (\ref{sss.cocart}). Quitte à remplacer $X$ par une
composante connexe (\ref{sss.cocarteb}), \ops $X$ intègre.

Choisissons $w: W\to X$ étale dominant \tq les $w^*L_\lbd$ \ver ent
\ref{cond.mr}. On applique \ref{lemm.fact} à $w$ pour trouver un
ouvert dense $G$-stable $U$ de $X$, et $(g,\beta) : (Z,G')\to (U,G)$
avec $\beta$ surjectif faisant de $Z$ un $(\Ker \beta)$-torseur sur
$U$, \tsq les $L_\lbd|Z$ \ver ent \ref{cond.mr}. D'après
\ref{sss.modif}, il suffit de montrer que $(Rj_*j^* L_\lbd)_\ldI$ et
$(R(fj)_* j^* L_\lbd)_\ldI$ sont \Ecomps, où $j: U\hra X$. D'après
\ref{prop.quot} (a), $j^*L_\lbd \simeq (g,\beta)_* (g,\beta)^* j^*
L_\lbd$. Donc, quitte à remplacer $X$ par $Z$, \ops
\begin{cond}
(a) les $L_\lbd$ \ver ent \ref{cond.mr}.
\end{cond}
D'après \ref{lemm.compact}, $(f,\al)=(p,\al)(j,\Id)$, où $j$ est
une immersion ouverte et $p$ est propre. D'après $(A')$, $R(p,\al)_*$
préserve la \Ecompte. Donc \ops (a) et
\begin{cond}
(b) $f$ est une immersion ouverte et $\al=\Id$.
\end{cond}
Ces hypothèses seront toujours conservées.

Le cas $d=1$ résulte alors de (B). En effet, on s'intéresse aux
traces situées en un point $y\in Y-X$. Quitte à remplacer $G$ par
$G_d(y)$ et $Y$ par un voisinage affine $G_d(y)$-stable de~$y$, \ops
de plus $Y$ affine. En vertu de \ref{prop.DL}, \ops $G=\{1\}$. Soit
$g: Y'\to Y$ le normalisé de $Y_\red$. Formons le carré cartésien
\[\xymatrix{X'\ar[d]^{g_X}\ar@{^{(}->}[r]^{f'} & Y'\ar[d]^g\\
X\ar@{^{(}->}[r]^{f} & Y}
\]
Le morphisme $L_\lbd \to g_{X*} g_X^*L_\lbd$ est un \iso sur le lieu
de normalité de $X_\red$, qui est un ouvert dense. D'après
\ref{sss.modif}, il suffit donc de voir que les $Rf_* g_{X*} g_X^*
L_\lbd = g_* Rf'_* g_X^* L_\lbd$ forment un système \Ecomp. Quitte
à remplacer $f$ par $f'$, \ops de plus $Y$ normal. En appliquant
(B), on achève la démonstration de $(D_1)$.

Dans la suite, $d\ge 2$. Soit $K_1$ une extension galoisienne de $K$
\teq les composantes \irrs de $Y\times_\eta \eta_1$ soient \geoqt
\irrs \cite[4.5.10]{EGAIV}, où $\eta_1=\Spec K_1$. Quitte à faire
le changement de base $\eta_1\to \eta$ et remplacer $G$ par
$G\times\Gal(K_1/K)$, \ops qu'on a un morphisme $G$-\eqrt $Y\to
\eta_1$ et que les composantes \irrs de $Y$ sont \geoqt \irrs sur
$\eta_1$. Comme dans l'alinéa précédent, on peut supposer de plus
$Y$ normal. Quitte à remplacer $Y$ par une composante connexe, on
peut supposer $Y$ intègre, donc \geoqt \irr sur $\eta_1$. On peut
supposer $\dim Y= d$ et $X$ non vide.

On applique \ref{lemm.deJong} à $Y\to \eta$ et l'ouvert $X$ de $Y$.
On obtient une altération galoisienne $(v,\beta) : (Y',G')\to
(Y,G)$ avec $Y'$ \reg et un \dcn $G'$-strict $D$ de $Y'$ contenant
$Y'-v^{-1}(X)$. Soit $U=Y'-D$ et formons le \diag à carré
cartésien
\[
\xymatrix{(U,G') \ar@{^{(}->}[r]^{(j,\Id_{G'})} &(v^{-1}(X),
G')\ar@{^{(}->}[r]^(.6){(f',\Id_{G'})} \ar[d]^{(v_X,\beta)}
& (Y',G')\ar[d]^{(v,\beta)}\\
&(X,G)\ar@{^{(}->}[r]^(.6){(f,\Id_G)} & (Y,G)}
\]
Le morphisme $L_\lbd \to R(v_{X},\beta)_* v_X^*L_\lbd$ est un \iso
sur un ouvert dense de $X$ (\ref{prop.quot} (a)). D'après
\ref{sss.modif} et $(A')$, il suffit donc de voir que les $Rf_*
R(v_{X},\beta)_* v_X^* L_\lbd = R(v,\beta)_* Rf'_* v_X^* L_\lbd$
forment un système \Ecomp. Par $(A')$, il suffit alors de voir que
$(Rf'_* v_X^* L_\lbd)_\ldI$ est \Ecomp. Comme $\dim (v^{-1}(X)-U)\le
d-1$, il suffit, d'après \ref{sss.modif}, de voir que $(Rj_*
j^*v_X^* L_\lbd)_\ldI$ et $(R(f'j)_* j^*v_X^* L_\lbd)_\ldI$ sont
\Ecomps. Donc \ops
\begin{cond}
(c) $Y$ est \reg irréductible de dimension $d$.
\end{cond}
et $X$ est le complémentaire d'un \dcn $G$-strict de $Y$.

Posons $Y-X = \sum_{i\in J} D_i$. Pour $i\in J$, notons $X_{\ip} = Y
- \bigcup_{j\in J-\{i\}} D_j$, et formons le diagramme à carré
cartésien
\[\xymatrix{&D_i^\circ \ar@{^{(}->}[r]^{j'_\ip}\ar[d]_{\iota'_i} & D_i \ar[d]^{\iota_i}\\
X\ar@{^{(}->}[r]^{j^\ip} & X_\ip \ar@{^{(}->}[r] & Y}
\]
D'après (a), $L_\lbd$ est à faisceaux de cohomologie lisses sur
$X$, \mrs sur $Y$. D'après \cite[(3.7.1)]{Zheng-inte}, le morphisme
de changement de base
\[\iota_i^*Rf_* L_\lbd \to R(j_\ip')_* (\iota'_i)^* Rj^\ip_* L_\lbd\]
est un \iso. Cet \iso est $G_d(D_i)$-\eqrt. En appliquant
$(D_{d-1})$ à $j'_{\ip}$, il suffit de prouver la
$E$-compatibilité de $(Rj^\ip_*L_\lbd)_\ldI$. On est ramené à
prouver la stabilité par $Rf_*$ sous les hypothèses additionnelles
(c) et
\begin{cond}
(d) $X$ est le complémentaire d'un diviseur régulier $D$ de~$Y$.
\end{cond}
Refaisant un dévissage comme dans la démonstration de $(D_1)$,
\ops de plus
\begin{cond}
(e) $Y$ affine.
\end{cond}

En vertu de \ref{prop.DL}, \ops $G=\{1\}$. Pour $y\in |D|$, on prend
un sous-schéma régulier $C$ de $Y$ de dimension~$1$ tel que $C\cap D =
\{y\}$. Alors d'après \cite[3.7 (ii)]{Zheng-inte}, $(Rf_* L_\lbd)|C
\simeq Rf'_* (L_\lbd|C\cap X)$, où $f': C\cap X \hra C$. Il suffit
donc d'appliquer l'hypothèse (B). Ceci achève la démonstration de
$(D_d)$.
\end{proof}
\end{ssect}

La démonstration de $(D_d)$, $d\ge 2$ peut aussi être achevée par la
méthode de la démonstration de \cite[Th.\ finitude, 2.4]{SGA4d}.
Sous les hypothèses (a) à (e), on plonge $Y/G$ dans un espace
affine $\bA^n_\eta$. On note $g : Y \to \bA^n_\eta$. Pour tout $1\le
i\le n$, on note $p_i : \bA^n_\eta \to \bA^1_\eta$ la $i$-ième
projection et on considère
\[
\xymatrix{X\ar@{^{(}->}[r]^{f}\ar[rd] & Y\ar[d]^{p_i g} \\
&\bA^1_\eta}
\]
D'après le théorème de changement de base générique (\ibid[1.9
(ii)] pour le cas des coefficients de torsion ; le cas $l$-adique
résulte de ce cas par l'argument usuel : passer aux gradués pour
la filtration $l$-adique sur un modèle entier), il existe un ouvert
dense $U_i$ de $\bA^1_\eta$ \tq $Rf_*$ commute à tout changement de
base $S'\to U_i$. Quitte à rétrécir $U_i$, \ops que les fibres de
$(p_i g)^{-1} (U_i) \to U_i$ sont de dimension $\le d-1$. Alors
$(D_{d-1})$ implique que
\[(Rf_* L_\lbd|(p_i g)^{-1}(U_i))_\ldI\]
est \Ecomp, pour tout $i$. Notons que $E= Y- \bigcup_{i=1}^{n} (p_i
g)^{-1} (U_i)$ est un ensemble fini, donc $E\cap (Y-X)$ est rare
dans $Y-X$, car $Y-X$ est purement de dimension $d-1\ge 1$. Comme
les $Rf_* L_\lbd|{Y-X}$ sont à faisceaux de cohomologie lisses
\cite[3.7 (i)]{Zheng-inte}, $(Rf_* L_\lbd|Y-X)_\ldI$ est \Ecomp en
vertu de la proposition suivante (où $m=d-1$).

\begin{prop}\label{prop.generique}
Soit $m \ge 0$. Faisons l'hypothèse $(D_m)$. Soient $(X,G)$ \tf sur
$\eta$ avec $X$ \reg, de dimension $\le m$, $(L_\lbd)_\ldI \in
\prod_{\ldI} \Dbc(X,G,\Q_{l_\lbd})$ à \faisx de cohomologie lisses.
Supposons qu'il existe $U\subset X$ un ouvert dense $G$-stable \tq
$(L_\lbd |U)_\ldI \in \Dbc(U,G,E)$. Alors $(L_\lbd)_\ldI$ est
\Ecomp.
\end{prop}

Rappelons que l'hypothèse $(D_m)$ sera satisfaite (une fois (A) et
(B) établis).

\begin{proof}
On a une filtration de $X$ par des ouverts $G$-stables
\[U = U_{0} \subset
U_{1} \subset \cdots \subset U_m = X\] \teq $U_{i} - U_{i-1}$ soit
\reg et purement de codimension $i$ dans $U_{i}$. On montre
$(L_\lbd|U_i)_\ldI \in \Dbc(U_i,G,E)$ par récurrence sur~$i$. Le
cas $i=0$ est une hypothèse de la proposition et le cas $i=m$
permet de conclure. Supposons cela vrai pour un $i$ \ver ant $0 \le
i \le m-1$. Soit $j_i : U_i \hra U_{i+1}$. On a la formule de
projection
\[
L_{\lbd}|U_{i+1} \otimes Rj_{i*} \Q_{l_\lbd} \simto Rj_{i*} (L_\lbd
| U_i).
\]
D'après l'\hyp de \rec et $(D_m)$, $Rj_{i*} (L_\lbd | U_i) \in
\Dbc(U_{i+1},G,E)$. On a un \iso $G$-\eqrt
\[R^q j_{i*} \Q_l =
\begin{cases}
\Q_{l,U_{i+1}} & \text{si $q=0$,}\\
\Q_{l,U_{i+1}-U_i}(-i-1) & \text{si $q=2i+1$,}\\
0 & \text{sinon.}
\end{cases}
\]
Donc pour tout $x\in |U_{i+1}-U_i|$, tout $\xb \to x$ et tout
$(g,\phi)\in G_d(x,\xb)$ avec $\rho(\phi) = F_0^n$, $n>0$ un entier
(on a utilisé les notations de \ref{ss.compte}), on a
\[\Tr((g,\phi),(Rj_{i*}\Q_l)_\xb )= 1- p^{n(i+1)}.\]
Donc on a $(L_{\lbd}|U_{i+1})_\ldI \in \Dbc(U_{i+1},G,E)$.
\end{proof}

\section{Indépendance de $l$ des cycles proches et fin de la démonstration de
\ref{theo.6op}}\label{s4}
\begin{ssect}\label{ss.DS}
Soient $S$ un trait hensélien, $l$ un nombre premier inversible sur
$S$. Soit $\bD_S$ la catégorie des couples $(X,S_1)$ où $S_1$ est
un trait fini sur~$S$ et $X$ est un schéma \tf sur $S_1$. Un
morphisme $(X_1,S_1) \to (X_2,S_2)$ dans $\bD_S$ est un couple
$(f,g)$ de flèches s'insérant dans un carré commutatif
\[\xymatrix{X_1\ar[r]\ar[d]^{f} &S_1\ar[d]^{g} \\ X_2\ar[r] & S_2}\]
On note par $s$ le point fermé de $S$ et par $\eta$ le point générique de $S$. Soit $\bD_{s,S}$ la sous-catégorie pleine de $\bD_S$ formée des
couples $(X,S_1)$ avec $X_\eta = \vide$.

À un objet $(X,S_1)$ de $\bD_{s,S}$, on associe le topos
$X\times_{s_1} \eta_1$ \cite[XIII 1.2.4]{SGA7}, où $s_1$ est le
point fermé de $S_1$, $\eta_1$ est le point générique de $S_1$.
Pour $(X,S_1,G)\in \Eq(\bD_{s,S})$, on définit comme dans
\ref{ss.Dbc} les catégories $\Mod_c(X\times_{s_1} \eta_1,G,\Qlb)$,
$\Dbc(X\times_{s_1} \eta_1,G,\Qlb)$, les six opérations et la
dualité. Soit $\lbar{s_1} \to s_1$ un point géométrique
algébrique. Un faisceau d'ensembles $\cF$ sur
$(X\times_{s_1}\eta_1,G)$ est un faisceau sur $X_{\lbar{s_1}}$, muni
d'une action continue (\cite[XIII 1.1.2]{SGA7}) de
\[G_{\lbar{s_1},\eta_1} = G\times_{\Gal(k_1/k_1^G)} \Gal(\lbar{K_1}/K_1^G), \]
compatible à l'action de $G_{\lbar{s_1},\eta_1}$ sur
$X_{\lbar{s_1}}$ (via $G_{\lbar{s_1},\eta_1} \to G_{\lbar{s_1}} =
G\times_{\Gal(k_1/k_1^G)} \Gal(\lbar{k_1}/k_1^G) $). Ici $k_1 =
\kappa(s_1)$, $\lbar{k_1} = \kappa(\lbar{s_1})$, $K_1 =
\kappa(\eta_1)$, $\lbar{K_1}$ est une clôture séparable de $K_1$.

Pour $(x,S_1,G)\in \Eq(\bD_{s,S})$ avec $x=\Spec \kappa(x)$, et
$\xb$ un point géométrique algébrique au-dessus de $x$, on pose
\[
G_{\xb,\eta_1} = G\times_{\Gal(\kappa(x)/\kappa(x)^G)}
\Gal(\kappa(\xb)/\kappa(x)^G) \times_{\Gal(\kappa(\xb)/k_1^G)}
\Gal(\lbar{K_1}/K_1^G).
\]
En d'autres termes, un élément de $G_{\xb,\eta_1}$ est un triplet
$(g,\phi,\psi)$ où $g\in G$, $\phi\in
\Gal(\kappa(\xb)/\kappa(x)^G)$, $\psi\in \Gal(\lbar{K_1}/K_1^G)$ \tq
$\phi$ soit la réduction de $\psi$ et que $\phi$ et $g$ induisent
le même automorphisme sur $\kappa(x)$. On a une équivalence de
\cats
\begin{align*}
\Mod_c(x\times_{s_1} \eta_1, G,\Qlb) &\to \Rep(G_{\xb,\eta_1},\Qlb)\\
L &\mapsto L_\xb.
\end{align*}
Pour $(X,S_1,G)\in \Eq(\bD_{s,S})$, $x$ un point fermé de $X$, on
pose
\[G_d(x,\xb,\eta_1) = G_d(x)_{\xb,\eta_1}.\]
Pour $L\in
\Dbc(X\times_{s_1} \eta_1, G,\Qlb)$, on pose $L_x = i_x^* L$, où
$i_x : (x,S_1,G_d(x)) \to (X,S_1,G)$.
\end{ssect}

On a la variante suivante de \ref{prop.DL1}.

\begin{prop}\label{prop.DL1eta}
  Soient $(X,S_1,G)\in \Eq(\bD_{s,S})$, $S'$ un trait fini étale
  sur $S$ d'extension résiduelle galoisienne, $T$ une composante de
  $S_1\times_S S'$. Notons $t$ le point fermé de $T$, $\eta_T$ le point générique de $T$,
  $H=G_d(T)$, $Y=X\times_{s_1} t$. Considérons $(Y,T,H) \to
  (X,S_1,G)$. Soient $y$ un point fermé de $Y$ au-dessus d'un point fermé $x$ de $X$, $\yb \to y$ au-dessus de $\xb \to x$.
  Alors
  l'image de l'homomorphisme induit $i: H_d(y,\yb,\eta_T) \to G_d(x,\xb, \eta_1)$ est
  \[G'= G_d(x)\times_{\Gal(\kappa(x)/\kappa(x)^{G_d(x)})} \Gal(\kappa(\xb)/\kappa(x)^{G_d(x)})
    \times_{\Gal(\kappa(\xb)/k_1^{G_d(x)})} \Gal(\lbar{K_1}/K_1^{G_d(x)}\cdot K'),
  \]
  où $K'$ est le corps des fractions de $S'$.
\end{prop}

\begin{proof}
  Il est évident que $\Img i \subset G'$. Soit $(g,\phi,\psi)\in
  G'$. Posons $\Gamma= \Gal(K'/K)$. Formons le diagramme commutatif à carrés cartésiens
  \[\xymatrix{
    (Y,H) \ar[d] \ar[r] & (Y,(G\times \Gamma)_d(T)) \ar[d] \ar[r] & (X\times_s s', G\times \Gamma) \ar[d] \ar[r] & (X,G) \ar[d] \\
    (t,H) \ar[r] & (t,(G\times \Gamma)_d(T)) \ar[r] & (s_1\times_s s', G\times \Gamma) \ar[d] \ar[r] & (s_1,G) \ar[d] \\
    && (s',\Gamma) \ar[r] & (s,\{1\})   }
  \]
  D'après  une variante de \ref{prop.quot2} (b), l'homomorphisme
  $(G\times \Gamma)_d(y,\yb,\eta_T) \to G_d(x,\xb,\eta_1)$ est un
  isomorphisme. Soit $(g,\gamma,\phi',\psi') \in (G\times
  \Gamma)_d(y,\yb,\eta_T)$ l'image inverse de $(g,\phi,\psi)$. Alors
  $\phi'|k'=\Id_{k'}$. Comme $(\gamma,\phi')\in G_d(s',\sbar)$, on a
  $\gamma=1$, \ie $(g,\phi,\psi) = i(g,\phi',\psi')$.
\end{proof}

\begin{ssect}\label{ss.Cmgeta}
  Soit $(S_{(m)})_{m\ge 1}$ un système comme dans \ref{ss.Cmg}.
  On désigne par $\Eq^\adm(\bD_{s,S})$ la sous-catégorie pleine de
  $\Eq(\bD_{s,S})$ formée des triples $(X,S_1,G)$ où $G$ agit sur
  $X$ de façon admissible.
  Pour $(X,S_1,G) \in \Eq^\adm(\bD_{s,S})$, on applique \ref{ss.Cmg}
  à $(X,G) \to (S_1,G)$ au-dessus de $S$. Pour $m\ge 1$, $g\in G$,
  choisissons une composante $T_{m,g}$ de $S_1^{(m,g)}$ et posons
  $Y_{m,g} = X^{(m,g)} \times_{S_1^{(m,g)}} T_{m,g}$. On va définir un
  foncteur
  \[C_{m,g} : \Dbc(X\times_{s_1} \eta_1,G,\Qlb) \to \Dbc(Y_{m,g}
  \times_{t_{m,g}} \eta_{T_{m,g}},\Qlb),
  \]
  analogue du $C_{m,g}$ de \ref{ss.Cmg}. Le groupe $\Z/n_g\Z$ agit diagonalement sur $S_1\times_S S_{(n_g m)}$.
  Choisissons $S'$ une
  composante de $S_1\times_S S_{(n_g m)}$ au-dessus de $T_{m,g}$, $X'= X\times_{S_1}
  S'$. On note
  \[\xymatrix{(X,S_1,G) & \ar[l]_-{d} (X',S',(\Z/n_g\Z)_d(S'))\ar[r]^-{e} &(Y_{m,g},T_{m,g},\{1\}). } \]
  On pose $C_{m,g} = e_* d^*$, qui ne dépend pas du choix de $S'$.
\end{ssect}

\begin{ssect}\label{ss.Ecompeta}
Supposons dorénavant, sauf en \ref{ss.RP},
que le corps résiduel $k$ de $S$ soit $\F_{p^f}$. Alors
$G_d(x,\xb,\eta_1)= G_d(x) \times_{\Gal(\kappa(x)/k)}
\Gal(\lbar{K_1}/K_1^{G_d(x)})$. Soient $E$, $I$, $\gamma$ comme dans
\ref{ss.compte}. Pour $(X,S_1,G)\in \Eq(\bD_{s,S})$, on dit que
$(L_\lbd)_\ldI \in \prod_{\ldI} \Dbc(X\times_{s_1}
\eta_1,G,\lbar{\Q_{l_\lbd}})$ est $(E,I,\gamma)$-\emph{compatible}
(ou \Ecomp s'il n'y a pas de confusion à craindre) si pour tout
$x\in |X|$, tout $\xb \to x$ et tout $(g,\phi,\psi) \in
G_d(x,\xb,\eta_1)$ avec $\rho(\phi)$ une puissance entière de
$F_0^f$, le système $(\Tr((g,\phi,\psi),(L_\lbd)_\xb))_\ldI$ est
$(E,I,\gamma)$-compatible. Ici $\rho :
\Gal(\kappa(\xb)/\kappa(x)^{G_d(x)}) \hra \Aut(\kappa(\xb))$ est
l'inclusion, $F_0\in\Aut(\kappa(\xb))$ désigne le Frobenius
géométrique absolu $a\mapsto a^{1/p}$.

Comme en \ref{prop.pos}, la définition ne changera pas si on se
restreint aux triplets $(g,\phi,\psi) \in G_d(x,\xb,\eta_1)$ avec
$\rho(\phi) = F_0^{fm}$, où $m$ est un entier $\ge N(x)$, $N(x)$
est une constante qui ne dépend que de $x$.
\end{ssect}

\begin{prop}\label{prop.Xseta}
La \Ecompte ainsi définie est stable par les six opérations et la
dualité.
\end{prop}

\begin{proof}
  Cela découle du cas fini de \ref{theo.6op} et \ref{coro.D}.
  À titre d'exemple, prouvons la
  stabilité par $D$. Soient $(X,S_1,G)\in\Eq(\bD_{s,S})$,
  $(L_\lbd)_\ldI \in \Dbc(X\times_{s_1} \eta_1,G,E)$. Il suffit
  de montrer que pour $x\in |X|$, $\xb \to x$, $(g,\phi,\psi) \in
  G_d(x,\xb,\eta_1)$ avec $\rho(\phi) = F_0^{mf}$, $m\ge 1$,
  $(\Tr((g,\phi,\psi),(DL_\lbd)_\xb))_\ldI$ est $E$-compatible.
  Soient $s'= \Spec \F_{p^{fm}}$, $S'$ un trait étale sur~$S$ de point fermé $s'$.
  Appliquons \ref{prop.DL1eta}. On a  $(g,\phi,\psi)\in
  H_d(y,\yb,\eta_T)$. L'unique homomorphisme continu
  $\Gal(\kappa(\bar{t})/\kappa(t)^H) \to
  \Gal(\kappa(\lbar{\eta_T})/\kappa(\eta_T)^H)$ qui envoie $\phi$
  sur $\psi$ définit un foncteur $F: \Dbc(Y\times_t \eta_T, H,\Qlb) \to
  \Dbc(Y,H,\Qlb)$ de sorte que
  \[\Tr((g,\phi,\psi),(DL_\lbd)_\xb) =
  \Tr((g,\phi),(F(DL_\lbd)_{t})_\yb).
  \]
  Comme $F(DL_\lbd)_{t} \simeq
  DF(L_\lbd)_{t}$, il suffit d'appliquer le cas fini de
  \ref{coro.D}.
\end{proof}

Pour $m\ge 1$, soit $S_{(m)}$ un trait étale sur $S$ de corps
résiduel $\F_{p^{fm}}$. L'isomorphisme \eqref{eq.Z/mZ} induit une
action de $\Z/m\Z$ sur $S_{(m)}$. Appliquons \ref{ss.Cmgeta} au
système $(S_{(m)})_{m \ge 1}$.

\begin{prop}\label{prop.DLeta}
  Soient $(X,S_1,G)\in \Eq^\adm(\bD_{s,S})$, $(L_\lbd)_\ldI \in \prod_{\ldI} \Dbc(X\times_{s_1}\eta_1,G,\Qllb)$.
  Pour que $(L_\lbd)_\ldI$ soit $E$-compatible, il faut et il suffit
  que pour tout $m\ge 1$ et tout $g\in G$,
  \[(C_{m,g} L_\lbd)_{\ldI} \in \prod_\ldI \Dbc(Y_{m,g}\times_{t_{m,g}} \eta_{T_{m,g}},\Qllb) \]
  soit $E$-compatible.
\end{prop}

\begin{proof}
  La nécessité est claire. Soit $(L_\lbd)_\ldI \in \prod_{\ldI} \Dbc(X\times_{s_1} \eta_1,
  \Qllb)$ tel que pour tout $m\ge 1$ et tout $g\in G$,
  $(C_{m,g}L_\lbd)_\ldI$ soit $E$-compatible. Comme $e^* e_* =\Id$ (variante de \ref{prop.quot}
  (b)), on a $(d^* L_\lbd) \in \Dbc(X'\times_{s'} \eta', \Z/n_g\Z
  ,E)$. Il suffit de montrer que pour tout $x\in |X|$, tout $\xb \to x$ et tout $(g,\phi,\psi)\in
  G_d(x,\xb,\eta_1)$ avec $\rho(\phi) = F_0^{fm}$, $m\ge 1$,
  $(\Tr((g,\phi,\psi),(L_\lbd)_\xb))_\ldI$ est $E$-compatible.
  Notons $T$ la composante de $S_1\times_S S_{(m)}$ dominée par
  $S'$, $H = G_d(T)$, $Y=X\times _{s_1} t$. Soient $y\in |Y|$ au-dessus de
  $x$, $\yb \to y$ au-dessus de $\xb\to x$.
  D'après \ref{prop.DL1eta}, on a $(g,\phi,\psi) \in H_d(y,\yb,\eta_T)$, correspondant à
  $(\lbar{1},\phi,\psi) \in (\Z/n_g\Z)_d(y,\yb,\eta_T)$.
  Formons le diagramme commutatif à carrés cartésiens
  \[\xymatrix@C=1em{
  (X',(\Z/n_g\Z)_d(S'))\ar[r]^-{(\Id,\Delta)}\ar[d] & (X',(\Z/n_g\Z\times \Z/n_g\Z)_d(S'))\ar[r]\ar[d]
  & (Y\times_{s_{(m)}} s_{(n_g m)}, \Z/n_g\Z\times \Z/n_g\Z)\ar[r]\ar[d] & (Y,\Z/n_g\Z)\ar[d]\\
  (S',(\Z/n_g\Z)_d(S'))\ar[r]^-{(\Id,\Delta)} & (S',(\Z/n_g\Z\times \Z/n_g\Z)_d(S'))\ar[r]
  & (t\times_{s_{(m)}} s_{(n_g m)}, \Z/n_g\Z\times \Z/n_g\Z)\ar[r]\ar[d] & (t,\Z/n_g\Z)\ar[d]\\
  && (s_{(n_g m)},\Z/n_g\Z)\ar[r] & (s_{(m)},\{1\})}
  \]
  Soient $x'\in |X'|$ au-dessus de $y$, $\lbar{x'}\to x'$ au-dessus de $\yb \to y$.
  D'après  une variante de
  \ref{prop.quot2} (b), l'homomorphisme
  $(\Z/n_g \Z \times \Z/n_g \Z)_d(x',\lbar{x'}, \eta') \to (\Z/n_g \Z)_d(y,\yb,\eta_T)$
  est un isomorphisme. Soit $(\bar{1},b,\phi,\psi)\in (\Z/n_g \Z \times \Z/n_g \Z)_d(x',\lbar{x'}, \eta')$
  l'image inverse de $(\bar{1},\phi,\psi)$. Alors $(b,\phi)\in (\Z/n_g\Z)_d(\eta_{n_g
  m},\etab)$, donc $b=\bar{1}$. Bref, $d_{\lbar{x'},\eta'} : (\Z/n_g\Z)_d(x',\lbar{x'},\eta') \to
  G_d(x,\xb,\eta_1)$ envoie $(\bar{1},\phi,\psi)$ sur
  $(g,\phi,\psi)$. Donc les
  \[\Tr((g,\phi,\psi),(L_\lbd)_\xb ) =\Tr((\bar{1},\phi,\psi), (d^* L_\lbd)_{\lbar{x'}}) \]
  forment  un système $E$-compatible.
\end{proof}

\begin{ssect}\label{ss.RP}
Soit $(X,S_1,G)\in \Eq(\bD_S)$.  On a le foncteur des cycles proches
\[R\Psi_{X/S_1} : \Dbc(X_{\eta_1},G,\Qlb) \to \Dbc(X_{s_1}\times_{s_1} \eta_1,G,\Qlb).\]
Un morphisme dans $\Eq(\bD_S)$ de la forme $(f,\Id,\al): (X,S_1,G)
\to (Y,S_1,H)$ induit un morphisme $R\Psi_{Y/S_1}
R(f_{\eta_1},\al)_* \to R(f_{s_1},\al)_* R\Psi_{X/S_1}$, qui est un
isomorphisme lorsque $f$ est propre.

Pour $(X,S_1,G)\in \Eq(\bD_{S})$, $L\in \Dbc(X_{\eta_1},G,\Q_l)$, on
a, avec les notations de \ref{ss.Cmg} et \ref{ss.Cmgeta},
\beq\label{eq.CmgRP}
C_{m,g}R\Psi_{X/S_1} L \simeq
R\Psi_{Y_{m,g}/T_{m,g}}(C_{(X_{\eta_1},G),m,g} L)_{\eta_{T_{m,g}}}.
\eeq
\end{ssect}

Le résultat principal de ce \S\ est le suivant.

\begin{theo}\label{theo.RP}
Soient $S_1$ un trait fini sur $S$, $X$ un schéma \tf sur $S_1$,
$G$ un groupe fini agissant sur $X\to S_1$ par des
$S$-automorphismes. Alors $R\Psi_{X/S_1}$ préserve la \Ecompte.
\end{theo}

Pour $S$ d'égale \cara (ou du moins, pour $S$ le hensélisé en un
point fermé d'une courbe lisse sur $\F_{p^f}$), le résultat était
connu de Gabber, mais non publié.

Pour $(g,\phi,\psi)\in G_{\lbar{s_1},\eta_1}$ avec $\rho(\phi) =
F_0^{fn}$, $n\in \Z$, \ref{theo.RP} combiné avec \cite[2.4,
5.3]{Zheng-inte} implique que
\[\sum_{i\in \Z} (-1)^i \Tr((g,\phi,\psi), H^i_c(X_{\sbar},(R\Psi_{X/S_1} \Q_l)_{\sbar}))
\in p^{fd\min\{0,n\}} \Z
\]
est indépendant de~$l$, où $d=\dim X_\eta$. On retrouve ainsi le
cas (très particulier) de \cite[6.1.3]{Mieda} où $\Gamma$ est le
graphe transposé de $a_g$.



La démonstration de \ref{theo.RP} sera donnée en \ref{ss.demoRP}.
Commençons par les deux cas particuliers suivants.

\begin{lemm}\label{ss.RPtriv}
  Supposons de plus $X$ quasi-fini sur $S_1$. Alors $R\Psi_{X/S_1}$ préserve la \Ecompte.
\end{lemm}

\begin{lemm}\label{prop.RPRj}
Soient $G$ un groupe fini agissant sur un trait
$S_1$ fini sur $S$, $(X,Z)$ un couple semi-stable $G$-strict sur
$S_1$, $u : X-Z \hra X_{\eta_1}$, $(L_\lbd)_\ldI\in \Dbc(X-Z,G,E)$
avec $L_\lbd$ \ver ant \ref{cond.mr}. Alors on a
\[(R\Psi_{X/S_1} Ru_* L_\lbd)_\ldI \in \Dbc(X_{s_1}\times_{s_1}
\eta_1,G,E).\]
\end{lemm}

Pour la notion de couple semi-stable $G$-strict, on renvoie aux
définitions qui précèdent \ref{lemm.deJong}.

\begin{proof}[Démonstration de \ref{ss.RPtriv}]
\Ops $X$ réduit. Soit $f: X'\to X$ le normalisé de $X$. Alors $
f_\eta = \Id$ et $R\Psi_{X/S_1} \simeq Rf_{s*} R\Psi_{X'/S_1}$. Donc
\ops $X$ normal de sorte que $X$ est somme disjointe de spectres
d'extensions finies de $s_1$, de spectres d'extensions finie de
$\eta_1$ et de traits finis sur $S_1$. \Ops que $X=S_3$ est un trait
fini sur $S_1$. Soient $K_2$ une sous-extension maximale \nre de
$K_3/K_1$, $S_2$ le normalisé de $S_1$ dans $K_2$, $g : S_3 \to
S_2$. Alors $R\Psi_{S_3/S_1} \simeq R\Psi_{S_2/S_1} Rg_{\eta_1*}$.
Comme $R\Psi_{S_2/S_1}: \Dbc(\eta_2,G,\Qlb) \to
\Dbc(s_2\times_{s_1}\eta_1,G,\Qlb)$ s'identifie à l'identité, il
résulte de \ref{rq.th} (iv) que $R\Psi_{S_3/S_1}$ préserve la
\Ecompte.
\end{proof}


\begin{proof}[Démonstration de \ref{prop.RPRj}]
Soit $D$ une composante de $X_s$. On note par $Z^{(D)}$ la réunion
des autres composantes de $Z$. Formons le diagramme commutatif à
carrés cartésiens
\[
  \xymatrix{X-Z \ar@{^{(}->}[dd]_{u}\ar@{^{(}->}[r] & X-Z^{(D)}\ar[dd] &
        D^\circ\ar[l]\ar@{^{(}->}[d]^j\\
    && D\ar@{^{(}->}[d]\\
    X_\eta\ar@{^{(}->}[r] & X& X_s \ar[l]}
\]
Comme $L_\lbd$ est \mr sur $X$, on a un isomorphisme
\cite[(5.6.2)]{Zheng-inte}
\[(R\Psi_{X/S_1} Ru_* L_\lbd)|D \simeq Rj_* R\Psi_{(X-Z^{(D)})/S_1} L_\lbd. \]
Cet isomorphisme est $G_d(D)$-équivariant. Appliquant
\ref{theo.6op} dans le cas fini (variante \ref{prop.Xseta}) à
$Rj_*$, on est ramené à prouver la proposition au cas où
$Z=X_{s}$ est un diviseur \reg. Alors $X\to S_1$ est lisse. D'après
\ref{prop.DLeta} et \eqref{eq.CmgRP}, \ops $G=\{1\}$. Soit $x\in
|X_s|$. Il existe un sous-schéma fermé $Y$ de $X$ qui est un trait
fini étale sur $S'$ de point fermé~$x$. D'après \ibid[5.6 (iii)],
$(R\Psi_{X/S_1} L_\lbd)_x \simeq R\Psi_{Y/S_1} (L_\lbd)_{\eta_Y}$.
Il suffit donc d'appliquer \ref{ss.RPtriv} à $Y\to S_1$.
\end{proof}

\begin{ssect}\label{ss.demoAB}
\begin{proof}[Fin de la démonstration de \ref{theo.6op}]
Pour finir la démonstration, il reste à traiter le cas des
courbes. Plus précisément, il reste à prouver les énoncés (A)
et (B) de \ref{prop.sousAB} pour $K$ un corps local.

Prouvons d'abord (B). 
On peut supposer $X$ géométriquement irréductible. On applique
\ref{lemm.deJong} à $G=\{1\}$ et $X\to S$ et on obtient un groupe
fini $G'$ et un diagramme commutatif $G'$-équivariant
\[
\xymatrix{U'\ar@{^{(}->}[r]^{j'}\ar[d]^{v_{U}}\ar@{}[rd]|\Box
& X'\ar[d]^{v}\ar@{^{(}->}[r]^{j_1} & Y\ar[r]^{g} & S'\ar[d]^{u}\\
U \ar@{^{(}->}[r]^{j} & X \ar[rr] & & S}
\]
avec $(Y,Y-U')$ un couple semi-stable $G'$-strict sur $S'$,
$Y-U'=Y_{s'} \cup H$, où $H$ est la réunion d'un nombre fini de
sections de $g$. Quitte à rétrécir $Y$, on peut supposer que $X'=
Y_{\eta'}$. Notons que $v$ est fini. Le cône de $L_\lbd \to
(v_{U},\al)_* v_U^* L_\lbd$ est à support de dimension $0$
(\ref{prop.quot} (a)), où $\al : G'\to \{1\}$. D'après \ref{rq.th}
(iv), $((v_U,\al)_* v_U^* L_\lbd)_\ldI$ est $E$-compatible, donc il
suffit de \ver er que les $Rj_* (v_{U},\al)_* v_U^* L_\lbd =
(v,\al)_* Rj'_* v_U^* L_\lbd$ forment un système \Ecomp. Comme $(v,
\al)_*$ préserve la $E$-compatibilité, il suffit donc de vérifier
que $(Rj'_*L'_\lbd)_\ldI$ est $E$-compatible, où $L'_\lbd = v_U^*
L_\lbd$.
D'après \cite[5.6 (i)]{Zheng-inte}, on a $i_{s'}^* R\Psi_{Y/S'}
Rj'_* L'_\lbd \simeq R\Psi_{H/S'} i_{\eta'}^* Rj'_* L'_\lbd$, où
$i_{s'} : H\cap Y_{s'} \to Y_{s'}$, $i_{\eta'} : H\cap Y_{\eta'} \to
Y_{\eta'}$. Cet isomorphisme est $G'$-\eqrt. Comme $R\Psi_{H/S'} :
\Dbc(H_{\eta'},G',\Qlb) \to \Dbc(H_{s'}\times_{s'}\eta', G',\Qlb)$
s'identifie à l'identité, il suffit alors d'appliquer
\ref{prop.RPRj}.

Prouvons (A). Plus généralement, prouvons que si $G$ est un groupe
fini agissant sur $a_X$ comme dans (A) par des $\eta$-automorphismes,
et $(L_\lbd)_\ldI \in \Dbc(X,G,E)$ avec tous les $L_\lbd$ vérifiant
\ref{cond.mr}, alors $(Ra_{X!}L_\lbd)_\ldI \in \Dbc(\eta_1,G,E)$.
\Ops $X$ \geoqt \irr sur $\eta_1$. Soit $S_1$ le normalisé de $S$
dans $K_1$. L'action de $G$ sur $\eta_1$ se prolonge en une action
sur $S_1$. On applique \ref{lemm.deJong} à $X\to S_1$ et obtient un
homomorphisme surjectif $\al : G' \to G$ et un diagramme commutatif
$G'$-équivariant
\[\xymatrix{X'\ar[d]^{v}\ar@{^{(}->}[r]^{j} & Y\ar[r]^{g} & S'\ar[d]^{u}\\
 X \ar[rr] & & S_1}
\]
avec $(Y,Y-X')$ un couple semi-stable $G'$-strict sur $S'$ et $g$
propre. Le cône de $L_\lbd \to (v,\al)_* v^* L_\lbd$ est à
support de dimension $0$. Donc il suffit de \ver er que $(Ra_{X!}
(v,\al)_* v^* L_\lbd)_\ldI$ appartient à  $\Dbc(\eta_1,G,E)$. Soit
$L'_\lbd = v^* L_\lbd$. On a $Ra_{X!} (v,\al)_* L'_\lbd
 = (u_{\eta_1},\al)_* Rg_{\eta'*} Rj_{\eta' !} L'_\lbd$,
 \[R\Psi_{S'/S'} Rg_{\eta'*} Rj_{\eta' !} L'_\lbd \simeq
Rg_{s*} R\Psi_{Y/S'} Rj_{\eta'!} L'_\lbd.
\]
Comme $R\Psi_{S'/S'} : \Dbc(\eta',G,\Qlb) \to \Dbc(s'\times_{s'}
\eta',G,\Qlb)$ s'identifie à l'identité, il suffit de \ver er que
$(R\Psi_{Y/S'} Rj_{\eta'!} L'_\lbd)_\ldI$ appartient à
$\Dbc(Y_{s'},G,E)$. On a le \trdist
\[R\Psi_{Y/S'} Rj_{\eta'!} L'_\lbd \to
R\Psi_{Y/S'} Rj_{\eta'*} L'_\lbd \to R\Psi_{Y/S'} i_* i^*
Rj_{\eta'*} L'_\lbd \to,\] où $i: Y_{\eta'}-X'\to Y_{\eta'}$. Pour
le deuxième terme, on applique \ref{prop.RPRj}. Pour le troisième,
on applique (B), \ref{prop.DL} et \ref{ss.RPtriv}. On a achevé la
démonstration de \ref{theo.6op}.
\end{proof}
\end{ssect}

\begin{ssect}\label{ss.demoRP}
\begin{proof}[Démonstration de \ref{theo.RP}]
Supposons $\sharp I = 2$.  Il suffit de montrer pour tout $d\in \N$
:
\begin{cond}
$(*_d)$ Pour tout $(X,S_1,G)\in \Eq(\bD_S)$ et toute famille
$E$-compatible $(L_\lbd)_\ldI$ sur $X_{\eta_1}$, de support de
dimension $\le d$, la famille $(R\Psi_{X/S_1} L_\lbd)_\ldI$ est
$E$-compatible.
\end{cond}
On montre $(*_d)$ par \rec sur $d$. \Ops $\dim X_\eta \le d$ et $X$
réduit. Si $d= 0$, alors $X$ est quasi-fini sur $S_1$. Ce cas est
déjà traité en \ref{ss.RPtriv}. Soit $d\ge 1$ \tq $(*_{d-1})$
soit vrai. Prouvons $(*_d)$. Soit $w: U\hra X_\eta$ un ouvert dense
$G$-stable affine et normal. Comme $(Rw_*w^* L_\lbd)_\ldI$ est
$E$-compatible (le cas local de \ref{theo.6op}) et le cône
$M_\lbd$ de $L_\lbd \to Rw_*w^* L_\lbd$ est à support dans
$Y_\eta-U$ de dimension $\le d-1$, $(R\Psi_{X/S_1}M_\lbd)_\ldI$ est
$E$-compatible en vertu de l'hypothèse de récurrence $(*_{d-1})$.
Donc il suffit de montrer que $(R\Psi_{X/S_1}Rw_*w^*L_\lbd)_\ldI$
est $E$-compatible. Quitte à remplacer $U$ par un sous-schéma
ouvert fermé $G$-stable de~$U$, on peut supposer que $G$ agit
transitivement sur $\pi_0(U)$. Soit $C$ une composante connexe de
$U$. Alors $(C,G_d(C)) \to (U,G)$ est cocartésien
(\ref{sss.cocart}). D'après \ref{sss.cocarteb}, quitte à remplacer
$U$ par $C$, $X$ par l'adhérence de $C$ et $G$ par $G_d(C)$, on
peut supposer $X$ irréductible. En appliquant \ref{lemm.fact},
quitte à rétrécir $U$, \ops qu'il existe $(f,\al): (Z,G')\to
(U,G)$ faisant de $Z$ un $(\Ker\al)$-torseur sur $U$ et \tq les $f^*
L_\lbd$ \ver ent \ref{cond.mr}. On applique \ref{lemm.compact} au
morphisme composé $Z\sto{f} U \hra X$ et trouve un \diagcomm
$G'$-\eqrt
\[\xymatrix{Z\ar[d]^{f}\ar@{^{(}->}[r]^j & Y\ar[d]^{g}\\
U\ar@{^{(}->}[r] & X}
\]
où $g$ est propre et $j$ est une immersion ouverte. Notons $L_\lbd'
= w^* L_\lbd$. On a $L_\lbd' \simeq R(f,\al)_* f^* L_\lbd'$
(\ref{prop.quot} (a)), donc
\begin{multline*}
R\Psi_{X/S_1} Rw_* L'_\lbd \simeq R\Psi_{X/S_1} Rw_* R(f,\al)_* f^*
L'_\lbd =
R\Psi_{X/S_1} R(g_{\eta_1},\al)_* Rj_{\eta_1 *} f^* L'_\lbd\\
\simeq R(g_{s_1},\al)_* R\Psi_{Y/S_1} Rj_{\eta_1 *} f^* L'_\lbd.
\end{multline*}
Quitte à changer les notations, il suffit de montrer l'assertion
suivante
\begin{cond}
  $(\dag_d)$ Soient $(X,S_1,G)\in\Eq(\bD_S)$ \tq $\dim X_{\eta_1} \le
  d$, $w: U\hra X_{\eta_1}$ ouvert $G$-équivariant, $(L_\lbd)_\ldI\in
  \Dbc(U,G,E)$ \tq les $L_\lbd$ vérifient \ref{cond.mr}.
  On a $(R\Psi_{X/S_1} Rw_* L_\lbd)_\ldI \in
  \Dbc(X_{s_1},G,E)$.
\end{cond}

Soit $K'$ une extension finie galoisienne de $K_1^G$ contenant $K_1$
telle que les composantes \irrs de $U\otimes_{K_1} K' $ soient
\geoqt \irrs. Soient $u : S'\to S_1$ le normalisé de $S_1$ dans
$K'$, $X' = X\times_{S_1} S'$, $w'=w\times_{S_1} S'$,  $G' =
G\times_{\Gal(K_1/K_1^G)}\Gal(K'/K_1^G)$, $\al: G'\to G$ la projection. Alors $R\Psi_{X/S_1} Rw_*
\simeq (u_{X_{s_1}}, \al)_* u_{X_{s_1}}^* R\Psi_{X/S_1} Rw_* \simeq
(u_{X_{s_1}}, \al)_* R\Psi_{X'/S'} Rw'_* u_{U}^*$. Donc \ops que les
composantes \irrs de $U$ sont \geoqt \irrs. Comme au début de la
démonstration, quitte à remplacer $X$ par $X_\red$ et à
rétrécir $U$, \ops en plus $U$ normal. Quitte à remplacer $U$ par
une composante connexe $C$, $X$ par l'adhérence de $C$, et $G$ par
$G_d(C)$ (\ref{sss.cocart}, \ref{sss.cocarteb}), \ops $X$ intègre
et $X_{\eta_1}$ \geoqt \irr. Le problème étant local, \ops en plus
$X$ \sep (voire affine). Comme $d\ge 1$, on applique
\ref{lemm.deJong} et obtient $\beta :G'\to G$ et un \diagcomm
$G'$-\eqrt
\[\xymatrix{
    V\ar@{^{(}->}[r]^{j} &
    U'\ar@{^{(}->}[r]^{w'}\ar[d]^{v_{1U}}\ar@/_2pc/[dd]_(.3){v_U} &
    X_\eta'\ar[d]\ar@{^{(}->}[r] & X' \ar[rd]\ar[d]^{v_1}\ar@/_2pc/[dd]_(.3){v}\\
    & U\times_{\eta_1} \eta' \ar@{^{(}->}[r]^{w_{S'}}\ar[d]^{u_{U}} &
    X_\eta\times_{\eta_1} \eta'\ar[d]\ar@{^{(}->}[r] & X\times_{S_1} S' \ar[r]\ar[d]^{u_X} & S'\ar[d]^{u}\\
    & U\ar@{^{(}->}[r]  & X_\eta\ar@{^{(}->}[r] & X \ar[r] & S_1}
\]
à carrés cartésiens \tsq $(u,\beta)$ et $(v,\beta)$  soient des
altérations galoisiennes, et que $(X',X'-V)$ soit un couple
semi-stable $G'$-strict sur $S'$. Le $V$ du diagramme est le
$X'-Z_d$ de \ref{lemm.deJong}. Soit $G_{S'} = G'/\Ker \beta \cap
\Ker (G'\to \Aut(S'))$. Alors $\beta$ se factorise en
$G'\sto{\beta_1} G_{S'} \sto{\al} G$, $G_{S'}$ agit sur $S'$ et le
morphisme $((v_{1U})_\red,\beta_1) : (U',G') \to ((U\times_{\eta_1}
{\eta'})_\red,G_{S'})$ est une altération galoisienne. D'après
\ref{prop.quot} (a),
\[R\Psi_{X/S_1}
Rw_* L_\lbd \simeq (u_{X_s},\al)_* (u_{X_s},\al)^* R\Psi_{X/S_1} Rw_*
L_\lbd \simeq (u_{X_{s}},\al)_* R\Psi_{X\times_{S_1}S'/S'} Rw_{S'*}
(u_{U},\al)^* L_\lbd,
\]
où $u_{X_s}$ est le changement de base de $u_X$ par $X_s \to X$, et
le cône de
\[(u_{U},\al)^* L_\lbd \to R(v_{1U},\beta_1)_*
(v_{1U},\beta_1)^* (u_{U},\al)^* L_\lbd
\]
est à support de dimension $\le d-1$. Grâce au cas local de
\ref{theo.6op} et l'hypothèse de récurrence $(*_{d-1})$, il suffit
donc de voir que les
\begin{multline*}
(u_{X_{s}},\al)_* R\Psi_{X\times_{S_1}S'/S'} Rw_{S'*}
R(v_{1U},\beta_1)_* (v_{1U},\beta_1)^* (u_{U},\al)^*
L_\lbd\\
\simeq R(v_{X_{s}},\beta)_* R\Psi_{X'/S'} Rw'_* (v_U,\beta)^* L_\lbd
\end{multline*}
forment un système \Ecomp. Ici $v_{X_s}$ est le changement de base
de $v$ par $X_s\to X$. Quitte à remplacer $v_U^* L_\lbd$ par $Rj_*
j^* v_U^* L_\lbd$, on est donc ramené à montrer $(\dag_d)$ dans le
cas où $(X,X-U)$ est un couple semi-stable $G$-strict. Ceci est
déjà fait en \ref{prop.RPRj}.
\end{proof}
\end{ssect}

La définition suivante généralise \ref{def.Ecomp}.

\begin{defn}
  Soit $(X,G)\in \Eq/S$. Un système $(L_\lbd)_\ldI \in \prod_\ldI
  \Dbc(X,G,\Qllb)$ est \emph{$E$-compatible} si $(j^*L_\lbd)_\ldI$
  et $(i^*L_\lbd)_\ldI$ sont $E$-compatibles, où $j: X_\eta \hra
  X$, $i: X_s \to X$.
\end{defn}

Pour $(X,S_1,G)\in \Eq(\bD_S)$, le foncteur des cycles évanescents
\[R\Phi_{X/S_1} : \Dbc(X,G,\Qlb)\to \Dbc(X_{s_1}\times_{s_1} \eta_1,G,\Qlb) \]
préserve la $E$-compatibilité. C'est une conséquence immédiate
de \ref{theo.RP}.

\begin{exem}
  Soient $T$ et $T'$ deux traits henséliens d'égale
  caractéristique, à corps résiduel $k=\F_{p^f}$, munis
  d'uniformisantes. Soit $(E,I,\gamma)$ comme plus haut avec $E$
  contenant les racines $p$-ièmes de l'unité. Soit $\psi : \F_p \to
  E^\times$ un caractère. Pour $\ldI$, on note
  $\psi_\lbd = \iota_\lbd \circ \psi : \F_p \to {\Qllb}^\times$. Alors les
  transformations de Fourier locales \cite[2.3.2.4]{Laumon}
  \[\cF_\psi^{(0,\infty')}, \cF_\psi^{(\infty,0')}, \cF_\psi^{(\infty,\infty')} :
  \Mod_c(\eta_T,\Qlb) \to \Mod_c(\eta_{T'},\Qlb)
  \]
  préservent la $E$-compatibilité. Plus précisément, pour
  $(V_\lbd)_\ldI \in \prod_\ldI \Mod_c(\eta_T,\Qllb)$,
  $E$-compatible, les systèmes
  \[(\cF_{\psi_\lbd}^{(0,\infty')} V_\lbd)_\ldI,
    (\cF_{\psi_\lbd}^{(\infty,0')} V_\lbd)_\ldI,
    (\cF_{\psi_\lbd}^{(\infty,\infty')} V_\lbd)_\ldI
  \]
  sont $E$-compatibles.
  Cela résulte de \ibid[2.4.2.1], de la stabilité de la
  $E$-compatibilité par $R\Phi$, et de \ref{prop.saeta}.
\end{exem}

La proposition suivante généralise \ref{theo.6op} et \ref{coro.D}.

\begin{prop}\label{prop.6op}
  La $E$-compatibilité sur $S$ est stable par les six opérations
  et la dualité.
\end{prop}

\begin{proof}
  La stabilité par $\otimes$ et $(f,\al)^*$ est triviale. Il suffit
  de prouver la stabilité par $R(f,\al)_*$ et $D$.

  On montre d'abord que pour $(X,G)\in \Eq/S$, $(L_\lbd)_\ldI\in
  \Dbc(X_\eta,G,E)$, on a $(Rj_* L_\lbd)_\ldI \in \Dbc(X,G,E)$, où $j: X_\eta \hra
  X$. La suite spectrale de
  Hochschild-Serre donne
  \[E_2^{pq} = H^p(I_\eta,R^q \Psi_{X/S} L_\lbd) \Rightarrow i^* R^{p+q} j_* L_\lbd, \]
  où $I_\eta=\Ker(\Gal(\etab/\eta)\to \Gal(\sbar/s))$, $i : X_s \to X$. Comme $(R\Psi_{X/S}
  L_\lbd)_\ldI$ est $E$-compatible (\ref{theo.RP}), il suffit de montrer
  que pour $(Y,G)\in \Eq/s$, $(M_\lbd)_\ldI \in \Dbc(Y\times_s\eta,G,E)$, on a $(R\Gamma(I_\eta,M_\lbd))_\ldI \in
  \Dbc(Y,G,E)$. Pour cela, on peut supposer que $Y$ est fini sur
  $s$, $\sharp I=2$, $[M_\lambda]=[\cF_\lambda]-[\cG_\lambda]$, $[\cF_\lambda], [\cG_\lambda]\in \mathrm{Mod}_c(Y\times_s\eta, G,\overline{\mathbb{Q}_{l_\lambda}})$ et $I_\eta$ agit sur $\cF_\lambda$ et $\cG_\lambda$ par un quotient fini $Q$.
  Alors
  \[\left([R\Gamma(I_\eta,M_\lbd)]\right)_\ldI = \left([\cF_\lambda^Q]-[\cF_\lambda^Q(-1)]-[\cG_\lambda^Q]+[\cG_\lambda^Q(-1)]\right)_\ldI \in K(Y,G,E). \]

  Prouvons la stabilité par $R(f,\al)_*$, où $(f,\al) : (X,G) \to
  (Y,H)$ est un morphisme dans $\Eq/S$. Soit $(L_\lbd)_\ldI \in
  \Dbc(X,G,E)$. On a les triangles distingués
  \begin{gather*}
    i_* Ri^! L_\lbd \to L_\lbd \to Rj_* j^* L_\lbd \to,\\
    R(fi,\al)_* Ri^! L_\lbd \to R(f,\al)_* L_\lbd \to R(fj,\al)_*
    j^* L_\lbd \to.
  \end{gather*}
  Comme $(Rj_* j^* L_\lbd)_\ldI$ est $E$-compatible,
  $(Ri^!L_\lbd)_\ldI$ l'est aussi. Il suffit donc d'appliquer
  \ref{theo.6op} pour $R(f_s,\al)_*$ et $R(f_\eta,\al)_*$.

  Prouvons la stabilité par $D$. Soit $(L_\lbd)_\ldI \in
  \Dbc(X,G,E)$. On a le triangle distingué
  \[i_* D_{(X_s,G)}(i^* L_\lbd) \to D_{(X,G)} L_\lbd \to Rj_* D_{(X_\eta,G)} j^* L_\lbd \to. \]
  Il suffit donc d'appliquer \ref{coro.D} pour $D_{(X_s,G)}$ et
  $D_{(X_\eta,G)}$.
\end{proof}

\begin{remq}
  Si on remplace les faisceaux usuels par les faisceaux de Weil
  \cite[1.1.10]{WeilII},
  on peut définir la \Ecompte de la même manière. Tous les résultats
  concernant la \Ecompte (notamment \ref{prop.6op}, \ref{theo.RP}, \ref{prop.generique},
  \ref{prop.saeta}) restent valables.
\end{remq}

\section{Appendice : Indépendance de $l$ sur les champs
algébriques}\label{s5}

\newcommand{\fV}{{\mathcal{V}}}
\newcommand{\cX}{{\mathcal{X}}}
\newcommand{\cY}{{\mathcal{Y}}}
\newcommand{\Ch}{\mathbfit{Ch}}
\newcommand{\Four}{\mathrm{Four}}
\newcommand{\fS}{\mathfrak{S}}
\newcommand{\G}{\mathbb{G}}
\newcommand{\Ftilde}{\widetilde{F}}
\newcommand{\pib}{\bar{\pi}}
\newcommand{\omi}{\omicron}

\newcommand{\DM}{de Deligne-Mumford\xspace}

Soient $K$ un corps fini de \cara $p$ ou un corps local de \cara
résiduelle~$p$ (\ref{ss.compte}), $\eta = \Spec K$. On note
$\Ch/\eta$ la \cat des $\eta$-champs algébriques
\cite[4.1]{Laumon-MB} \tf. Soit $l$ un nombre premier $\neq p$. Pour
$\cX \in \Ch/\eta$, on note $\Mod_c(\cX,\Qlb)$ la catégorie des
$\Qlb$-faisceaux constructibles sur le site lisse-étale de $\cX$
\ibid[12.1 (i)]. On dispose, par \cite{Laszlo-Olsson2}, d'une
catégorie triangulée $\Dbc(\cX,\Qlb)$ munie d'une $t$-structure
usuelle de c\oe ur $\Mod_c(\cX,\Qlb)$ et d'un formalisme de six
opérations. Pour $f: \cX\to \cY$ un morphisme de $\eta$-champs
algébriques \tf, on a
\begin{align*}
D_\cX : \Dbc(\cX,\Qlb)^\op &\to \Dbc(\cX,\Qlb),\\
-\otimes- : \Dbc(\cX,\Qlb) \times \Dbc(\cX,\Qlb) &\to \Dbc(\cX,\Qlb),\\
R\cHom_\cX(-,-) : \Dbc(\cX,\Qlb)^\op \times \Dbc(\cX,\Qlb) &\to \Dbc(\cX,\Qlb),\\
f^*, Rf^! : \Dbc(\cY,\Qlb) &\to \Dbc(\cX,\Qlb).
\end{align*}
Rappelons que $f$ est dit \emph{relativement de Deligne-Mumford}
\cite[7.3.3]{Laumon-MB}, si pour tout schéma affine $Y$ et tout
morphisme $Y\to \cY$, le produit fibré $Y\times_\cY \cX$ est un
$\eta$-champ de Deligne-Mumford. Dans ce cas, on a
\[
Rf_*, Rf_! : \Dbc(\cX,\Qlb) \to \Dbc(\cY,\Qlb).
\]
Si $\cX$ est un $\eta$-champ \DM \tf, $\Mod_c(\cX,\Qlb)$ s'identifie
à la \cat des $\Qlb$-faisceaux constructibles sur le site
\emph{étale} de $\cX$ \ibid[12.1 (ii)].

Le foncteur
\begin{align*}
F : \Eq/\eta &\to \Ch/\eta\\
(X,G) &\mapsto [X/G] = [X/G/\eta]\quad \text{\ibid[2.4.2, 3.4.1]}
\end{align*}
n'est ni fidèle ni plein. Notons $S$ l'ensemble des flèches
$(f,\al) :(X,G)\to (Y,H)$ dans $\Eq/\eta$ telles que $G$ soit
isomorphe à $H\times H'$, $\al$ s'identifie à la projection
$H\times H' \to H$, et $f$ fasse de $X$ un $H'$-torseur sur $Y$.
Notons $Q: \Eq/\eta \to (\Eq/\eta)(S^{-1})$ le foncteur de
localisation.

\begin{prop}\label{prop.ch}
  La catégorie localisée $(\Eq/\eta)(S^{-1})$ admet un calcul de
  fractions à droite au sens faible (\ie pour toute flèche $t$ dans $(\Eq/\eta)(S^{-1})$,
  il existe des flèches $s,t'$ dans $\Eq/\eta$ avec $s\in S$ telles que $t=
  Q(t')Q(s)^{-1}$). Le foncteur $F$ se factorise en
  \[(\Eq/\eta) \sto{Q} (\Eq/\eta)(S^{-1}) \sto{\Ftilde} \Ch/\eta, \]
  où $\Ftilde$ est un foncteur plein (\ie $\Ftilde$ envoie
  $\Hom_{(\Eq/\eta)(S^{-1})} (A,B)$ \emph{sur} $\Hom_{\Ch/\eta}(\Ftilde A, \Ftilde B)$,
  pour tous $A, B \in (\Eq/\eta)(S^{-1})$).
\end{prop}

\begin{proof}
  On vérifie que
  \begin{cond}
    (i) Pour tout $(X,G)\in \Eq/\eta$, $\Id_{(X,G)} \in S$.

    (ii) Pour tout $s_1,s_2\in S$, composables, $s_1\circ s_2\in S$.

    (iii) Tout diagramme dans $\Eq/\eta$ de la forme
    \[\xymatrix{&\ar[d]^s\\ \ar[r]_t &} \]
    où $s\in S$,
    peut être complété en un diagramme commutatif dans $\Eq/\eta$
    \[\xymatrix{\ar[r]\ar[d]_{s'} & \ar[d]^s\\ \ar[r]_t &} \]
    où $s'\in S$.
  \end{cond}
  Donc $(\Eq/\eta)(S^{-1})$ admet un calcul de fractions à droite
  au sens faible.

  Pour tout $s\in S$, $F(s)$ est un isomorphisme. Donc $F$ se
  factorise à travers $Q$. Montrons que $\Ftilde$ est plein.
  Soient $(X,G), (Y,H) \in \Eq/\eta$, $\al :[X/G] \to [Y/H]$ dans
  $\Ch/\eta$. Formons le diagramme $2$-commutatif de $\eta$-champs
  à carrés $2$-cartésiens
  \[\xymatrix{
    X' \ar[d] \ar[r] & X \ar[d] \\
    \cX \ar[d] \ar[r] & [X/G] \ar[d]^{\al} \\
    Y \ar[r] & [Y/H]   }
  \]
  Notons $s:(X',G\times H) \to (X,G)$, $t:(X',G\times H)\to (Y,H)$.
  Alors $s\in S$ et $Ft = \al \circ Fs$, \ie $\al=
  \Ftilde(ts^{-1})$.
\end{proof}


\begin{prop}\label{prop.omi}
  Pour $(X,G)\in \Eq/\eta$, il existe une
  équivalence de \cats canonique $\omi_{X,G}$
  qui rend $2$-commutatif le diagramme suivant
  \[\xymatrix{
  \Dbc([X/G],\Qlb) \ar[r]^{\omi_{X,G}}\ar[rd]_{p^*_{X,G}}
  &\Dbc(X,G,\Qlb)\ar[d]^{\omega_{X,G}} \\
  &\Dbc(X,\Qlb)}
  \]
  où $p_{X,G} : X\to [X/G]$ est la projection, $\omega_{X,G}$ est le
  foncteur oubli de l'action de $G$. Cette équivalence commute
  aux six opérations et  à la dualité. Plus précisément, on a
  \begin{gather}
    \label{eq.e.Hom}
    \omi_{X,G} R\cHom_{[X/G]}(-,-) \simeq
    R\cHom_{(X,G)}(\omi_{X,G}-,\omi_{X,G}-),\\
    \label{eq.e.tD}
    \omi_{X,G}(-\otimes -) \simeq (\omi_{X,G}-) \ot (\omi_{X,G}-), \quad \omi_{X,G}
    D_{[X/G]} \simeq D_{(X,G)} \omi_{X,G},
  \end{gather}
  et, pour $(f,\al): (X,G)\to (Y,H)$ un morphisme dans $\Eq/\eta$,
  on a
  \begin{gather}
    \label{eq.e.eh}
    \omi_{X,G} [f/\al]^* \simeq (f,\al)^* \omi_{Y,H},\\
    \label{eq.e.eb}
    \omi_{Y,H} R[f/\al]_* \simeq R(f,\al)_* \omi_{X,G},
  \end{gather}
  et, lorsque $f$ est séparé, on a
  \begin{gather}
    \label{eq.e.!h}
    \omi_{X,G} R[f/\al]^! \simeq R(f,\al)^! \omi_{Y,H},\\
    \label{eq.e.!b}
    \omi_{Y,H} R[f/\al]_! \simeq R(f,\al)_! \omi_{X,G}.
  \end{gather}
\end{prop}

\begin{proof}
  Au niveau des topos étales, on a un diagramme $2$-commutatif de foncteurs
  images inverses
  \cite[12.3.3, 12.4.6]{Laumon-MB}
  \[\xymatrix{[X/G]\sptilde\ar[r]^\sim\ar[rd] & (X,G)\sptilde\ar[d]\\ &X\sptilde} \]
  d'où l'existence de $\omi_{X,G}$.
  Les isomorphismes \eqref{eq.e.Hom} à \eqref{eq.e.eb} sont clairs. Les isomorphismes \eqref{eq.e.!h} et \eqref{eq.e.!b} en déduisent par la dualité.
\end{proof}

\begin{remqe}
Le diagramme
\[\xymatrix{\Dbc([Y/H],\Qlb) \ar@{}[rdd]|{(*)} \ar[dd]^{[f/\al]^*}
\ar[r]^{\omi_{Y,H}}\ar[rrd]_(.4){p^*_{Y,H}}
&\Dbc(Y,H,\Qlb) \ar[rd]^{\omega_{Y,H}} \ar@{-->}[dd]^{(f,\al)^*}\\
&&\Dbc(Y,\Qlb) \ar[dd]^{f^*} \\
\Dbc([X/G],\Qlb) \ar@{-->}[r]^{\omi_{X,G}}\ar[rrd]_(.4){p^*_{X,G}}
&\Dbc(X,G,\Qlb)\ar@{-->}[rd]^{\omega_{X,G}} \\
&&\Dbc(X,\Qlb)}
\]
dont le carré $(*)$ correspond à \eqref{eq.e.eh}, est $2$-commutatif. Des diagrammes analogues existent pour \eqref{eq.e.Hom},
\eqref{eq.e.tD} et \eqref{eq.e.!h}.
\end{remqe}

\begin{ssect}\label{ss.omi}
  En vertu de \ref{prop.omi}, on peut interpréter certains cas de
  \ref{prop.cb} et \ref{prop.quot} en termes des champs. Si $(f,\al) : (X,G) \to
  (Y,Q)$ est un morphisme de $\Eq/\eta$ vérifiant les hypothèses
  de \ref{prop.quot} (b) ($\al$ est surjectif et $f$ fait de $X$ un $(\Ker \al)$-torseur sur
  $Y$), alors $[f/\al] :[X/G] \to [Y/Q]$ est un isomorphisme, d'où
  les conclusions de \ref{prop.quot} (a) et (b) (les flèches d'adjonction $\Id \to (f,\al)_* (f,\al)^*$
  et $(f,\al)^*(f,\al)_* \to \Id$ sont des isomorphismes). La
  proposition \ref{prop.cb} (cas $S=\eta$) s'interprète comme un théorème de changement de base pour les champs algébriques (\cf \cite[\S~12]{Laszlo-Olsson2}) en vertu de la proposition suivante.
\end{ssect}

\begin{prop}\label{prop.Mackeychamps}
  Soient $(f,\al) :(X,G) \to (Y,H)$, $(g,\beta) : (Y',H')\to (Y,H)$
  deux morphismes de même but dans $\Eq/\eta$. Pour $r\in H$,
  formons le carré cartésien
  \beq\label{eq.dc.cb1}
  \xymatrix{(X', G')_r \ar[rr]^{(h,\gamma)_r} \ar[d]^{(f',\al')_r}
  && (X, G)\ar[d]^{(f,\al)}\\
  (Y',H')\ar[r]^{(g,\beta)} & (Y,H)\ar[r]^{T_r} &(Y,H)}
  \eeq
  où $T_r = (a_r, h \mapsto r^{-1} h r)$ est comme dans
  \ref{ss.cat}. Le $\eta$-champs $[X'/G']_r = F((X',G')_r)$ ne
  dépend, à isomorphisme près, que de la double classe $(\Img \beta) r (\Img \al) \subset
  H$ et le carré
  \beq\label{eq.dc.cb2}
  \xymatrix{
    \coprod_r[X'/G']_r \ar[d]_{[f'/\al']_r} \ar[r]^{[h/\gamma]_r} & [X/G] \ar[d]^{[f/\al]} \\
    [Y'/H'] \ar[r]^{[g/\beta]} & [Y/H]   }
  \eeq
  est $2$-cartésien, où $r$
  parcourt un système de représentants de $\Img \beta\backslash
  H/\Img \al$.
\end{prop}

\begin{proof}
  Comme $F(T_r) \simeq \Id_{[Y/H]}$, \eqref{eq.dc.cb2} est
  $2$-commutatif. Il suffit de montrer que pour tout $\eta$-schéma
  affine connexe non vide $U$, le foncteur
  \[\coprod_r([X'/G']_r)_U
  \sto{\left(([f'/\al']_r)_U, ([h/\gamma]_r)_U\right)} [Y'/H']_U \times_{[Y/H]_U}[X/G]_U
  \]
  est une équivalence de catégories et que l'image essentielle de $([X'/G']_r)_U$
  ne dépend que de la double classe de~$r$. Construisons un quasi-inverse $v$ du foncteur.
  Un objet de $ [Y'/H']_U
  \times_{[Y/H]_U}[X/G]_U$ est un triplet $(Z,W,k)$ où $Z\in
  [Y'/H']_U$, $W\in [X/G]_U$, $k$ est un isomorphisme $[g/\beta]_U(Z) \simto
  {}[f,\al]_U(W)$ dans $[Y/H]_U$. Rappelons que $Z$ est un
  $H'$-torseur sur $U$ muni d'un $\eta$-morphisme
  $H'$-équivariant $Z\to Y'$ et le morphisme $(Z,H') \to
  ([g/\beta]_U(Z),H)$ est cocartésien (\ref{ss.cat}). Notons $\beta_* Z
  =[g/\beta]_U(Z)$. De même pour $W$ et $[f,\al]_U(W)$ : $\al_* W =
  [f,\al]_U(W)$. L'isomorphisme $k$ est un $Y$-morphisme
  $H$-équivariant. Pour $r\in H$, formons le carré cartésien
  au-dessus de \eqref{eq.dc.cb1}
  \[\xymatrix@C=5em{V_r \ar[d]\ar[rr] && W\ar[d]\\ Z\ar[r] & \beta_* Z\ar[r]^{a_r\circ k = k \circ a_r}_{\sim} & \al_* W} \]
  Notons $C=\ensdr{r\in G}{V_r\neq \vide}$. Alors $C\in \Img \beta\backslash H/\Img
  \al$. Pour $r\in C$, $V_r$ est un $G'$-torseur sur~$U$. On pose
  $v(Z,W,k)=V_r$.
\end{proof}

Soient $E$, $I$, $\gamma$ comme dans \ref{ss.compte}.

\begin{defn}
Soient $\cX \in \Ch/\eta$, $(L_\lbd)_\ldI\in \prod_{\ldI}
\Dbc(\cX,\lbar{\Q_{l_\lbd}})$. On dit que $(L_\lbd)_\ldI$ est
$(E,I,\gamma)$-compatible (ou $E$-compatible s'il n'y a pas de
confusion à craindre) si pour tout morphisme $i : x \to \cX$ dans
$\Ch/\eta$ où $x$ est le spectre d'une extension finie de $K$,
$(i^* L_\lbd)_\ldI$ appartient à $\Dbc(x,E)$ (\ref{def.Ecomp}).
\end{defn}

Les systèmes $E$-compatibles sur $\cX$ forment une sous-\cat
triangulée de $\prod_{\ldI} \Dbc(\cX,\lbar{\Q_{l_\lbd}})$, notée
$\Dbc(\cX,G,E)$. Lorsque $\cX = X$ est un schéma, la définition
ci-dessus coïncide avec \ref{def.Ecomp} : $\Dbc(\cX,E) \simeq
\Dbc(X,E) = \Dbc(X,\{1\},E)$. Plus généralement, on a

\begin{prop}\label{prop.Ecompo}
  Pour $(X,G)\in \Eq/\eta$ et $(L_\lbd)_\ldI\in \prod_{\ldI}
  \Dbc([X/G],\lbar{\Q_{l_\lbd}})$, $(L_\lbd)_\ldI$ appartient à
  $\Dbc([X/G],E)$ \ssi $(\omi_{X,G}L_\lbd)_\ldI$ appartient à
  $\Dbc(X,G,E)$.
\end{prop}

\begin{proof}
  La suffisance est claire. En effet, soient $(L_\lbd)_\ldI \in \prod_\ldI
  \Dbc([X/G],\Qllb)$ tel que $(\omi_{X,G}L_\lbd)_\ldI \in
  \Dbc(X,G,E)$, $i: x\to [X/G]$ dans $\Ch/\eta$ où $x$ est le spectre d'une extension finie de $K$.
  D'après \ref{prop.ch}, il existe des morphismes $s: (Y,H)\to (x,\{1\})$, $t : (Y,H) \to (X,G)$
  dans $\Eq/\eta$ vérifiant $s\in S$, tels que $i=F(t)F(s)^{-1}$.
  D'après \ref{prop.omi},
  \[i^* L_\lbd \simeq F(s)_* F(t)^* L_\lbd \simeq s_*t^*\omi_{X,G}L_\lbd.\]
  Donc $(i^*L_\lbd)_\ldI$ est $E$-compatible, en vertu de
  \ref{rq.th} (iv).

  Prouvons la nécessité. Soient $(L_\lbd)_\ldI \in \Dbc([X/G],E)$,
  $x\in |X|$. Alors, d'après \ref{prop.omi}, $i_x^* \omi_{X,G}L_\lbd \simeq \omi_{x,G_d(x)} F(i_x)^*
  L_\lbd$, où $i_x : (x,G_d(x)) \to (X,G)$. Quitte à
  remplacer $(X,G)$ par $(x,G_d(x))$, $L_\lbd$ par $F(i_x)^*
  L_\lbd$, on peut supposer que $X$ est le spectre d'une extension finie de $K$.
  On applique alors la construction \ref{ss.Cmg} : pour $g\in G$, $m\ge
  1$, on a construit
  \[\xymatrix{(X^{(m,g)},\{1\}) & \ar[l]_-{e_{m,g}}(X_{\eta_{n_g m}}, \Z/n_g\Z) \ar[r]^-{d_{m,g}} & (X,G),} \]
  donnant lieu à
  \[\xymatrix{X^{(m,g)} & \ar[l]_-{Fe_{m,g}}^-\sim [X_{\eta_{n_g m}}/ (\Z/n_g\Z)] \ar[r]^-{Fd_{m,g}} & [X/G].} \]
  D'après \ref{prop.omi},
  \[
  C_{m,g} \omi_{X,G} L_\lbd = e_{m,g *} d_{m,g}^* \omi_{X,G} L_\lbd \simeq (Fe_{m,g})_*
  (Fd_{m,g})^* L_\lbd \simeq ((Fd_{m,g})(Fe_{m,g})^{-1})^* L_\lbd.
  \]
  Donc on a $(C_{m,g} \omi_{X,G} L_\lbd)_\ldI \in \Dbc(X^{(m,g)},E)$.
  Par conséquent, on a $(\omi_{X,G} L_\lbd)_\ldI \in \Dbc(X,G,E)$ en vertu de \ref{prop.DL}.
\end{proof}

\newcommand{\vieuxlemme}{%
Le lemme suivant est une variante de \cite[6.3]{Laumon-MB}.

\begin{lemm} Soient $S$ un schéma \qs, $\cX$ un $S$-champ algébrique
(\qs), $k$ un $S$-corps, $x=\Spec k$, $H$ un groupe fini, $i: B(H/x)
\to \cX$ un morphisme. Alors il existe un diagramme $2$-commutatif
\[\xymatrix{&[X/G]\ar[d]^\phi \\
B(H/x) \ar[ur]\ar[r]^i & \cX}
\]
où $X$ est un schéma affine, $G$
un groupe fini agissant sur $X$ à droite et $\phi$ un morphisme
lisse.
\end{lemm}

\begin{proof} On procède comme dans \ibid[6.7]. On prend un
diagramme $2$-commutatif à carrés $2$-cartésiens
\[\xymatrix{T_2\ar[r]\ar[d] & B_2 \ar[r]\ar[d] & X_2\ar[d]^{\pi}\\
x \ar[r] & B(H/x) \ar[r]^{i} & \cX}
\]
où $X_2$ est un schéma affine, $\pi$ est lisse et $T_2$ est non
vide. Le morphisme $T_2 \to B_2$ est un $G$-torseur et $T_2$ est un
espace algébrique lisse sur $x$. Il existe une immersion fermée
$G$-équivariante $\Spec L \to T_2$, où $L$ est une $k$-algèbre
finie étale non nulle. Posons $d=[L:K]$. On a ainsi une section
$G$-équivariante de $\ET_d(T_2/x) \to x$, ou, ce qui revient au
même, un diagramme $2$-commutatif
\[\xymatrix{&&\ET_d(X_2/\cX)\ar[d] \\
x\ar[rru]^{x_1}\ar[r] & B(H/x)\ar[ru]\ar[r]^i & \cX}
\]
Quitte à remplacer $\cX$ par $\ET_d (X_2/\cX)$, \ops que $\cX$ est
un $S$-champ \DM. Dans ce cas, on peut prendre $\pi$ étale. On a
$\ET_d(X_2/\cX) = [\SEC_d(X_2/\cX)/\fS_d/S]$ et $\SEC_d(X_2/\cX)$
est un schéma quasi-affine. Le point $x_1$ correspond à une
$\fS_d$-trajectoire de $\SEC_d(X_2/\cX)$, qui est donc contenue dans
un ouvert affine $\fS_d$-\eqrt, d'où le résultat en prenant pour
$X$ ledit ouvert et pour $G$ le groupe $\fS_d$.
\end{proof}}

\begin{prop}\label{prop.chlisse}
  Pour $\cX \in \Ch/\eta$,
$(L_\lbd)_\ldI \in \prod_{\ldI}\Dbc(\cX,\lbar{\Q_{l_\lbd}})$ est
$E$-compatible \ssi pour tout morphisme \emph{lisse} $\phi : X\to
\cX$, où $X$ un schéma affine, on a $(\phi^* L_\lbd)_\ldI \in
\Dbc(X,E)$.
\end{prop}

\begin{proof}
  La nécessité est claire. La suffisance découle de
  \cite[6.3]{Laumon-MB}.
\end{proof}

\begin{prop}\label{prop.chcomp}
La $E$-\compte est stable par la dualité, $\otimes$, $R\cHom$,
$f^*$, $Rf^!$, et, lorsque $f$ est relativement \DM, par $Rf_*$ et
$Rf_!$.
\end{prop}

\begin{proof}
La stabilité par $\otimes$ et par $f^*$ résulte de la définition.
Pour la dualité, on applique \ref{prop.chlisse}. Soient
$(L_\lbd)_\ldI \in \Dbc(\cX,E)$,  $X$ un schéma affine, $\phi :
X\to \cX$ un morphisme lisse. \Ops que $\phi$ est purement de
dimension $n$. Alors
\[\phi^* D_{\cX} L_\lbd \simeq D_{X} R\phi^! L_\lbd
\simeq (D_X \phi^* L_\lbd)(-n)[-2n],
\]
donc $(\phi^*D_{\cX}L_\lbd)_\ldI$ est \Ecomp.

Il reste à vérifier la stabilité par $Rf_!$ pour $f: \cX\to \cY$
relativement \DM. D'après le théorème de changement de base \cite[12.5.3]{Laszlo-Olsson2}, \ops
que $\cY=y$ est le spectre d'une extension finie de $K$. Alors $\cX$
est un champs \DM. D'après \cite[6.1.1]{Laumon-MB}, il existe un
ouvert non vide $\cU \simeq [X/G]$, où $X$ est un schéma affine
\tf sur $y$, $G$ un groupe fini agissant à droite sur $X$. Soient
$j: \cU\hra \cX$, $i$~l'immersion fermée complémentaire. Alors on
a le triangle distingué
\[R(fj)_! j^* L_\lbd \to Rf_! L_\lbd \to R(fi)_! i^* L_\lbd \to.\]
D'après \ref{theo.6op}, \ref{prop.ch}, \ref{prop.omi} et \ref{prop.Ecompo}, $R(fj)_!$
préserve la \Ecompte. On conclut par récurrence noethérienne.
\end{proof}

\begin{exem}
  Soit $\pi: V\to \eta$ un fibré vectoriel de rang constant~$r$.
  Alors $\G_m = \G_{m,\eta}$ agit par homothétie sur~$V$ et
  on peut former le $\eta$-champ quotient $\pib : \fV = [V/\G_m] \to
  \eta$. Soit $\pi\spcheck : V\spcheck \to \eta$ le fibré vectoriel
  dual de $\pi$. Comme précédemment, on forme le $\eta$-champ
  quotient $\pib\spcheck : \cV\spcheck = [V\spcheck/\G_m] \to
  \eta$. Alors la transformation de Fourier homogène
  \cite[1.5]{Laumon-hom} $\Four_{\fV/\eta} : \Dbc(\fV,\Qlb) \to \Dbc(\fV \spcheck,
  \Qlb)$ préserve la $E$-compatibilité.

  En effet, comme dans la démonstration de \ibid[1.9], il suffit de
  montrer que pour $(L_\lbd)_\ldI \in \Dbc(\fV, E)$, les systèmes
  \begin{gather*}
    \Four_{\cV/\eta}(i_* i^* L_\lbd) \simeq (\pi \spcheck)^* \sigma^*
    i^* L_\lbd[r],\\
    i_* \spcheck (i \spcheck)^* \Four_{\fV/\eta} (j_! j^* L_\lbd)
    \simeq i_* \spcheck \sigma^* R\pi_! j_! j^* L_\lbd [r],\\
    j_! \spcheck (j \spcheck)^* \Four_{\fV/\eta} (j_! j^* L_\lbd)
    \simeq j_! \spcheck R\pr_! \spcheck (\pr^* j^* L_\lbd \otimes RJ_* \Qllb)
    [r-1]
  \end{gather*}
  sont $E$-compatibles. Ici $i: B(\G_m) \to \fV$, $i \spcheck : B(\G_m) \to \fV
  \spcheck$, $j : \bP(V) \hra \fV$, $j \spcheck : \bP(V \spcheck) \hra
  \fV$, $\pi=[\pi/\G_m] : \fV \to B(\G_m)$, $\pi \spcheck = [\pi \spcheck/\G_m]:
  \fV \spcheck \to B(\G_m)$, $\sigma : B(\G_m) \to B(\G_m)$ envoie
  $\fL$ sur $\fL^{\otimes -1}$, $\pr \spcheck$ et $\pr$ sont les
  deux projections canoniques de $\bP(V \spcheck) \times \bP(V)$, $J$
  est une immersion ouverte \ibid[1.5]. Ces morphismes étant tous
  relativement de Deligne-Mumford, voire représentables
  \cite[3.9]{Laumon-MB}, il suffit d'appliquer \ref{prop.chcomp}.
\end{exem}

Pour les champs algébrique de type fini sur un corps fini, on
dispose, par \cite[\S~9]{Laszlo-Olsson3}, d'un formalisme de poids.
On a des analogues de \ref{ss.j!*} et de \ref{prop.poids}.

\renewcommand{\refname}{Bibliographie}
\bibliographystyle{smfalpha}
\bibliography{rat1e}

\providecommand{\SortNoop}[1]{}
\providecommand{\bysame}{\leavevmode ---\ }
\providecommand{\og}{``}
\providecommand{\fg}{''}
\providecommand{\smfandname}{\&}
\providecommand{\smfedsname}{\'eds.}
\providecommand{\smfedname}{\'ed.}
\providecommand{\smfmastersthesisname}{M\'emoire}
\providecommand{\smfphdthesisname}{Th\`ese}
\begin{thebibliography}{{\SortNoop{zz}}{\SortNoop{}}{\sigle{EGAIV}}67}

\bibitem[BBD82]{BBD}
{\scshape A.~A. Be{\u\i}linson, J.~Bernstein {\normalfont \smfandname}
  P.~Deligne} -- {\og Faisceaux pervers\fg}, in \emph{Analyse et topologie sur
  les espaces singuliers (I)}, Astérisque, vol. 100, Soc.\ math.\ France, 1982,
  p.~5--171.

\bibitem[BLR90]{BLR}
{\scshape S.~Bosch, W.~Lütkebohmert {\normalfont \smfandname} M.~Raynaud} --
  \emph{\eng{Néron models}}, Ergeb.\ Math.\ Grenzgeb.\ (3), vol.~21,
  Springer-Verlag, 1990.

\bibitem[Con97]{Conrad}
{\scshape B.~Conrad} -- {\og \eng{Deligne's notes on Nagata
  compactifications}\fg}, prépublication, 1997.

\bibitem[Del73]{constantes}
{\scshape P.~Deligne} -- {\og Les constantes locales des équations
  fonctionnelles des fonctions {$L$}\fg}, in \emph{Modular functions of one
  variable II}, Lecture Notes in Math., vol. 349, Springer-Verlag, 1973,
  p.~501--597.

\bibitem[Del80]{WeilII}
\bysame , {\og La conjecture de {Weil : II}\fg}, \emph{Publ.\ math.\ IHÉS}
  \textbf{52} (1980), p.~137--252.

\bibitem[dJ96]{deJong1}
{\scshape A.~J. de~Jong} -- {\og \eng{Smoothness, semi-stability and
  alterations}\fg}, \emph{Publ.\ math.\ IHÉS} \textbf{83} (1996), p.~51--93.

\bibitem[dJ97]{deJong2}
{\scshape A.~J. de~Jong} -- {\og Families of curves and alterations\fg},
  \emph{Ann.\ Inst.\ Fourier} \textbf{47} (1997), no.~2, p.~599--621.

\bibitem[DL76]{Deligne-Lusztig}
{\scshape P.~Deligne {\normalfont \smfandname} G.~Lusztig} -- {\og
  Representations of reductive groups over finite fields\fg}, \emph{Ann.\ of
  Math.\ (2)} \textbf{103} (1976), p.~103--161.

\bibitem[Eke90]{Ekedahl}
{\scshape T.~Ekedahl} -- {\og \eng{On the adic formalism}\fg}, in \emph{The
  Grothendieck Festschrift, vol.\ II}, Progr.\ Math., vol.~87, Birkhäuser,
  1990, p.~197--218.

\bibitem[Fuj02]{Fujiwara}
{\scshape K.~Fujiwara} -- {\og \eng{Independence of $l$ for intersection
  cohomology (after Gabber)}\fg}, in \emph{\eng{Algebraic geometry 2000,
  Azumino}}, Adv. Stud. Pure Math., vol.~36, Math.\ Soc.\ Japan, 2002,
  p.~145--151.

\bibitem[Ill06]{IllMisc}
{\scshape L.~Illusie} -- {\og \eng{Miscellany on traces in $l$-adic cohomology:
  a survey}\fg}, \emph{Jpn. J.~Math.\ (3)} \textbf{1} (2006), p.~107--136.

\bibitem[Kat83]{Katz}
{\scshape N.~M. Katz} -- {\og \eng{Wild ramification and some problems of
  ``independence of~$l$''}\fg}, \emph{Amer.\ J. Math.} \textbf{105} (1983),
  p.~201--227.

\bibitem[Lau87]{Laumon}
{\scshape G.~Laumon} -- {\og Transformation de {F}ourier, constantes
  d'équations fonctionnelles et conjecture de {W}eil\fg}, \emph{Publ.\ math.\
  IHÉS} \textbf{65} (1987), p.~131--210.

\bibitem[Lau03]{Laumon-hom}
\bysame , {\og Transformation de {F}ourier homogène\fg}, \emph{Bull.\ Soc.\
  math.\ France} \textbf{131} (2003), p.~527--551.

\bibitem[LMB00]{Laumon-MB}
{\scshape G.~Laumon {\normalfont \smfandname} L.~Moret-Bailly} -- \emph{Champs
  algébriques}, Ergeb.\ Math.\ Grenzgeb.\ (3), vol.~39, Springer, 2000.

\bibitem[LO08]{Laszlo-Olsson2}
{\scshape Y.~Laszlo {\normalfont \smfandname} M.~Olsson} -- {\og \eng{The six
  operations for sheaves on Artin stacks II: adic coefficients}\fg},
  \emph{Publ.\ math.\ IHÉS} \textbf{107} (2008), p.~169--210.

\bibitem[LO09]{Laszlo-Olsson3}
\bysame , {\og \eng{Perverse $t$-structure on Artin stacks}\fg},
  \emph{Math.~Z.} \textbf{261} (2009), p.~737--748.

\bibitem[Mie07]{Mieda}
{\scshape Y.~Mieda} -- {\og \eng{On $l$-independence for the étale cohomology
  of rigid spaces over local fields}\fg}, \emph{Compos.\ Math.} \textbf{143}
  (2007), p.~393--422.

\bibitem[Mor08]{Morel-pondere}
{\scshape S.~Morel} -- {\og Complexes pondérés sur les compactifications de
  {B}aily-{B}orel : le cas des variétés de {S}iegel\fg}, \emph{J.~Amer.\ Math.\
  Soc.} \textbf{21} (2008), p.~63--100.

\bibitem[Och99]{Ochiai}
{\scshape T.~Ochiai} -- {\og \eng{$l$-independence of the trace of
  monodromy}\fg}, \emph{Math.\ Ann.} \textbf{315} (1999), p.~321--340.

\bibitem[Sai03]{Saito-weight}
{\scshape T.~Saito} -- {\og Weight spectral sequences and independence of
  {$l$}\fg}, \emph{J.~Inst.\ math.\ Jussieu} \textbf{2} (2003), p.~583--634.

\bibitem[Ser98]{Serre}
{\scshape J.-P. Serre} -- \emph{Représentations linéaires des groupes finis},
  5\ieme \smfedname, Hermann, 1998.

\bibitem[Vid04]{Vidal}
{\scshape I.~Vidal} -- {\og Théorie de {B}rauer et conducteur de {S}wan\fg},
  \emph{J.~Algebraic Geom.} \textbf{13} (2004), p.~349--391.

\bibitem[Zhe08]{Zheng-inte}
{\scshape W.~Zheng} -- {\og Sur la cohomologie des faisceaux {$l$}-adiques
  entiers sur les corps locaux\fg}, \emph{Bull.\ Soc.\ math.\ France}
  \textbf{136} (2008), p.~465--503.

\bibitem[{\SortNoop{zz}}{\SortNoop{}}{\sigle{EGAIV}}67]{EGAIV}
{\scshape A.~{\SortNoop{zz}}{\SortNoop{}}{\sigle{EGAIV}}Grothendieck} -- {\og
  Éléments de géométrie algébrique : {IV. É}tude locale des schémas et des
  morphismes de schémas\fg}, \emph{Publ.\ math.\ IHÉS} \textbf{20, 24, 28, 32}
  (1964--1967).

\bibitem[{\SortNoop{zz}}{\SortNoop{}}{\sigle{SGA1}}03]{SGA1}
\emph{Revêtements étales et groupe fondamental} -- Séminaire de géométrie
  algébrique du Bois-Marie 1960--1961, dirigé par A.~Grothendieck, Doc. math.,
  vol.~3, Soc.\ math.\ France, 2003.

\bibitem[{\SortNoop{zz}}{\SortNoop{}}{\sigle{SGA4}}73]{SGA4}
\emph{Théorie des topos et cohomologie étale des schémas} -- Séminaire de
  géométrie algébrique du Bois-Marie 1963--1964, dirigé par M.~Artin,
  A.~Grothendieck, \mbox{J.-L.}~Verdier, Lecture Notes in Math., vol. 269, 270,
  305, Springer-Verlag, 1972--1973.

\bibitem[{\SortNoop{zz}}{\SortNoop{}}{\sigle{SGA4\textonehalf}}77]{SGA4d}
{\scshape P.~{\SortNoop{zz}}{\SortNoop{}}{\sigle{SGA4\textonehalf}}Deligne} --
  \emph{Cohomologie étale}, Lecture Notes in Math., vol. 569, Springer-Verlag,
  1977.

\bibitem[{\SortNoop{zz}}{\SortNoop{}}{\sigle{SGA7}}73]{SGA7}
\emph{Groupes de monodromie en géométrie algébrique} -- Séminaire de géométrie
  algébrique du Bois-Marie 1967--1969, I, dirigé par A.~Grothendieck, II, par
  P.~Deligne, N.~Katz, Lecture Notes in Math., vol. 288, 340, Springer-Verlag,
  1972--1973.

\end{thebibliography}

\vspace{2ex}

{\footnotesize Weizhe Zheng

Columbia University MC4406, 2990 Broadway, New York, NY 10027, USA

\texttt{zheng@math.columbia.edu}}
\end{document}